\documentclass[3p,times]{elsarticle}

\usepackage{ecrc}

\volume{}

\firstpage{1}

\journalname{Preprint}
\jid{Inerst Preprint}

\runauth{Zeifang and Sch\"utz}

\usepackage{amssymb}
\usepackage{bibref}

\usepackage{amsmath,amssymb,graphicx,stmaryrd}
\usepackage{amsthm}
\usepackage[dvipsnames]{xcolor}
\usepackage{tikz}
\usepackage{hyperref}
\hypersetup{
	colorlinks,
	linkcolor={red!50!black},
	citecolor={blue!50!black},
	urlcolor={blue!80!black}
}
\usepackage{pgfplots}
\usepackage{pgfplotstable}
\usepackage{placeins}
\usepackage{color}
\usepackage{pdfpages}
\usetikzlibrary{shapes,arrows,calc,external,matrix,positioning,patterns}
\usepackage{rotating}
\usepackage{ifthen}

\usepackage[mathlines]{lineno} 
\usepackage{etoolbox} 
\newcommand*\linenomathpatch[1]{%
	\cspreto{#1}{\linenomath}%
	\cspreto{#1*}{\linenomath}%
	\csappto{end#1}{\endlinenomath}%
	\csappto{end#1*}{\endlinenomath}%
}
\newcommand*\linenomathpatchAMS[1]{%
	\cspreto{#1}{\linenomathAMS}%
	\cspreto{#1*}{\linenomathAMS}%
	\csappto{end#1}{\endlinenomath}%
	\csappto{end#1*}{\endlinenomath}%
}
\expandafter\ifx\linenomath\linenomathWithnumbers
\let\linenomathAMS\linenomathWithnumbers
\patchcmd\linenomathAMS{\advance\postdisplaypenalty\linenopenalty}{}{}{}
\else
\let\linenomathAMS\linenomathNonumbers
\fi

\linenomathpatch{equation}
\linenomathpatchAMS{gather}
\linenomathpatchAMS{multline}
\linenomathpatchAMS{align}
\linenomathpatchAMS{alignat}
\linenomathpatchAMS{flalign}

\renewcommand{\vec}[1]{\mathbf{#1}}
\DeclareMathOperator{\R}{\mathbb{R}}

\newcommand{\II}{\mathcal I}

\newcommand{\eps}{\varepsilon}
\newcommand{\epsFD}{\varepsilon_{\text{FD}}}

\newcommand{\dt}{\Delta t}
\newcommand{\Rt}{\vec{R}^{(1)}}
\newcommand{\Rtt}{\vec{R}^{(2)}}
\newcommand{\ww}{\mathbf{W}}

\DeclareMathOperator{\Id}{Id}

\DeclareMathOperator{\diag}{diag}

\newcommand{\Tend}{T_{end}}

\newcommand{\km}{k_{\max}}

\DeclareMathOperator{\Sr}{Sr}

\newtheorem{algorithm}{Algorithm}

\newtheorem{remark}{Remark}

\newcommand{\algs}[4]{#1^{#2,[#3],{#4}}}

\newcommand{\method}[2]{\text{HBPC}(#1,#2)}

\DeclareMathAlphabet\mathbfcal{OMS}{cmsy}{b}{n}
\newcommand{\Jgeo}{J_\text{\tiny geo}}

\pgfplotstableset{col sep=comma}
\pgfplotscreateplotcyclelist{epslist}{%
	{red,thick,mark=o},
	{orange,thick,mark=pentagon},
	{green!60!black!90!,thick,mark=square},
	{blue!90!white!50,thick,mark=triangle},
	{blue!60!black!80,thick,mark=x},
	{violet	,thick	,mark=diamond},	
	{red,thick,dashed,mark=o},
	{orange,thick,dashed,mark=pentagon},
	{green!60!black!90!,thick,dashed,mark=square},
	{blue!90!white!50,thick,dashed,mark=triangle},
	{blue!60!black!80,thick,dashed,mark=x},
	{violet	,thick	,dashed,mark=diamond},
	{red,thick,dotted,mark=o},
	{orange,thick,dotted,mark=pentagon},
	{green!60!black!90!,thick,dotted,mark=square},
	{blue!90!white!50,thick,dotted,mark=triangle},
	{blue!60!black!80,thick,dotted,mark=x},
	{violet	,thick	,dotted,mark=diamond}}
\pgfplotscreateplotcyclelist{green}{%
	{green!60!black!10!,thick,mark=o},
	{green!60!black!20!,thick,mark=pentagon},
	{green!60!black!30!,thick,mark=square},
	{green!60!black!40!,thick,mark=triangle},
	{green!60!black!50!,thick,mark=x},
	{green!60!black!60!,thick,mark=otimes},
	{green!60!black!70!,thick,mark=o},
	{green!60!black!80!,thick,mark=pentagon},
	{green!60!black!90!,thick,mark=square},
	{green!60!black!100!,thick,mark=triangle},
	{blue!60!black!10!,thick,mark=o},
	{blue!60!black!20!,thick,mark=pentagon},
	{blue!60!black!30!,thick,mark=square},
	{blue!60!black!40!,thick,mark=triangle},
	{blue!60!black!50!,thick,mark=x},
	{blue!60!black!60!,thick,mark=otimes},
	{blue!60!black!70!,thick,mark=o},
	{blue!60!black!80!,thick,mark=pentagon},
	{blue!60!black!90!,thick,mark=square},
	{blue!60!black!100!,thick,mark=triangle}}
	
\definecolor{color1}{RGB}{         0  189.2578  255.0000}
\definecolor{color2}{RGB}{         0  117.5391  255.0000}
\definecolor{color3}{RGB}{         0   45.8203  255.0000}
\definecolor{color4}{RGB}{   25.8984         0  255.0000}
\definecolor{color5}{RGB}{   97.6172         0  255.0000}
\definecolor{color6}{RGB}{  169.3359         0  255.0000}
\definecolor{color7}{RGB}{  241.0547         0  255.0000}
\definecolor{color8}{RGB}{  255.0000         0  197.2266}
\definecolor{color9}{RGB}{  255.0000         0  125.5078}
\definecolor{color10}{RGB}{  255.0000         0   53.7891}

\pgfplotscreateplotcyclelist{hierarchy}{%
	{color1,thick,mark=o},
	{color2,thick,mark=pentagon},
	{color3,thick,mark=square},
	{color4,thick,mark=triangle},
	{color5,thick,mark=x},
	{color6,thick,mark=otimes},
	{color7,thick,mark=o},
	{color8,thick,mark=pentagon},
	{color9,thick,mark=square},
	{color10,thick,mark=triangle}}

\definecolor{col1}{RGB}{         0  135.4688  255.0000}
\definecolor{col2}{RGB}{         0    9.9609  255.0000}
\definecolor{col3}{RGB}{  115.5469         0  255.0000}
\definecolor{col4}{RGB}{  241.0547         0  255.0000}
\definecolor{col5}{RGB}{  255.0000         0  143.4375}
\definecolor{col6}{RGB}{  255.0000         0   17.9297}

\pgfplotscreateplotcyclelist{hierarchy6}{%
	{col1,thick,mark=o},
	{col2,thick,mark=pentagon},
	{col3,thick,mark=square},
	{col4,thick,mark=triangle},
	{col5,thick,mark=x},
	{col6,thick,mark=otimes}}

\newboolean{compilefromscratch}
\setboolean{compilefromscratch}{true}
\setboolean{compilefromscratch}{false}

\begin{document}

	\begin{frontmatter}
		
		\title{Two-derivative deferred correction time discretization for the discontinuous Galerkin method}%
		
		\author[1]{Jonas Zeifang}
		\ead{jonas.zeifang@uhasselt.be}
		\author[1]{Jochen Sch\"utz}
		\ead{jochen.schuetz@uhasselt.be}

		\address[1]{Faculty of Sciences \& Data Science Institute, Hasselt University, Agoralaan Gebouw D, BE-3590 Diepenbeek, Belgium}

		\begin{abstract}
			In this paper, we use an implicit two-derivative deferred correction time discretization approach and combine it with a spatial discretization of the discontinuous Galerkin spectral element method to solve (non-)linear PDEs. The resulting numerical method is high order accurate in space and time. As the novel scheme handles two time derivatives, the spatial operator for both derivatives has to be defined. This results in an extended system matrix of the scheme. We analyze this matrix regarding possible simplifications and an efficient way to solve the arising (non-)linear system of equations.
			It is shown how a carefully designed preconditioner and a matrix-free approach allow for an efficient implementation and application of the novel scheme. For both, linear advection and the compressible Euler equations, up to eighth order of accuracy in time is shown.
			Finally, it is illustrated how the method can be used to approximate solutions to the compressible Navier-Stokes equations.
		\end{abstract}
		
		
	\end{frontmatter}
	
	
	\section{Introduction}\label{sec:Introduction}
	We aim for solving hyperbolic PDEs that can be cast into flux formulation
\begin{align}\label{eq:hyperbolic}
	\vec{w}_t+\nabla_x\cdot\vec{F}(\vec{w})=0,
\end{align}
with state vector $\vec{w}$ and flux $\vec{F}(\vec{w})$. 
Obviously, upon defining 
\begin{align}\label{eq:first_derivative}
	\Rt(\vec{w}):=-\nabla_x\cdot\vec{F}(\vec{w}),
\end{align}
this can be cast as an ODE in some infinite-dimensional function space, 
\begin{align}\label{eq:ODEtypde}
 \vec{w}_t = \Rt(\vec{w}). 
\end{align}
%
While standard time discretization methods (for example Runge-Kutta or BDF methods) only use the information of the first time derivative $\vec{w}_t$ for time stepping, two-derivative methods make use of the second time derivative $\vec w_{tt}$, computed through
\begin{align}\label{eq:second_derivative}
	\vec{w}_{tt}=\left(-\nabla_x\cdot\vec{F}(\vec{w})\right)_t = \nabla_x\cdot\left(-\vec{F}(\vec{w})\right)_t=\nabla_x\cdot\left(-\frac{\partial \vec{F}}{\partial\vec{w}}\vec{w}_t\right)=\nabla_x\cdot\left(-\frac{\partial\vec{F}}{\partial\vec{w}}\Rt(\vec{w})\right)=:\Rtt(\vec{w},\Rt(\vec{w})).
\end{align}
Using this additional derivative introduces more flexibility when designing time stepping methods. This allows to achieve high order of accuracy while using less temporal quadrature nodes than conventional schemes.
The class of two-derivative methods belongs to the general class of multistep-multistage-multiderivative methods~\cite{HaWa73}. Important multiderivative methods are the Taylor methods which are typically referred to the class of Lax-Wendroff methods when considering PDEs~\cite{LaxWend1960}; and multiderivative Runge-Kutta schemes~\cite{KastlungerWanner1972}.
PDE discretizations with Lax-Wendroff methods have been widely addressed in literature, see e.g. the by far not exhaustive list~\cite{Christlieb2015,JiangShuZhang13,Moe2017,Qiu2005,ZorioEtAl}. Also the arbitrary high order derivative Riemann problem (ADER) approach~\cite{Titarev2002} is based on this idea, see~\cite{Busto2020} for a recent overview on this technique.
Explicit two-derivative Runge-Kutta methods have also been used in the context of PDEs, see e.g.~\cite{Seal13} where WENO and DG schemes have been combined with a multiderivative Runge-Kutta methods. Other examples for different PDE discretizations with explicit multiderivative Runge-Kutta methods can e.g. be found in~\cite{Ji2018,He2020,PanXuLiLi2016,Chouchoulis2021}.
A semi-implicit two-derivative Runge-Kutta method has been used as a generalization of a Lax-Wendroff finite difference approach in~\cite{TCW14}. 

Implicit two-derivative methods have only rarely been used in the context of PDE discretizations.
In~\cite{MultiDerHDG2015} an implicit Lax-Wendroff approach for a two-derivative two-point method has been introduced. For the spatial discretization a hybridized discontinuous Galerkin (HDG) method has been used. The authors showed that even if the underlying two-derivative implicit method is A-stable, a severe timestep restriction for the PDE discretization can be observed.
In a follow-up publication, this timestep restriction is overcome by the introduction of a systematic approach to introduce an additional solution variable per time derivative~\cite{SSJ2017}. In that work, the authors show that the discretization of a linear PDE preserves the stability properties of the ODE solver, corresponding to the method of lines approach. The authors illustrate this with a DG method with a two- and three-derivative two-point time discretization. More details can also be found in~\cite{DissAlex}. This idea will lay the foundation of the current work.

In this work, differently to~\cite{SSJ2017}, higher order time discretization is achieved by using a deferred correction time discretization method.
In~\cite{SealSchuetz19}, a fourth order accurate two-derivative deferred correction time discretization method has been presented. Recently, it has been combined with the idea of using multiple stages to obtain higher orders of accuracy~\cite{Schutz2021}. There it has been shown that this class of schemes can be used to solve stiff ODEs while obtaining - at least in principle - arbitrary orders of accuracy. 
Here, these ODE solvers will be used as time discretization method to solve PDEs, while the spatial discretization is done with the discontinuous Galerkin spectral element method (DGSEM)~\cite{kopriva2009}. The novel scheme is constructed such that it can handle non-linear PDEs, which will be illustrated with the compressible Euler and Navier-Stokes equations. One of the main drawbacks of the approach outlined in~\cite{SSJ2017} is that one obtains an extended and hence larger system matrix. In this work, we show how the linearized system arising from this approach can be solved efficiently. Moreover, our novel approach allows for a straight-forward parallelization of the spatial domain as we introduce a matrix-free discretization and a relatively simple preconditioning strategy.
Summing up, the main advances with respect to previous works are
\begin{itemize}
	\item the use of implicit two-derivative high-order deferred correction methods for a PDE discretization, especially with the DGSEM;
	\item the handling of non-linear hyperbolic-parabolic PDEs with the ansatz described in~\cite{SSJ2017};
	\item the introduction of a carefully designed preconditioner which is combined with a matrix-free discretization. This eases the implementation of two-derivative schemes in already existing numerical codes and allows for a straight-forward parallelization of the spatial domain.
\end{itemize}
The work is structured as follows:
In Sec.~\ref{sec:HBPC_Method} the two-derivative deferred correction method is introduced in a semi-discrete setting. In the following section (Sec.~\ref{sec:MD-DGSEM}), the fully discrete scheme is derived. After deriving the operators for the first and second temporal derivative $\Rt$ and $\Rtt$, two approaches how to improve the solution procedure of the arising extended linear system are presented. Both approaches are compared in Sec.~\ref{sec:NumRes} for the linear scalar advection and the compressible Euler equations. Furthermore it is shown how the implementation can be pursued matrix-free and the high order temporal accuracy of the novel scheme is illustrated. The extension of the method to the Navier-Stokes equations and some illustrative applications are presented in Sec.~\ref{sec:Applications}. Finally, a conclusion is drawn and an outlook is given in Sec.~\ref{sec:Conclusion}.

	\section{Semi-Discrete Formulation}\label{sec:HBPC_Method}
	In this section, the semi-discrete in time formulation of the novel scheme is reviewed. For that purpose, we briefly recall the serial algorithm from~\cite{SealSchuetz19} and~\cite{Schutz2021} and slightly modify it. The algorithm describes a \textbf{p}redictor-\textbf{c}orrector approach with $\km$ correction steps to approximate a two-derivative \textbf{H}ermite-\textbf{B}irkhoff Runge-Kutta method of order $q$ and is therefore labeled as $\method{q}{\km}$.
As we do not use an IMEX splitting in this paper as it is done in~\cite{SealSchuetz19,Schutz2021}, the predictor is modified such that a fourth-order two-point Hermite-Birkhoff Runge-Kutta method is successively used to obtain the predicted solution. Please note that for convenience, we stick with the name $\method{q}{\km}$ although we use a different predictor than in the original publication~\cite{Schutz2021}.

\subsection{Hermite-Birkhoff Predictor-Corrector Time Discretization}
The spatial operators to calculate $\vec{w}_t$ and $\vec{w}_{tt}$, viz. $\Rt$ and $\Rtt$ are used to define the semi-discrete in time formulation. In principle, those operators can be discretized by any spatial discretization such as e.g. a discontinuous Galerkin, a finite volume or a finite difference discretization. In this work, we choose the DGSEM, for which the discrete formulation of $\Rt(\vec{w})$ and $\Rtt(\vec{w},\Rt(\vec{w}))$ will be described in Sec.~\ref{sec:MD-DGSEM}. 

\begin{remark}
 $\Rtt$ depends on two arguments, namely $\vec{w}$ and $\Rt(\vec{w})$. We have explicitly chosen this notation, as it will become important in following sections, where we define an auxiliary variable $\boldsymbol{\sigma}$ to equal the discrete form of $\Rt(\vec{w})$. To keep the notation short in the following, however, we will abuse notation in this section only and set $\Rtt(\vec{w}) := \Rtt(\vec{w},\Rt(\vec{w}))$ for the sake of a better readability. 
\end{remark}

The method to be presented below relies on two-derivative Butcher tableaux consisting of matrices $B^{(1)}$, $B^{(2)} \in \R^{s \times s}$ and vector $c\in \R^s$. They define the limiting Hermite-Birkhoff Runge-Kutta scheme and are given in the appendix, Eq.~\eqref{eq:butcher4}-\eqref{eq:butcher8}. More details can be found in~\cite{Schutz2021}. 
These Butcher tableaux define a quadrature formula $\II_l$ of order $q$ through
\begin{align*} 
	\II_l:=
	\dt \sum_{j=1}^s B^{(1)}_{lj} \Rt(\algs{\vec{w}}{n}{k}{j})  + \dt^2 \sum_{j=1}^s B^{(2)}_{lj} \Rtt(\algs{\vec{w}}{n}{k}{j})
\end{align*}
for every stage $1\leq l\leq s$.
Then, the temporal discretization of a PDE of type~\eqref{eq:hyperbolic} with the $\method{q}{\km}$ method and the spatial operators $\Rt$ and $\Rtt$ is given by
\begin{algorithm}[$\method{q}{\km}$]\label{alg:mdpde}
	Solve the following expression for $\algs {\vec{w}} {n} {0} {l}$ and $1 \leq l \leq s$: 
	First, the initial conditions are filled with $\algs{\vec{w}}{-1}{k}{s}:=w_0$.
	To advance the solution to Eq.~\eqref{eq:hyperbolic} from time level $t^{n}$ to time level $t^{n+1}$, fill the values $\algs {\vec{w}} {n} {0} {l}$ using an implicit fourth order predictor
	\begin{enumerate}
		\item \textbf{Predict.} Solve the following expression for $\algs {\vec{w}} {n} {0} {l}$ and $2 \leq l \leq s$: 
		\begin{align}\label{eq:predictor}
		\begin{split}
		\algs {\vec{w}} {n} {0} {1} &:= \algs{\vec{w}}{n-1}{\km}{s}\\
		\algs {\vec{w}} {n} {0} {l} &:= \algs{\vec{w}}{n}{0}{l-1} + \frac{\Delta c_l\dt}{2}\left(\Rt(\algs{\vec{w}}{n}{0}{l-1}) + \Rt(\algs {\vec{w}} n 0 l) \right) + \frac{\left(\Delta c_l\dt\right)^2}{12} \left(\Rtt(\algs {\vec{w}} n 0 {l-1}) - \Rtt(\algs {\vec{w}} n 0 l)\right),
		\end{split}
		\end{align}
		with $\Delta c_l:=c_l-c_{l-1}$.\\
		Subsequently:
		\item \textbf{Correct.} Solve the following for $\algs {\vec{w}} {n} {k+1} {l}$, for each $2 \leq l \leq s$ and each $0 \leq k < k_{\max}$: 
		\begin{align}
		\begin{split}
		\label{eq:corrector}
		\algs {\vec{w}} {n} {k+1} 1 := \algs{\vec{w}}{n-1}{\km}{s}&, \\
		\algs {\vec{w}} {n} {k+1} l := \algs{\vec{w}}{n-1}{\km}{s}& 
		+ {\dt} \left( \Rt(\algs {\vec{w}} {n} {k+1} l) - \Rt(\algs {\vec{w}} n {k} l)\right)
		- \frac{\dt^2}{2}\left(\Rtt(\algs {\vec{w}} {n} {k+1} l)- \Rtt(\algs {\vec{w}} {n} {k} l)\right) + \II_l.
		\end{split}
		\end{align}
		\item \textbf{Update.} In order to preserve the first-same-as-last property, we put
		$\vec{w}^{n+1} := \algs {\vec{w}} n {k_{\max}} s$.
	\end{enumerate}
\end{algorithm}

\begin{remark}
 In \cite{Schutz2021}, it has been shown that the order of accuracy of the class of schemes is given by $\min(4+\km,q)$. That means that starting with the fourth order of the predictor, the schemes pick up one order of accuracy per correction step. This continues until the order of the used quadrature rule is reached.
\end{remark}

\subsection{Solving for the Stage Values}
To solve for the stage values of the predictor and corrector, see Eq.~\eqref{eq:predictor} and Eq.~\eqref{eq:corrector}, one can apply Newton's method. For both predictor and corrector, the resulting non-linear equations are very similar, they can be written in the generalized form
\begin{align}\label{eq:Newton_Residual}
\begin{split}
	\vec{G}(\ww)&:=\vec{g}(\ww)-\vec{b}=0,\\
	\text{with}\quad\vec{g}(\ww)&\hphantom{:}=\ww-\alpha_1\dt\Rt(\ww) + \frac{\alpha_2\dt^2}{2}  \Rtt(\ww).
\end{split}
\end{align}
For the predictor, $\ww$, $\vec{b}$, $\alpha_1$ and $\alpha_2$ are given by
\begin{align}\label{eq:nonlinear_predictor}
\begin{split}
	\ww=\algs{\vec{w}}{n}{0}{l},\quad \vec{b}&=\algs{\vec{w}}{n}{0}{l-1}+\alpha_1\dt \Rt(\algs{\vec{w}}{n}{0}{l-1}) + \frac{\alpha_2\dt^2}{2}  \Rtt(\algs{\vec{w}}{n}{0}{l-1}),\quad \alpha_1=\frac{\Delta c_l}{2},\quad\text{and}\quad\alpha_2=\frac{\Delta c_l^2}{6},
\end{split}
\end{align}
and for the corrector steps the quantities are given by
\begin{align}\label{eq:nonlinear_corrector}
\begin{split}
	\ww=\algs{\vec{w}}{n}{k+1}{l},\quad \vec{b}=\algs{\vec{w}}{n-1}{\km}{s}-\alpha_1\dt \Rt(\algs {\vec{w}} n k l) + \frac{\alpha_2\dt^2}{2} \Rtt(\algs {\vec{w}} n k l),\quad\text{and}\quad\alpha_1=\alpha_2=1.
\end{split}
\end{align}
Consequently, the method consists of subsequent solves of equations of type~\eqref{eq:Newton_Residual}. Hence, two important building blocks remain to be defined to obtain the fully discrete scheme: A discretization for $\Rt$ and $\Rtt$ has to be set up and an efficient method to solve equations of the type $\vec{G}(\ww)=0$ has to be found. Both building blocks will be described in the following section.

	\section{Fully Discrete Scheme}\label{sec:MD-DGSEM}
	\subsection{The Discontinuous Galerkin Spectral Element Method}
The discontinuous Galerkin spectral element method has been introduced in~\cite{kopriva2009}. Based on the discretization of the domain $\Omega$ with $n_E$ quadrangular (2d) or hexahedral elements (3d) $\Omega_e$, it utilizes high order nodal polynomials to represent the solution inside each element. As it is characteristic for discontinuous Galerkin methods, see e.g.~\cite{ReHi73,ShuReviewDG}, discontinuities are allowed across cell boundaries.

\subsubsection{Calculation of the First Temporal Derivative with the DGSEM}
In this subsection, the discrete formulation of the DGSEM is briefly recalled. We closely follow Hindenlang et al.~\cite{MunzDG12}, and refer to the corresponding equations in their work where appropriate.
The DGSEM is derived through the weak formulation of Eq.~\eqref{eq:hyperbolic}
\begin{align}\label{eq:weakform}
	\sum_{e=1}^{n_E}\quad\left(\vec{w}_t,\phi\right)_{\Omega_e}-\left(\vec{F}(\vec{w}),\nabla_x\phi\right)_{\Omega_e}+\left\langle\vec{F}^*(\vec{w}^L,\vec{w}^R)\cdot\vec{n},\phi\right\rangle_{\partial\Omega_e}=0,\qquad\forall\phi\in\Pi_N,
\end{align}
with the test functions $\phi$ taken from the - for now - not nearer specified set of polynomials $\Pi_N$. Element-wise integration over the domain $\Omega$ is denoted by the scalar product $(\cdot,\cdot)$, and $\langle\cdot,\cdot\rangle$ denotes the integral over the cell edges $\partial\Omega_e$. On the cell edges with normals $\vec{n}$, the flux is substituted by a numerical flux $\vec{F}^*$, which depends on ''left'' values $\vec{w}^L$ and ''right'' values $\vec{w}^R$ stemming from two neighboring elements.
This equation is transformed into reference space $\boldsymbol{\xi}=(\xi_1,\xi_2)$ and the cells are mapped onto the reference unit element. For the ease of presentation, we restrict ourselves to two spatial dimensions.
The steps of this transformation can be found in~\cite[Eqs.~(2)-(16)]{MunzDG12}. After this transformation, one obtains the transformed fluxes $\mathbfcal{F}^m$ in direction $m \in \{1,2\}$ and the solution $\vec{w}$ in reference space.
One important building block of the DGSEM is to represent both with the tensor-product of one dimensional Lagrange polynomials $\ell$, each of degree $N$
\begin{align}\label{eq:Interpolation}
	\vec{w}(\boldsymbol{\xi})\approx\vec{w}_h(\boldsymbol{\xi}):=\sum_{i,j=0}^{N}\vec{\hat{w}}_{ij}\ell_i(\xi_1)\ell_j(\xi_2),\quad\text{and}\quad
	\mathbfcal{F}^m(\boldsymbol{\xi})\approx\mathbfcal{F}_h^m(\boldsymbol{\xi}):=\sum_{i,j=0}^{N}\mathbfcal{\hat{F}}^m(\vec{\hat{w}}_{ij})\ell_i(\xi_1)\ell_j(\xi_2),\quad m=1,2,
\end{align}
where $(\hat{\bullet})_{ij}$ denotes the $(i,j)-$th basis coefficient of the polynomial representation. 
As interpolation points, the ${N}+1$ nodes of the Gauss-Legendre quadrature are chosen.
Following, collocation is performed, i.e. the quadrature rule to approximate the integrals in Eq.~\eqref{eq:weakform} uses the same nodes as they are used for the polynomial basis of the solution (Eq.~\eqref{eq:Interpolation}). Additionally, the test functions $\phi$ are chosen to be the same as the ansatz functions, i.e. $\Pi_N$ is the set of tensor-products of the one dimensional Lagrange polynomials of degree $N$.

After some algebraic manipulations, see~\cite[Eqs.~(23)-(37)]{MunzDG12}, the spatial DGSEM operator $\Rt_h(\vec{w})$ (which is an approximation of $\vec{w}_t$) in two dimensions for one element is given by
\begin{align}\label{eq:DGSEM_wt}
	\begin{split}
		\left(\vec{\hat{w}}_{ij}\right)_t=\left(\Rt_h(\vec{w}_h)\right)_{ij}:=-\frac{1}{{\Jgeo}_{ij}}\bigg[
		&\sum_{\lambda=0}^{{N}}\hat{\mathbfcal{F}}^1_{\lambda j}\hat{D}_{i\lambda}
		+\left([\vec f^*\hat{s}]_{j}^{+{\xi_1}}\hat{\ell}_i(1)+[\vec f^*\hat{s}]_{j}^{-{\xi_1}}\hat{\ell}_i(-1)\right)\\
		+&\sum_{\mu=0}^{{N}}\hat{\mathbfcal{F}}^2_{i\mu}\hat{D}_{j\mu}
		+\left([\vec f^*\hat{s}]_{i}^{+{\xi_2}}\hat{\ell}_j(1)+[\vec f^*\hat{s}]_{i}^{-{\xi_2}}\hat{\ell}_j(-1)\right)
		\bigg]\quad\forall~i,j.
	\end{split}
\end{align}
with the abbreviations
\begin{align*}
	\hat{\ell}_i:=\frac{\ell_i}{\omega_i},\quad\text{and}\quad\hat{D}_{ij}:=\left.-\frac{\omega_i}{\omega_j}\frac{\partial\ell_i(\xi)}{\partial\xi}\right|_{\xi=\xi_j},
\end{align*}
and the Jacobian of the geometrical transformation $\Jgeo$, see~\cite[Eqs.~(5), (32) and (38)]{MunzDG12}. Note that $\omega_{i}$ denotes the quadrature weight of the Gauss-Legendre quadrature at position $i$.
At the four edges of the element, which are denoted by $-\xi_1$, $+\xi_1$, $-\xi_2$ and $+\xi_2$, the transformed numerical flux $[\vec{f}^*\hat{s}]$ is evaluated, see~\cite[Eqs.~(31)-(33)]{MunzDG12} for the definition of the transformation and the definition of the surface element $\hat{s}$.
The numerical flux depends on the left and right values at the edge and the normal vector of the edge, i.e. $\vec{f}^*=\vec{f}^*(\vec{w}^L,\vec{w}^R,\vec{n})$.
Here, we use a global Lax-Friedrichs flux
\begin{align}\label{eq:Riemann}
	\vec{f^*}(\vec{w}^L,\vec{w}^R,\vec{n})=\left(\frac{1}{2}\left(\vec{F}(\vec{w}^L)+\vec{F}(\vec{w}^R)\right)+\lambda\left(\vec{w}^L-\vec{w}^R\right)\right)\cdot\vec{n},
\end{align}
where $\lambda$ is a globally constant value. The left and right values $\vec{w}^L$, $\vec{w}^R$ are obtained by evaluation of the solution polynomial at the cell edges, see~\cite[Eq.~(40)]{MunzDG12}.
For more details on the derivation of the DGSEM, see~\cite{kopriva2009,MunzDG12}. The novel scheme described in the current paper has been implemented in the open source code FLEXI\footnote{www.flexi-project.org, GNU GPL v3.0}; for a recent overview on this code see~\cite{Krais2020}.

\subsubsection{Calculation of the Second Temporal Derivative with the DGSEM}
For the considered two derivative time discretization methods, an approximation of $\vec{w}_{tt}$ is additionally required, see Eq.~\eqref{eq:second_derivative}.
In~\cite{SSJ2017}, the artificial quantity
\begin{align}\label{eq:Sigma}
	\boldsymbol{\sigma}:=\Rt_h(\vec{w})
\end{align}
has been introduced in order to facilitate the calculation of the second derivative term.
Here, we follow this idea and start from the discretized equation Eq.~\eqref{eq:weakform} that we differentiate with respect to time $t$ (see also Eq.~\eqref{eq:second_derivative} for comparison)
\begin{align}\label{eq:Euler_wtt}
	\sum_{e=1}^{n_E}\quad\left(\vec{w}_{tt},\phi\right)_{\Omega_e}-\left(\frac{\partial\vec{F(\vec{w})}}{\partial\vec{w}}\cdot\boldsymbol{\sigma},\nabla_x\phi\right)_{\Omega_e}+\left\langle\frac{\partial\vec{F}^*(\vec{w}^L,\vec{w}^R)}{\partial\vec{w}^L}\boldsymbol{\sigma}^L\cdot\vec{n}+\frac{\partial\vec{F}^*(\vec{w}^L,\vec{w}^R)}{\partial\vec{w}^R}\boldsymbol{\sigma}^R\cdot\vec{n},\phi\right\rangle_{\partial\Omega_e}=0,\qquad\forall\phi\in\Pi_N.
\end{align}
Following the same steps as outlined for $\Rt$, one obtains the DGSEM discretization of $\Rtt$ as
\begin{align}\label{eq:DGSEM_wtt}
	\begin{split}
		\left(\vec{\hat{w}}_{ij}\right)_{tt}=\left(\Rtt_h(\vec{w}_h,\boldsymbol{\sigma})\right)_{ij}:=-\frac{1}{{\Jgeo}_{ij}}\bigg[
		&\sum_{\lambda=0}^{{N}}\frac{\partial\hat{\mathbfcal{F}}^1_{\lambda j}}{\partial\vec{\hat{w}}}\hat{\boldsymbol{\sigma}}_{\lambda j}\hat{D}_{i\lambda}
		+\left([\vec{\widetilde f}^*\hat{s}]_{j}^{+{\xi_1}}\hat{\ell}_i(1)+[\vec{\widetilde f}^*\hat{s}]_{j}^{-{\xi_1}}\hat{\ell}_i(-1)\right)\\
		+&\sum_{\mu=0}^{{N}}\frac{\partial\hat{\mathbfcal{F}}^2_{i\mu}}{\partial\vec{\hat w}}\boldsymbol{\sigma}_{i\mu}\hat{D}_{j\mu}
		+\left([\vec{\widetilde f}^*\hat{s}]_{i}^{+{\xi_2}}\hat{\ell}_j(1)+[\vec{\widetilde f}^*\hat{s}]_{i}^{-{\xi_2}}\hat{\ell}_j(-1)\right)
		\bigg]\quad\forall~i,j,
	\end{split}
\end{align}
where the numerical flux $\widetilde{\vec{f}}^*=\widetilde{\vec{f}}^*(\vec{w}^L,\vec{w}^R,\boldsymbol{\sigma}^L,\boldsymbol{\sigma}^R,\vec{n})$ is given by
\begin{align*}
	\widetilde{\vec{f}}^*(\vec{w}^L,\vec{w}^R,\boldsymbol{\sigma}^L,\boldsymbol{\sigma}^R,\vec{n})=\left(\frac{1}{2}\left(\frac{\partial\vec{F}(\vec{w}^L)}{\partial\vec{w}}\boldsymbol{\sigma}^L+\frac{\partial\vec{F}(\vec{w}^R)}{\partial\vec{w}}\boldsymbol{\sigma}^R\right)+\lambda\left(\boldsymbol{\sigma}^L-\boldsymbol{\sigma}^R\right)\right)\cdot\vec{n}.
\end{align*}
Similar as it is done for $\vec{w}^{L/R}$, the values $\boldsymbol{\sigma}^{L/R}$ are obtained by evaluation of the solution polynomial at the cell edges.
This completes the spatial discretization and allows us now to formulate the non-linear equation system to be solved for the time stepping procedure.

\subsection{Solving the (Non-)Linear System of Equations}
\subsubsection{The Linear System for the Two-Derivative DGSEM}
As we have introduced the artificial quantity $\boldsymbol{\sigma}:=\Rt_h(\vec{w})$, we define an extended state vector $\boldsymbol{X}:=(\hat{\ww},\hat{\boldsymbol{\sigma}})^T$, where $\hat{\ww}$ and $\hat{\boldsymbol{\sigma}}$ contain the coefficients of the polynomial basis of $\ww$ and $\boldsymbol{\sigma}$. Note that for ease of presentation, we omit the hat symbol $(\hat\bullet)$ in the following. For this extended state vector, the residual~\eqref{eq:Newton_Residual} is modified and reads
\begin{align*}
	\mathbfcal{G}(\boldsymbol{X})=\begin{pmatrix}\mathbfcal{G}_1(\boldsymbol{X})\\\mathbfcal{G}_2(\boldsymbol{X})\end{pmatrix}:=\begin{pmatrix}\ww\\\boldsymbol{\sigma}\end{pmatrix}-\begin{pmatrix}\vec{b}\\0\end{pmatrix}-\begin{pmatrix}\vec{g}(\ww)\\\Rt_h(\ww)\end{pmatrix}\overset{!}{=}0.
\end{align*}
For the definitions of $\ww$, $\vec{b}$ and $\vec{g}(\ww)$ see Eq.~\eqref{eq:nonlinear_predictor} and Eq.~\eqref{eq:nonlinear_corrector} -- of course with $\Rt$ and $\Rtt$ replaced by their discrete counterparts $\Rt_h$ and $\Rtt_h$, respectively.
To solve for $\boldsymbol{X}$, Newton's method is applied which consists of $\mathtt{r}$ iterative solves
\begin{align}\label{eq:Newton}
\begin{split}
	\frac{\partial{\mathbfcal{G}}(\boldsymbol{X}^\mathtt{r})}{\partial\boldsymbol{X}}\cdot\Delta\boldsymbol{X}&=-{\mathbfcal{G}}(\boldsymbol{X}^\mathtt{r})\\
	\boldsymbol{X}^{\mathtt{r}+1}&=\boldsymbol{X}^\mathtt{r}+\Delta\boldsymbol{X},
\end{split}
\end{align}
with the Newton increment $\Delta\boldsymbol{X}:=(\Delta\ww, \Delta\boldsymbol{\sigma})^T$. To ease presentation, we will drop the superscript $\mathtt{r}$ in the following.
The matrix-vector product of the arising linear system for each Newton's iteration with the system matrix $\mathbfcal{J}$ is given by
\begin{align}\label{eq:nonlinear_eqsys}
\begin{split}
	\mathbfcal{J}\Delta\boldsymbol{X}:=\frac{\partial{\mathbfcal{G}}(\boldsymbol{X})}{\partial\boldsymbol{X}}\cdot\Delta\boldsymbol{X}&=\Delta\boldsymbol{X}-\frac{\partial}{\partial\boldsymbol{X}}\begin{pmatrix}\alpha_1\dt\Rt_h(\ww)-\frac{\alpha_2\dt^2}{2}\Rtt_h(\ww,\boldsymbol{\sigma})\\\Rt_h(\ww)\end{pmatrix}\Delta\boldsymbol{X}\\&=\begin{pmatrix}\Id-\alpha_1\dt\frac{\partial\Rt_h}{\partial\ww}+\frac{\alpha_2\Delta t^2}{2}\frac{\partial\Rtt_h}{\partial\ww}&\frac{\alpha_2\Delta t^2}{2}\frac{\partial\Rtt_h}{\partial\boldsymbol{\sigma}}\\-\frac{\partial\Rt_h}{\partial\ww}&\Id\end{pmatrix}\begin{pmatrix}\Delta\ww\\\Delta\boldsymbol{\sigma}\end{pmatrix}.
\end{split}
\end{align}
Due to the definition of $\boldsymbol{\sigma}$, see Eq.~\eqref{eq:Sigma}, two important observations can be made for Eq.~\eqref{eq:nonlinear_eqsys}:
\begin{itemize}
	\item A comparison of Eq.~\eqref{eq:DGSEM_wt} and Eq.~\eqref{eq:DGSEM_wtt} directly shows $\frac{\partial \Rtt_h}{\partial\boldsymbol{\sigma}}\equiv\frac{\partial\Rt_h}{\partial\ww}$.
	\item If the flux is linear, the Hessian contribution is zero, i.e. $\frac{\partial\Rtt_h}{\partial\ww}=0$.
\end{itemize}
Note that if the discretized system of equations~\eqref{eq:hyperbolic} has a linear flux, Newton's method is not required and one can directly solve the linear equation for $\boldsymbol{X}$ instead of $\Delta\boldsymbol{X}$. Nevertheless, utilizing Newton's method can have a beneficial influence on the achieved accuracy, see e.g.~\cite{Higham2002}. Therefore, we use Newton's method regardless if the system is linear or non-linear.

\subsubsection{Notes on the Practical Implementation}
For the implementation of the system matrix $\mathbfcal{J}$, one has to calculate the matrices $\frac{\partial\Rt_h}{\partial\ww}$ and $\frac{\partial\Rtt_h}{\partial\ww}$ and assemble them into $\mathbfcal{J}$. The former is the standard system matrix of the DGSEM, where a detailed derivation of the inner-element dependencies can be found in~\cite{Vangelatos,Zeifang2020diss}. For the calculation of the latter, all routines of the standard system matrix can be reused. The only difference is that the physical flux Jacobi matrices have to be substituted by the physical flux Hessian.
We solve the linear system with a GMRES method using the PETSc library~\cite{petsc3}.

\subsubsection{A Novel Preconditioner for the Extended Linear System: $\text{BJ}_\text{ext}$}
For a better convergence of the GMRES method, a preconditioner can be applied. Typical preconditioners are left ILU(0) preconditioning or an element-wise block-Jacobi preconditioner (BJ).
Alternatively, the special structure of the linear system, see Eq.~\eqref{eq:nonlinear_eqsys}, can be exploited to construct a problem-tailored preconditioner. In the following, two slightly different preconditioners are presented, which either take the Hessian contribution in $\mathbfcal{J}$ into account ($\text{BJ}^H_\text{ext}$) or neglect this contribution ($\text{BJ}_\text{ext}$).

\begin{figure}[ht]
	\centering
	\includegraphics[width=0.8\linewidth]{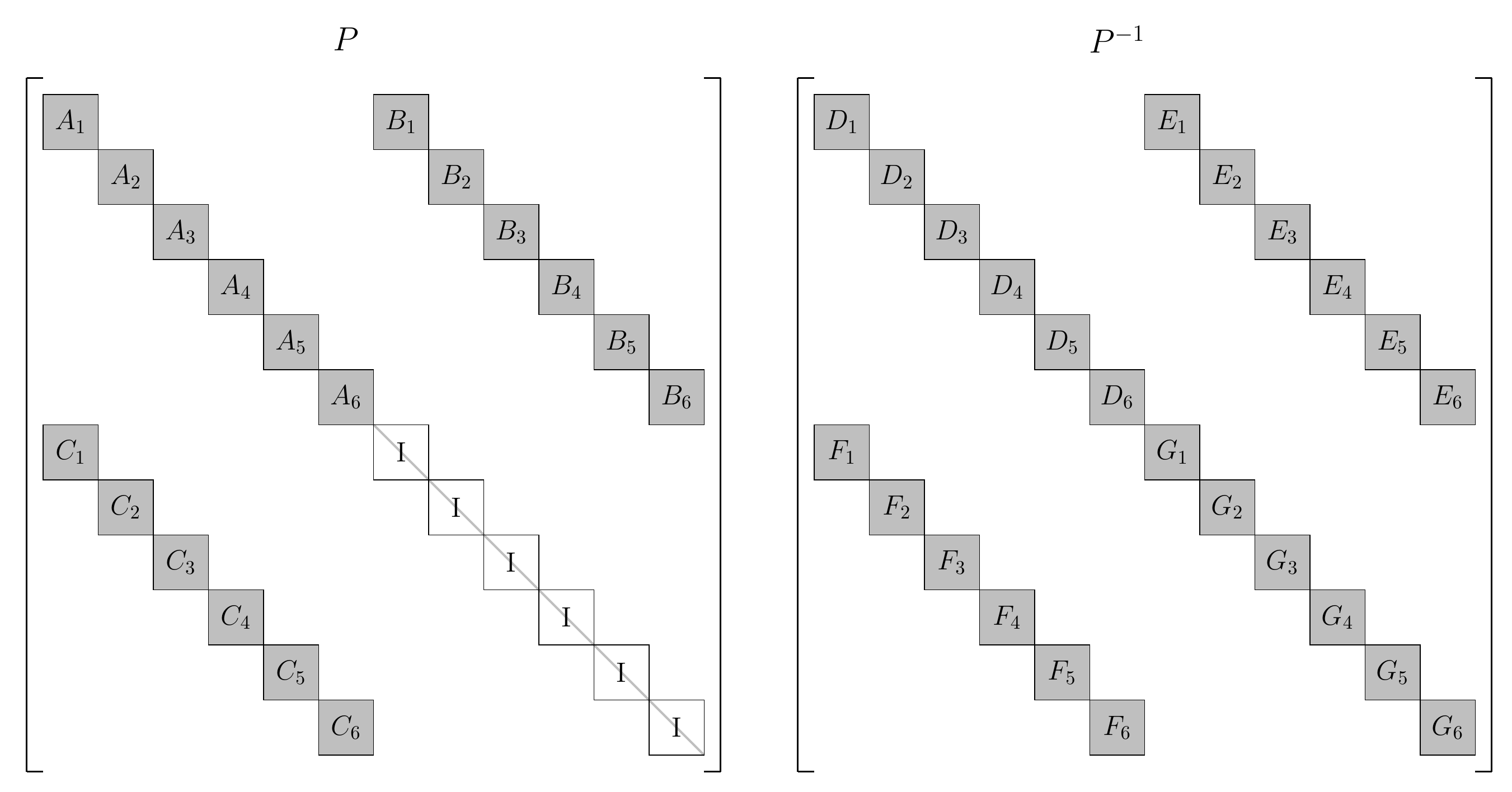}
	\caption{Sketch of extended block-Jacobi preconditioner matrix $P$ (left) and inverse of preconditioner matrix $P^{-1}$ (right) for a setup with six spatial discretization elements ($n_E=6$).}\label{fig:precond}
\end{figure}

\paragraph{$\text{BJ}^H_\text{ext}$ Preconditioner}
We start with the idea of using an element-wise block-Jacobi preconditioner and derive an extended block-Jacobi preconditioner.
That means that we build up the contributions to the system matrix without considering dependencies of neighboring elements, i.e. only the element-internal influence is taken into account. This preconditioner matrix $P$ is visualized in Fig.~\ref{fig:precond} on the left.
For the construction of $P$, the four different blocks are given by
\begin{align*}
	A_i=A_i^H:=\left.\Id-\alpha_1\dt\frac{\partial\Rt_h}{\partial\ww}+\frac{\alpha_2\dt^2}{2}\frac{\partial\Rtt_h}{\partial\ww}\right|_i,\quad 
	B_i:=\left.\frac{\alpha_2\dt^2}{2}\frac{\partial\Rtt_h}{\partial\boldsymbol{\sigma}}\right|_i=\frac{\alpha_2\dt^2}{2}\left.\frac{\partial\Rt_h}{\partial\ww}\right|_i,\quad
	C_i:=\left.-\frac{\partial\Rt_h}{\partial\ww}\right|_i,\quad i=1,\dots,{n_E},
\end{align*}
with ${n_E}$ denoting the number elements of the spatial discretization.
As we have neglected all element-neighbor dependencies and due to the special structure of $P$, we can directly find the inverse $P^{-1}$ which also consists of four different types of blocks, visualized in Fig.~\ref{fig:precond}, via
\begin{align*}
	D_i^H=\left(A_i-B_iC_i\right)^{-1},\quad
	E_i^H=-\left(A_i-B_iC_i\right)^{-1}\cdot B_i,\quad
	F_i^H=-C_i\cdot\left(A_i-B_iC_i\right)^{-1}~\text{and}\quad
	G_i^H=\Id+C_i\cdot\left(A_i-B_iC_i\right)^{-1}\cdot B_i,
\end{align*}
where the superscript $(\bullet)^H$ indicates that the Hessian contribution has been considered in $P$.
These block matrix contributions can be calculated independently for all elements. Moreover, one can see that they require only one matrix inversion of one element-sized matrix per element, viz.
\begin{align*}
	\left(A_i-B_iC_i\right)^{-1}=\left({\left.\Id-\alpha_1\dt\frac{\partial\Rt_h}{\partial\ww}\right|_i+\frac{\alpha_2\dt^2}{2}{\left.{\frac{\partial\Rtt_h}{\partial\ww}}\right|_i}}+\frac{\alpha_2\dt^2}{2}\left({\left.{\frac{\partial\Rt_h}{\partial\ww}}\right|_i}\right)^2 \right)^{-1}.
\end{align*}
As these matrices are relatively small, we calculate those inverses via LU-decomposition.

\paragraph{$\text{BJ}_{\text{ext}}$ Preconditioner}
Alternatively, we define the preconditioner $\text{BJ}_{\text{ext}}$ which neglects the Hessian in $P$, i.e.
\begin{align*}
	A_i:=\left.\Id-\alpha_1\dt\frac{\partial\Rt_h}{\partial\ww}\right|_i,\quad i=1,\dots,{n_E},
\end{align*}
This reveals that only $\left.\frac{\partial\Rt_h}{\partial\ww}\right|_i$ has to be calculated, which is the same as for standard one-derivative methods such as the implicit Euler method.
Again, the inverse of the preconditioner $P^{-1}$ consists of four different types of blocks which are now given by
\begin{align*}
	D_i=\left(A_i-B_iC_i\right)^{-1},\quad
	E_i=-B_i\cdot\left(A_i-B_iC_i\right)^{-1},\quad
	F_i=-C_i\cdot\left(A_i-B_iC_i\right)^{-1}~\text{and}\quad
	G_i=A_i\cdot\left(A_i-B_iC_i\right)^{-1}.
\end{align*}
For the calculation of those block matrices, we have exploited the fact that there holds $A_iB_i=B_iA_i$ and $B_iC_i=C_iB_i$, if the Hessian contribution in $A_i$ is neglected.
Again, only one inverse per element has to be calculated, which is
\begin{align*}
	\left(A_i-B_iC_i\right)^{-1}=\left({\left.\Id-\alpha_1\dt\frac{\partial\Rt_h}{\partial\ww}\right|_i+\frac{\alpha_2\dt^2}{2}\left({\left.{\frac{\partial\Rt_h}{\partial\ww}}\right|_i}\right)^2}\right)^{-1}.
\end{align*}
This preconditioner has two advantages compared to the $\text{BJ}^H_{\text{ext}}$ preconditioner: During the calculation of $E_i$ and $G_i$, the number of matrix-matrix multiplication can be reduced by one. This decreases the building costs of the preconditioner. Additionally, the implementation of the preconditioner is simpler as the Hessian does not have to be implemented. The drawback compared to the $\text{BJ}^H_{\text{ext}}$ preconditioner is that $P$ is a slightly worse approximation to $\mathbfcal{J}$.
Both alternatives will be compared in Sec.~\ref{sec:NumResEuler}.

\subsubsection{Reducing the Problem Size with the Schur Complement}
Another option to reduce the computational and implementational effort is to exploit the special structure of $\mathbfcal{J}$ by forming a Schur complement to reduce the problem size. Instead of solving Eq.~\eqref{eq:Newton}, the problem can be reduced to
\begin{align*}
	\mathbfcal{J}_{\text{Schur}}^\mathtt{r}\Delta\vec{{w}}&=-\mathbfcal{G}_1(\boldsymbol{X}^\mathtt{r})+\frac{\alpha_2\dt^2}{2}\frac{\partial\Rtt_h}{\partial\boldsymbol{\sigma}}\cdot\mathbfcal{G}_2(\boldsymbol{X}^\mathtt{r})\\
	\Delta\boldsymbol{\sigma}&=-\mathbfcal{G}_2(\boldsymbol{X}^\mathtt{r})+\frac{\partial\Rt_h}{\partial\ww}\Delta\vec{w}\\
	\boldsymbol{X}^{\mathtt{r}+1}&=\boldsymbol{X}^\mathtt{r}+\Delta\boldsymbol{X},
\end{align*}
where $\mathtt{r}$ again denotes the index of the Newton iterate, which will be omitted in the following.
The Schur complement system matrix $\mathbfcal{J}_{\text{Schur}}$ is given by
\begin{align}\label{eq:nonlinear_eqsys_Schur}
	\mathbfcal{J}_{\text{Schur}}:=\Id-\alpha_1\dt\frac{\partial\Rt_h}{\partial\ww}+\frac{\alpha_2\dt^2}{2}\frac{\partial\Rtt_h}{\partial\ww}+\frac{\alpha_2\dt^2}{2}\frac{\partial\Rtt_h}{\partial\boldsymbol{\sigma}}\frac{\partial\Rt_h}{\partial\ww}=\Id-\alpha_1\dt\frac{\partial\Rt_h}{\partial\ww}+\frac{\alpha_2\dt^2}{2}\frac{\partial\Rtt_h}{\partial\ww}+\frac{\alpha_2\dt^2}{2}\left(\frac{\partial\Rt_h}{\partial\ww}\right)^2.
\end{align}
The advantage of this approach is that while $\mathbfcal{J}\in\mathbb{R}^{2\cdot\text{size}(\vec{w})\times2\cdot\text{size}(\vec{w})}$, the Schur complement system matrix is only $\mathbfcal{J}_{\text{Schur}}\in\mathbb{R}^{\text{size}(\vec{w})\times\text{size}(\vec{w})}$. One drawback of this approach is that the square of DGSEM's system matrix has to be calculated.

In the following, we first consider the linear solver's convergence for the case with and without reducing the problem size with the Schur complement. This allows us to choose the better suited approach. Moreover, the effectiveness of the different preconditioning strategies is investigated. Finally, we illustrate the capabilities of the novel schemes by showing the experimental order of convergence.

	\section{Numerical Investigations}\label{sec:NumRes}
	In this section, some properties of the novel scheme are evaluated. We investigate if reducing the size of the system matrix with a Schur complement is beneficial and evaluate the best suited preconditioner. Moreover, it is shown how a matrix-free approach can be used to facilitate the implementation of the novel method. The linear scalar advection and the non-linear Euler equations of gas dynamics are taken as prototypical hyperbolic PDEs.

\subsection{Convergence of the Linear Solver}
\subsubsection{Linear Scalar Advection}\label{sec:DG_LinAdv_Convergence}
We start by considering the linear scalar advection equation
\begin{align*}
	w_t+\nabla_x\cdot(\vec{a}w)=0,
\end{align*}
with the advection velocity $\vec{a}$ being a constant vector.
The initial condition $w(\vec{x},t=0)=\sin(\pi\vec{x})$ is chosen, leading to the solution
\begin{align}\label{eq:linadv_ini}
	w(\vec{x},t) = \sin\left(\pi\left(\vec{x}-\vec{a} t\right)\right).
\end{align}
The advection velocity is set to $\vec{a}=(0.3,0.3)^T$ on the two dimensional domain $\Omega=[-1,1]^2$. For the numerical flux function, see Eq.~\eqref{eq:Riemann}, we choose $\lambda=|\vec{a}\cdot\vec{n}|$, with $\vec{n}$ being the normal vector of the element's face.

The domain is discretized with ${n_E}=16^2$ elements with a polynomial degree of the ansatz functions of ${N}=5$. To get a better understanding of the convergence behavior of the novel scheme, we consider the required amount of GMRES iterations with left ILU(0) preconditioning for three different cases and a varying timestep size:
\begin{itemize}
	\item $\mathbfcal{J}^{\text{predictor}}=\begin{pmatrix}\Id-\frac{\dt}{2}\frac{\partial\Rt_h}{\partial\ww}&\frac{\Delta t^2}{12}\frac{\partial\Rt_h}{\partial\ww}\\-\frac{\partial\Rt_h}{\partial\ww}&\Id\end{pmatrix}$, which arises from the $\method{q}{\km}$ predictor, see Eq.~\eqref{eq:nonlinear_eqsys},
	\item $\mathbfcal{J}^{\text{predictor}}_{\text{Schur}}=\Id-\frac{\dt}{2}\frac{\partial\Rt_h}{\partial\ww}+\frac{\dt^2}{12}\left(\frac{\partial\Rt_h}{\partial\ww}\right)^2$, which arises from the $\method{q}{\km}$ predictor, using the Schur complement, see Eq.~\eqref{eq:nonlinear_eqsys_Schur}, and
	\item $\mathbfcal{J}^{\text{Euler}}=\Id-\dt\frac{\partial\Rt_h}{\partial\ww}$, which arises from the first order implicit Euler method as a comparison.
\end{itemize}
Note that even though the problem is linear, we use Newton's method to solve the problem in order to have the same setup as for the non-linear case.
The relative convergence criteria are chosen to be $\eps_{\text{GMRES}}=10^{-3}$ and $\eps_{\text{Newton}}=10^{-8}$. The maximum number of Krylov subspaces of the GMRES method is chosen to be large enough such that no restart of the GMRES method is performed.

\pgfplotstableset{col sep=comma}
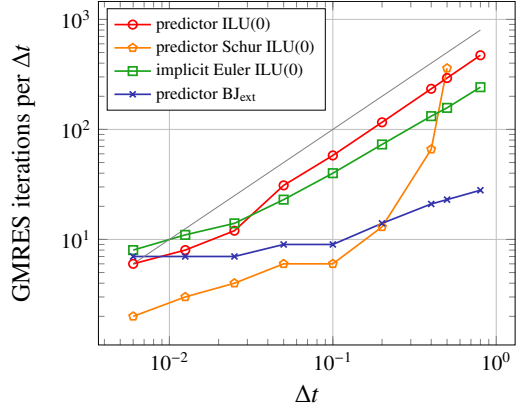
\begin{figure}[ht]
	\centering
        \ifthenelse{\boolean{compilefromscratch}}{
            \tikzsetnextfilename{linadv_condition_convergence_2d}
            \begin{tikzpicture}[scale=0.8]
            \begin{loglogaxis}[cycle list name=epslist,xlabel={$\dt$},ylabel={GMRES iterations per $\dt$} ,grid=major,legend style={at={(0.02,0.98)},anchor= north west,font=\footnotesize,fill opacity=1,text opacity=1},title={},label style={font=\large},title style={font=\large},legend cell align={left}]
            \addplot table[skip first n=1,x expr={\thisrowno{0}}, y expr={\thisrowno{2}}] {./figures/linadv/linadv_condition_convergence_2d.csv};
            \addplot table[skip first n=1,x expr={\thisrowno{0}}, y expr={\thisrowno{5}}] {./figures/linadv/linadv_condition_convergence_2d.csv};
            \addplot table[skip first n=1,x expr={\thisrowno{0}}, y expr={\thisrowno{8}}] {./figures/linadv/linadv_condition_convergence_2d.csv};
            \pgfplotsset{cycle list shift=1}
            \addplot table[skip first n=1,x expr={\thisrowno{0}}, y expr={\thisrowno{10}}]  {./figures/linadv/linadv_condition_convergence_2d.csv};
            \addplot[color=gray] coordinates {(0.006,6) (0.8,800)};
            \legend{predictor ILU(0),predictor Schur ILU(0), implicit Euler ILU(0),predictor $\text{BJ}_{\text{ext}}$}
            \end{loglogaxis}
            \end{tikzpicture}
        }
        {
            \includegraphics{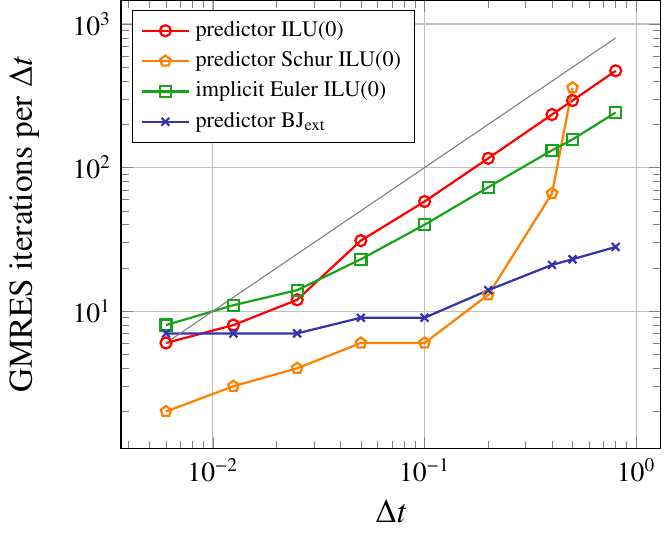}
        }
		\caption{Linear scalar advection equation: required iterations of linear solver to complete one timestep for varying timestep sizes; gray line indicates $\text{slope}=1$.}\label{fig:LinAdv_CondConv}
\end{figure}
Fig.~\ref{fig:LinAdv_CondConv} shows the required GMRES iterations per timestep for the three different schemes. The gray line is a reference which indicates slope one. That means, if the slope of a scheme is lower than the reference, one can expect an acceleration of the calculation when choosing a larger timestep size. If the slope is larger than the reference, a speed-up can only be achieved if the reduction in computational costs for setting up the Jacobian and the preconditioner counterbalances the increased costs due to the increased amount of iterations. First, we compare the implicit Euler method and the predictor of the $\method{q}{\km}$ method with ILU(0) preconditioning each. One can see that the slope for the predictor is slightly steeper than for the implicit Euler scheme and is very close to one. The slope and the number of iterations of the predictor can be reduced drastically by using the $\text{BJ}_{\text{ext}}$ preconditioner. Reducing the problem size with the Schur complement has only a favorable influence for small timesteps. For larger timesteps the iteration count suddenly increases very rapidly and leads to non-convergence for large timesteps. Note that we have observed the same qualitative behavior when using a standard element-wise BJ preconditioner instead of ILU(0) preconditioning for the Schur complement matrix. Summing up, choosing the $\text{BJ}_{\text{ext}}$ preconditioner for the system arising from the $\method{q}{\km}$ predictor without reducing the system size with a Schur complement seems to be the a good choice for a wide range of timestep sizes. The Schur complement discretization is only beneficial for relatively small timesteps.

\subsubsection{Euler Equations}\label{sec:NumResEuler}
Next, we consider the Euler equations of gas dynamics in two dimensions as an example for a non-linear equation system, which read
\begin{align}\label{eq:Euler}
	\ww_t+\nabla_x\cdot\vec{F}(\ww)=0,\quad\text{with}\quad\ww=\begin{pmatrix}\rho\\\rho\vec{v}\\E\end{pmatrix}\quad\text{and}\quad\vec{F}(\ww)=\begin{pmatrix}\rho\vec{v}\\\rho\vec{v}\otimes\vec{v}+\frac{1}{\eps^2}p\cdot\Id\\\vec{v}(E+p),
	\end{pmatrix},
\end{align}
with the density $\rho$, velocity $\vec{v}=(v_1,v_2)^T$, energy $E$ and the reference Mach number $\eps$. The pressure $p$ is calculated via the equation of state for a perfect gas
\begin{align*}
	p=(\gamma-1)\left(E-\frac{\eps^2}{2}\rho\|\vec{v}\|^2\right),
\end{align*}
with the isentropic coefficient $\gamma=1.4$. For the numerical flux function (Eq.~\eqref{eq:Riemann}), following~\cite{KS17,RSIMEXFullEuler} we choose $\lambda=\diag(\frac{1}{\eps},1,1,\frac{1}{\eps})$.
The considered test setup is an extension of the one for the linear scalar advection, see Eq.~\eqref{eq:linadv_ini}.
The initial conditions are hence given by
\begin{align}\label{eq:euler_convtest}
\begin{split}
	\rho(\vec{x},t)&=1+0.3\sin\left(\pi(\vec{x}-\vec{a}t)\right), \qquad
	\vec{v}=\vec{a}, \qquad
	p=1,
\end{split}
\end{align}
with $\vec{a}:=(0.3.0.3)^T$ on domain $\Omega=[-1,1]^2$. The domain is discretized with $n_E=16^2$ elements with a polynomial degree of the ansatz functions of ${N}=5$. The final time is set to $T_{\text{end}}=0.8$ and convergence tolerances are chosen to be $\eps_{\text{GMRES}}=10^{-3}$ and $\eps_{\text{Newton}}=10^{-8}$. The GMRES method performs a restart if no convergence has been reached after $700$ iterations.

Similar as for the linear case, see Sec.~\ref{sec:DG_LinAdv_Convergence}, we consider the three different system matrices $\mathbfcal{J}^{\text{predictor}}$, $\mathbfcal{J}^{\text{predictor}}_{\text{Schur}}$ and $\mathbfcal{J}^{\text{Euler}}$ to evaluate the convergence properties of the schemes for the non-linear case. Differently as for the linear case, the Hessian contribution is not zero (${\partial\Rtt_h}/{\partial\ww}\ne0$). We therefore additionally investigate the influence on the required iterations when neglecting this contribution in the system matrix.

\paragraph{Influence of Chosen System Matrix}
\pgfplotstableset{col sep=comma}

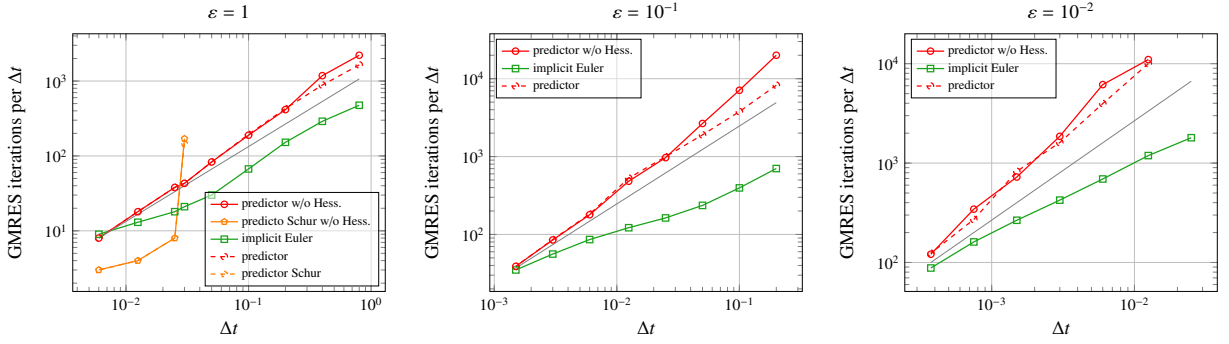
\begin{figure}[ht]
	\centering
	\begin{tabular}{ccc}
        \ifthenelse{\boolean{compilefromscratch}}
        {
            \tikzsetnextfilename{euler_condition_convergence_2d}
            \begin{tikzpicture}[scale=0.6]
            \begin{loglogaxis}[cycle list name=epslist,xlabel={$\dt$},ylabel={GMRES iterations per $\dt$} ,grid=major,legend style={at={(0.98,0.02)},anchor= south east,font=\footnotesize,fill opacity=0.5,text opacity=1},title={$\eps=1$},label style={font=\large},title style={font=\large},legend cell align={left}]
            \addplot table[skip first n=1,x expr={\thisrowno{0}}, y expr={\thisrowno{2}}] {./figures/euler/euler_condition_convergence_2d.csv};
            \addplot table[skip first n=1,x expr={\thisrowno{0}}, y expr={\thisrowno{8}}] {./figures/euler/euler_condition_convergence_2d.csv};
            \addplot table[skip first n=1,x expr={\thisrowno{0}}, y expr={\thisrowno{11}}] {./figures/euler/euler_condition_convergence_2d.csv};
            \pgfplotsset{cycle list shift=3}
            \addplot table[skip first n=1,x expr={\thisrowno{0}}, y expr={\thisrowno{5}}] {./figures/euler/euler_condition_convergence_2d.csv};
            \addplot table[skip first n=1,x expr={\thisrowno{0}}, y expr={\thisrowno{16}}] {./figures/euler/euler_condition_convergence_2d.csv};
            \addplot[color=gray] coordinates {(0.006,8) (0.8,1067)};
            \legend{predictor w/o Hess.,predicto Schur w/o Hess., implicit Euler, predictor, predictor Schur}
            \end{loglogaxis}
            \end{tikzpicture}&
            \tikzsetnextfilename{euler_condition_convergence_1e-1}
            \begin{tikzpicture}[scale=0.6]
            \begin{loglogaxis}[cycle list name=epslist,xlabel={$\dt$},ylabel={GMRES iterations per $\dt$} ,grid=major,legend style={at={(0.02,0.98)},anchor= north west,font=\footnotesize,fill opacity=0.5,text opacity=1},title={$\eps=10^{-1}$},label style={font=\large},title style={font=\large},legend cell align={left}]
            \addplot table[skip first n=1,x expr={\thisrowno{0}}, y expr={\thisrowno{1}}] {./figures/euler/euler_convergence_1e-1.csv};
            \pgfplotsset{cycle list shift=1}
            \addplot table[skip first n=1,x expr={\thisrowno{0}}, y expr={\thisrowno{3}}] {./figures/euler/euler_convergence_1e-1.csv};
            \pgfplotsset{cycle list shift=4}
            \addplot table[skip first n=1,x expr={\thisrowno{0}}, y expr={\thisrowno{7}}] {./figures/euler/euler_convergence_1e-1.csv};
            \addplot[color=gray] coordinates {(0.0015,37) (0.2,4933)};
            \legend{predictor w/o Hess., implicit Euler, predictor}
            \end{loglogaxis}
            \end{tikzpicture}&
            \tikzsetnextfilename{euler_condition_convergence_1e-2}
            \begin{tikzpicture}[scale=0.6]
            \begin{loglogaxis}[cycle list name=epslist,xlabel={$\dt$},ylabel={GMRES iterations per $\dt$} ,grid=major,legend style={at={(0.02,0.98)},anchor= north west,font=\footnotesize,fill opacity=0.5,text opacity=1},title={$\eps=10^{-2}$},label style={font=\large},title style={font=\large},legend cell align={left},ymax=2e4]
            \addplot table[skip first n=1,x expr={\thisrowno{0}}, y expr={\thisrowno{1}}] {./figures/euler/euler_convergence_1e-2.csv};
            \pgfplotsset{cycle list shift=1}
            \addplot table[skip first n=1,x expr={\thisrowno{0}}, y expr={\thisrowno{3}}] {./figures/euler/euler_convergence_1e-2.csv};
            \pgfplotsset{cycle list shift=4}
            \addplot table[skip first n=1,x expr={\thisrowno{0}}, y expr={\thisrowno{7}}] {./figures/euler/euler_convergence_1e-2.csv};
            \addplot[color=gray] coordinates {(0.000375,100) (0.025,6667)}; 
            \legend{predictor w/o Hess., implicit Euler, predictor}
            \end{loglogaxis}
            \end{tikzpicture}
        }
        {   
            \includegraphics{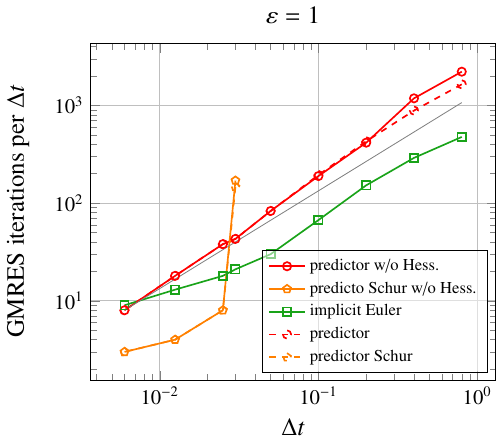} &
            \includegraphics{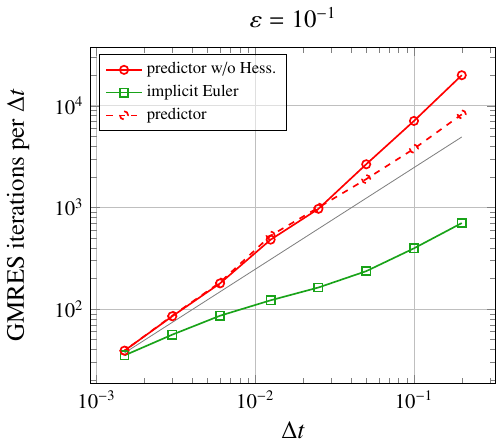}&
            \includegraphics{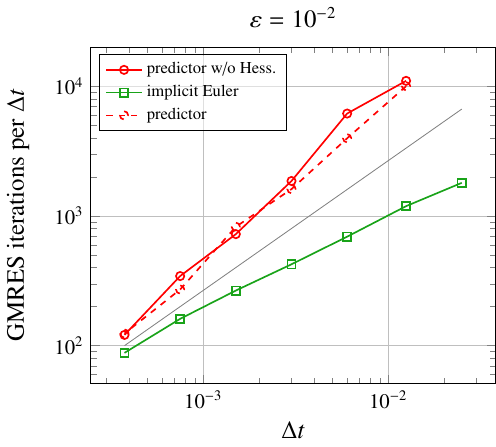} 
        }
	\end{tabular}
	\caption{Euler equations: required iterations of linear solver to complete one timestep for varying timestep sizes for different reference Mach numbers $\eps$ with and without Hessian contribution; gray line indicates $\text{slope}=1$. For all schemes an ILU(0) preconditioner is used. Note that missing data points indicate that no convergence of the linear solver could be obtained within a limit of $7000$ GMRES iterations per Newton step.}\label{fig:Euler_Hessian}
\end{figure}
The required linear iterations for those setups for a varying reference Mach number $\eps$ are visualized in Fig.~\ref{fig:Euler_Hessian}. One can see an almost linear behavior for the $\method{q}{\km}$ predictor and the implicit Euler scheme. Comparing the predictor and the implicit Euler scheme (both with ILU(0) preconditioner), shows that while the slope of the curve for the implicit Euler scheme remains below one that of the predictor is always slightly larger than one. Going to higher stiffnesses ($\eps=10^{-2}$) shows that the slope even increases for the predictor.

Reducing the system size with the Schur complement, reveals a strongly non-linear behavior of the required amount of iterations. Starting from a certain threshold, the required iterations increase rapidly such that no solution is obtained already for moderately large timesteps. In both cases, neglecting the Hessian contribution has only a minor influence on the required iterations for moderately large timesteps. For very large timesteps, the linear scaling of the iteration number can only be obtained when taking the Hessian contribution into account. This can consistently be observed for all considered stiffnesses. Therefore, the Hessian contribution in the system matrix is always taken into account in the remainder of this paper.

\paragraph{Influence of Preconditioner}
Additionally, we investigate the influence of the used preconditioner on the required iterations without using the Schur complement. For that purpose, we repeat the previous simulations with the Hessian contribution choosing different preconditioners. The preconditioners are rebuilt before each Newton iteration.
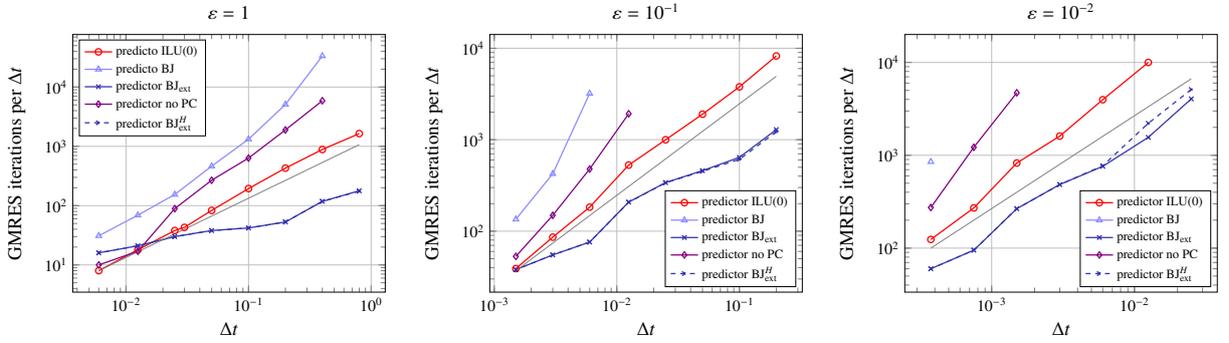
\begin{figure}[ht]
	\centering
	\begin{tabular}{ccc}
        \ifthenelse{\boolean{compilefromscratch}}
        {
            \tikzsetnextfilename{euler_iterations_1e0}
            \begin{tikzpicture}[scale=0.6]
            \begin{loglogaxis}[cycle list name=epslist,xlabel={$\dt$},ylabel={GMRES iterations per $\dt$} ,grid=major,legend style={at={(0.02,0.98)},anchor= north west,font=\footnotesize,fill opacity=1,text opacity=1},title={$\eps=1$},label style={font=\large},title style={font=\large},legend cell align={left}]
            \addplot table[skip first n=1,x expr={\thisrowno{0}}, y expr={\thisrowno{5}}] {./figures/euler/euler_condition_convergence_2d.csv};
            \pgfplotsset{cycle list shift=2}
            \addplot table[skip first n=1,x expr={\thisrowno{0}}, y expr={\thisrowno{13}}] {./figures/euler/euler_condition_convergence_2d.csv};
            \addplot table[skip first n=1,x expr={\thisrowno{0}}, y expr={\thisrowno{24}}] {./figures/euler/euler_condition_convergence_2d.csv};
            \addplot table[skip first n=1,x expr={\thisrowno{0}}, y expr={\thisrowno{26}}] {./figures/euler/euler_condition_convergence_2d.csv};
            \pgfplotsset{cycle list shift=6}
            \addplot table[skip first n=1,x expr={\thisrowno{0}}, y expr={\thisrowno{28}}] {./figures/euler/euler_condition_convergence_2d.csv};
            \addplot[color=gray] coordinates {(0.006,8) (0.8,1067)};
            \legend{predicto ILU(0),predicto BJ, predictor $\text{BJ}_{\text{ext}}$,predictor no PC, predictor $\text{BJ}^H_{\text{ext}}$}
            \end{loglogaxis}
            \end{tikzpicture}&
            \tikzsetnextfilename{euler_iterations_1e-1}
            \begin{tikzpicture}[scale=0.6]
            \begin{loglogaxis}[cycle list name=epslist,xlabel={$\dt$},ylabel={GMRES iterations per $\dt$} ,grid=major,legend style={at={(0.98,0.02)},anchor= south east,font=\footnotesize,fill opacity=1,text opacity=1},title={$\eps=10^{-1}$},label style={font=\large},title style={font=\large},legend cell align={left}]
            \addplot table[skip first n=1,x expr={\thisrowno{0}}, y expr={\thisrowno{7}}] {./figures/euler/euler_convergence_1e-1.csv};
            \pgfplotsset{cycle list shift=2}
            \addplot table[skip first n=1,x expr={\thisrowno{0}}, y expr={\thisrowno{5}}] {./figures/euler/euler_convergence_1e-1.csv};
            \addplot table[skip first n=1,x expr={\thisrowno{0}}, y expr={\thisrowno{13}}] {./figures/euler/euler_convergence_1e-1.csv};
            \addplot table[skip first n=1,x expr={\thisrowno{0}}, y expr={\thisrowno{15}}] {./figures/euler/euler_convergence_1e-1.csv};
            \pgfplotsset{cycle list shift=6}
            \addplot table[skip first n=1,x expr={\thisrowno{0}}, y expr={\thisrowno{17}}] {./figures/euler/euler_convergence_1e-1.csv};
            \addplot[color=gray] coordinates {(0.0015,37) (0.2,4933)};
            \legend{predictor ILU(0),predictor BJ,predictor $\text{BJ}_{\text{ext}}$, predictor no PC,predictor $\text{BJ}^H_{\text{ext}}$}
            \end{loglogaxis}
            \end{tikzpicture}&
            \tikzsetnextfilename{euler_iterations_1e-2}
            \begin{tikzpicture}[scale=0.6]
            \begin{loglogaxis}[cycle list name=epslist,xlabel={$\dt$},ylabel={GMRES iterations per $\dt$} ,grid=major,legend style={at={(0.98,0.02)},anchor= south east,font=\footnotesize,fill opacity=1,text opacity=1},title={$\eps=10^{-2}$},label style={font=\large},title style={font=\large},legend cell align={left},ymax=2e4]
            \addplot table[skip first n=1,x expr={\thisrowno{0}}, y expr={\thisrowno{7}}] {./figures/euler/euler_convergence_1e-2.csv};
            \pgfplotsset{cycle list shift=2}
            \addplot table[skip first n=1,x expr={\thisrowno{0}}, y expr={\thisrowno{5}}] {./figures/euler/euler_convergence_1e-2.csv};
            \addplot table[skip first n=1,x expr={\thisrowno{0}}, y expr={\thisrowno{13}}] {./figures/euler/euler_convergence_1e-2.csv};
            \addplot table[skip first n=1,x expr={\thisrowno{0}}, y expr={\thisrowno{15}}] {./figures/euler/euler_convergence_1e-2.csv};
            \pgfplotsset{cycle list shift=6}
            \addplot table[skip first n=1,x expr={\thisrowno{0}}, y expr={\thisrowno{17}}] {./figures/euler/euler_convergence_1e-2.csv};
            \addplot[color=gray] coordinates {(0.000375,100) (0.025,6667)}; 
            \legend{predictor ILU(0),predictor BJ,predictor $\text{BJ}_{\text{ext}}$, predictor no PC,predictor $\text{BJ}^H_{\text{ext}}$}	
            \end{loglogaxis}
            \end{tikzpicture}
        }
        {
            \includegraphics{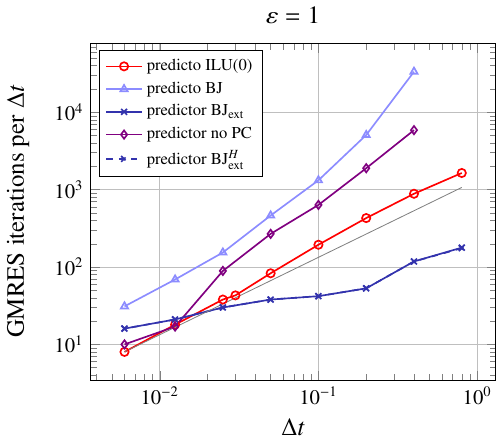} &
            \includegraphics{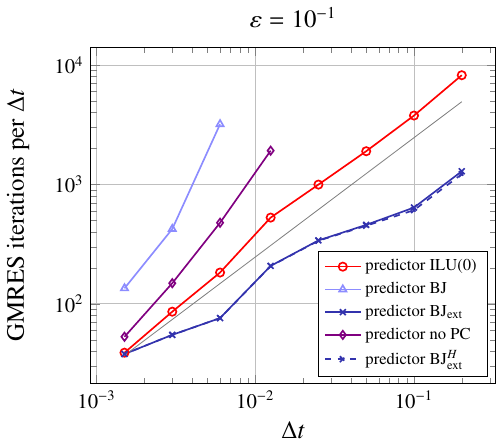}&
            \includegraphics{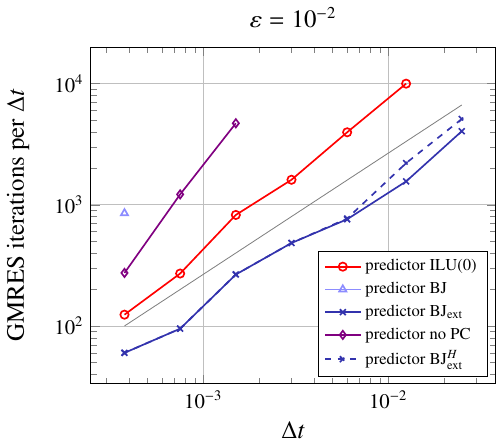}
        }
	\end{tabular}
	\caption{Influence of reference Mach number on required iterations to complete one timestep for different preconditioners for predictor of the $\method{q}{\km}$ scheme; gray line indicates $\text{slope}=1$. Note that missing data points indicate that no convergence of the linear solver could be obtained within a limit of $7000$ GMRES iterations per Newton step. The values for the predictor with $\text{BJ}_{\text{ext}}$ and the predictor with $\text{BJ}^H_{\text{ext}}$ are almost not distinguishable.}\label{fig:Euler_eps}
\end{figure}
We report the results of this series of simulations in Fig.~\ref{fig:Euler_eps}. The figure shows that the novel $\text{BJ}_{\text{ext}}$ and $\text{BJ}^H_{\text{ext}}$ preconditioners are best suited to reduce the required amount of iterations. They even perform better than the ILU(0) preconditioner, which performs better than using no preconditioner. Interestingly, choosing the standard element-wise BJ preconditioner has an unfavorable influence on the required iterations. One can observe that the slope of the curves increases for an increasing stiffness. A similar behavior has already been observed in~\cite{ZKBSM17}. To cure this, one could use adaptive tolerances for the GMRES and Newton's method. Another option is to use asymptotic preserving schemes as it is e.g. done in~\cite{ZKBSM17}.

\begin{figure}[ht]
	\centering
	\begin{tabular}{ccc}
        \ifthenelse{\boolean{compilefromscratch}}
        {
            \tikzsetnextfilename{euler_wallclock_1e0}
            \begin{tikzpicture}[scale=0.6]
            \begin{loglogaxis}[cycle list name=epslist,xlabel={$\dt$},ylabel={wallclocktime [$s$]} ,grid=major,legend style={at={(0.02,0.98)},anchor= north west,font=\footnotesize,fill opacity=1.0,text opacity=1},title={$\eps=1$},label style={font=\large},title style={font=\large},legend cell align={left}]
            \addplot table[skip first n=0,x expr={\thisrowno{6}}, y expr={\thisrowno{7}}] {./figures/euler/wallclocktime_precond.csv};
            \pgfplotsset{cycle list shift=2}
            \addplot table[skip first n=0,x expr={\thisrowno{12}}, y expr={\thisrowno{13}}] {./figures/euler/wallclocktime_precond.csv};
            \addplot table[skip first n=0,x expr={\thisrowno{0}}, y expr={\thisrowno{1}}] {./figures/euler/wallclocktime_precond.csv};
            \addplot table[skip first n=0,x expr={\thisrowno{18}}, y expr={\thisrowno{19}}] {./figures/euler/wallclocktime_precond.csv};
            \pgfplotsset{cycle list shift=6}
            \addplot table[skip first n=0,x expr={\thisrowno{24}}, y expr={\thisrowno{25}}] {./figures/euler/wallclocktime_precond.csv};
            \addplot[color=gray,forget plot] coordinates {(0.006,10) (0.8,1333)};
            \legend{predictor ILU(0),predictor BJ,predictor $\text{BJ}_{\text{ext}}$, predictor no PC,predictor $\text{BJ}^H_{\text{ext}}$}	
            \end{loglogaxis}
            \end{tikzpicture}&
            \tikzsetnextfilename{euler_wallclock_1e-1}
            \begin{tikzpicture}[scale=0.6]
            \begin{loglogaxis}[cycle list name=epslist,xlabel={$\dt$},ylabel={wallclocktime [$s$]} ,grid=major,legend style={at={(0.98,0.02)},anchor= south east,font=\footnotesize,fill opacity=1,text opacity=1},title={$\eps=10^{-1}$},label style={font=\large},title style={font=\large},legend cell align={left}]
            \addplot table[skip first n=0,x expr={\thisrowno{8}}, y expr={\thisrowno{9}}] {./figures/euler/wallclocktime_precond.csv};
            \pgfplotsset{cycle list shift=2}
            \addplot table[skip first n=0,x expr={\thisrowno{14}}, y expr={\thisrowno{15}}] {./figures/euler/wallclocktime_precond.csv};
            \addplot table[skip first n=0,x expr={\thisrowno{2}}, y expr={\thisrowno{3}}] {./figures/euler/wallclocktime_precond.csv};
            \addplot table[skip first n=0,x expr={\thisrowno{20}}, y expr={\thisrowno{21}}] {./figures/euler/wallclocktime_precond.csv};
            \pgfplotsset{cycle list shift=6}
            \addplot table[skip first n=0,x expr={\thisrowno{26}}, y expr={\thisrowno{27}}] {./figures/euler/wallclocktime_precond.csv};
            \addplot[color=gray,forget plot] coordinates {(0.0015,12) (0.2,1600)};
            \legend{predictor ILU(0),predictor BJ,predictor $\text{BJ}_{\text{ext}}$, predictor no PC,predictor $\text{BJ}^H_{\text{ext}}$}	
            \end{loglogaxis}
            \end{tikzpicture}&
            \tikzsetnextfilename{euler_wallclock_1e-2}
            \begin{tikzpicture}[scale=0.6]
            \begin{loglogaxis}[cycle list name=epslist,xlabel={$\dt$},ylabel={wallclocktime [$s$]} ,grid=major,legend style={at={(0.98,0.02)},anchor= south east,font=\footnotesize,fill opacity=0.5,text opacity=1},title={$\eps=10^{-2}$},label style={font=\large},title style={font=\large},legend cell align={left},ymax=2e4]
            \addplot table[skip first n=0,x expr={\thisrowno{10}}, y expr={\thisrowno{11}}] {./figures/euler/wallclocktime_precond.csv};
            \pgfplotsset{cycle list shift=2}
            \addplot table[skip first n=0,x expr={\thisrowno{16}}, y expr={\thisrowno{17}}] {./figures/euler/wallclocktime_precond.csv};
            \addplot table[skip first n=0,x expr={\thisrowno{4}}, y expr={\thisrowno{5}}] {./figures/euler/wallclocktime_precond.csv};
            \addplot table[skip first n=0,x expr={\thisrowno{22}}, y expr={\thisrowno{23}}] {./figures/euler/wallclocktime_precond.csv};
            \pgfplotsset{cycle list shift=6}
            \addplot table[skip first n=0,x expr={\thisrowno{28}}, y expr={\thisrowno{29}}] {./figures/euler/wallclocktime_precond.csv};
            \addplot[color=gray,forget plot] coordinates {(0.000375,14) (0.025,933)}; 
            \legend{predictor ILU(0),predictor BJ,predictor $\text{BJ}_{\text{ext}}$, predictor no PC,predictor $\text{BJ}^H_{\text{ext}}$}		
            \end{loglogaxis}
            \end{tikzpicture}
        }
        {
            \includegraphics{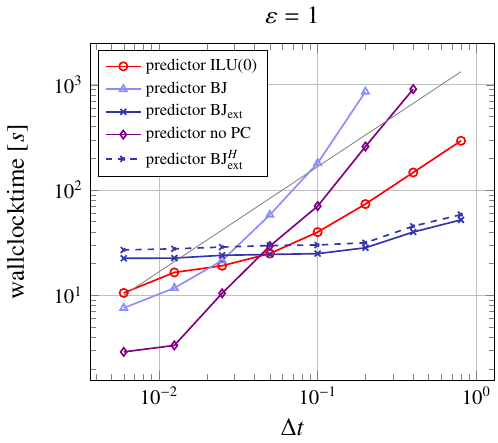}&
            \includegraphics{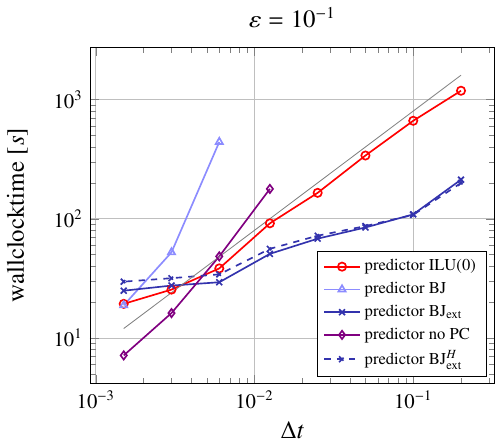}&
            \includegraphics{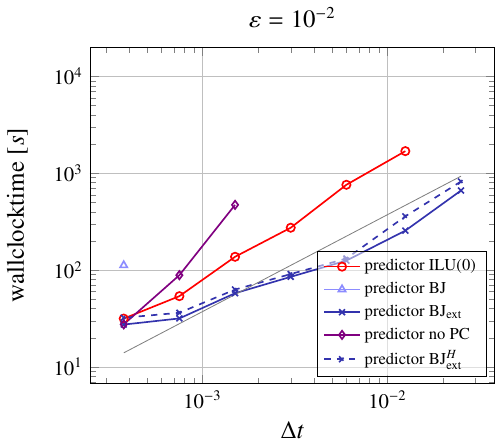}
        }
	\end{tabular}
	\caption{Influence of reference Mach number on required wallclocktime to complete one timestep for different preconditioners for the predictor of the $\method{q}{\km}$ scheme; gray line indicates $\text{slope}=1$. Note that missing data points indicate that no convergence of the linear solver could be obtained within a limit of $7000$ GMRES iterations per Newton step.}\label{fig:Euler_wallclock}
\end{figure}
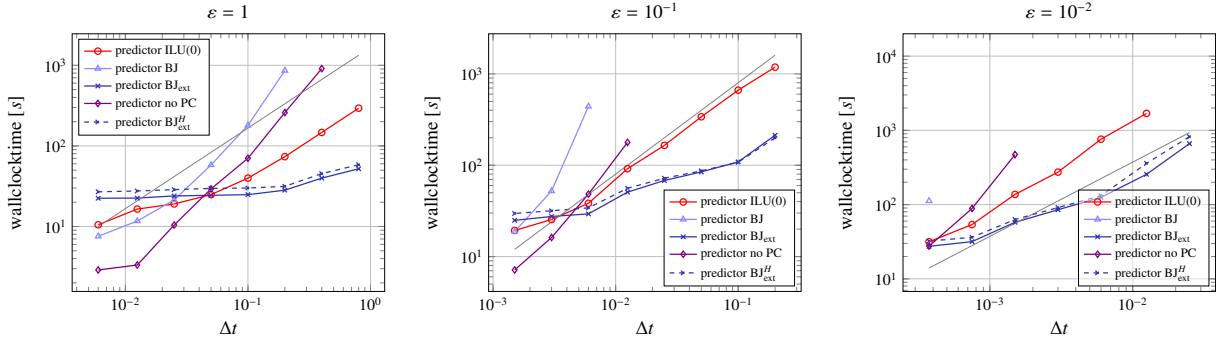

As the efficiency of a preconditioner cannot only be evaluated by its ability to reduce the number of iterations but also by its computational costs, we additionally consider the required wallclocktimes for this set of simulations, which are visualized in Fig.~\ref{fig:Euler_wallclock}.
One can see a very similar behavior as already observed for the linear iterations when considering the slope of the curves: Choosing the $\text{BJ}_{\text{ext}}$/$\text{BJ}^H_{\text{ext}}$ preconditioner gives the most efficient scheme. When the ILU(0) preconditioner is chosen, in only a few cases a speedup can be achieved by selecting a larger timestep size. For the $\text{BJ}_{\text{ext}}$/$\text{BJ}^H_{\text{ext}}$ preconditioner, a (super-) linear scaling is only observed for large timestep sizes.
Regarding the absolute values of the wallclocktime reveals the different costs of the preconditioners. Due to the missing costs of building and applying the preconditioner when choosing no PC, for relatively small timesteps this variant is the most efficient. Still, a good preconditioner such as the novel $\text{BJ}_{\text{ext}}$ preconditioner is required to obtain an efficient scheme for large timesteps.
The $\text{BJ}_{\text{ext}}$/$\text{BJ}^H_{\text{ext}}$ preconditioners are relatively expensive and hence for small timesteps and low stiffnesses are outperformed by the other preconditioners. In a practical application, the large building costs can be reduced by reusing the preconditioner for several Newton steps or even several timesteps.
Comparing the novel $\text{BJ}_{\text{ext}}$ and $\text{BJ}^H_{\text{ext}}$ preconditioner shows that neglecting the Hessian contribution increases the efficiency of the preconditioner. While this has only a small influence on the required iterations, see Fig.~\ref{fig:Euler_eps}, the reduced cost for building the preconditioner has a significant influence on the required wallclocktimes. Therefore, this approach is pursued in the remainder of this paper.

\subsection{Matrix-Free Discretization}
A matrix-free implementation allows to increase the flexibility of a solver and reduces the required memory consumption. Moreover, it has turned out that is is also beneficial in terms of computational time for high order discretizations of complex settings~\cite{Franciolini2017}. For the two-derivative DGSEM discretization it additionally facilitates a parallel distribution of the spatial domain on different processors.
\paragraph{Matrix-Free Approach}
Differently to the previously described matrix-based approach one does not explicitly form the system matrix $\mathbfcal{J}$ and multiply the vector $\Delta\boldsymbol{X}$. Instead, one approximates the matrix vector product $\mathbfcal{J}\Delta\boldsymbol{X}$ via a finite difference, see e.g.~\cite{Knoll2004} for an overview on matrix-free approaches. For the second derivative method, the matrix-vector product given in Eq.~\eqref{eq:nonlinear_eqsys} for a non-linear system is approximated by
\begin{align}\label{eq:MatrixFree_nonlinear}
\begin{split}
	\mathbfcal{J}\Delta\boldsymbol{X}&\approx\begin{pmatrix}
	\Delta\ww-\alpha_1\dt\frac{\Rt_h(\ww+\epsFD^w\Delta\ww)-\Rt_h(\ww)}{\epsFD^w}+\frac{\alpha_2\dt^2}{2}\frac{\Rtt_h(\ww+\epsFD^w\Delta\ww,\boldsymbol{\sigma})-\Rtt_h(\ww,\boldsymbol{\sigma})}{\epsFD^w}+\frac{\alpha_2\dt^2}{2}\frac{\Rtt_h(\ww,\boldsymbol{\sigma}+\epsFD^{\sigma}\Delta\vec{\boldsymbol{\sigma}})-\Rtt_h(\ww,\boldsymbol{\sigma})}{\epsFD^{\sigma}}\\
	-\frac{\Rt_h(\ww+\epsFD^w\Delta\ww)-\Rt_h(\ww)}{\epsFD^w}+\Delta\boldsymbol{\sigma}
	\end{pmatrix},
\end{split}
\end{align}
with $\epsFD^{(\bullet)}$ being a small user-defined value. Inspired by~\cite{Knoll2004}, we set
\begin{align*}
	\epsFD^{(\bullet)}:=\frac{\sqrt{\eps_{\text{machine}}}}{\eps\|\Delta{(\bullet)}\|_2},
\end{align*}
with $\eps_{\text{machine}}$ being approximately machine accuracy. Here, we use the fortran intrinsic \emph{epsilon} function, which gives approximately $\eps_{\text{machine}}\approx2\cdot10^{-16}$.
Additionally, both the reference Mach number and the perturbed quantity are taken into account, so ${(\bullet)}$ is either ${(\bullet)}=\ww$ or ${(\bullet)}=\boldsymbol{\sigma}$.
For a linear system, the matrix-vector product simplifies to 
\begin{align}\label{eq:MatrixFree_linear}
	\begin{split}
		\mathbfcal{J}\Delta\boldsymbol{X}&=\begin{pmatrix}
		\Delta\ww-\alpha_1\dt\frac{\Rt_h(\ww+\epsFD^w\Delta\ww)-\Rt_h(\ww)}{\epsFD^w}+\frac{\alpha_2\dt^2}{2}\frac{\Rt_h(\boldsymbol{\sigma}+\epsFD^{\sigma}\Delta\vec{\boldsymbol{\sigma}})-\Rt_h(\boldsymbol{\sigma})}{\epsFD^{\sigma}}\\
		-\frac{\Rt_h(\ww+\epsFD^w\Delta\ww)-\Rt_h(\ww)}{\epsFD^w}+\Delta\boldsymbol{\sigma}
		\end{pmatrix}=
		\begin{pmatrix}
		\Delta\ww-\alpha_1\dt\Rt_h(\Delta\ww)+\frac{\alpha_2\dt^2}{2}\Rt_h(\Delta\boldsymbol{\sigma})\\
		-\Rt_h(\Delta\ww)+\Delta\boldsymbol{\sigma}
		\end{pmatrix}.
	\end{split}
\end{align}
Eq.~\eqref{eq:MatrixFree_nonlinear} and Eq.~\eqref{eq:MatrixFree_linear} show that two or three different operator evaluations have to be performed per GMRES iteration for the linear or non-linear case, respectively.
\paragraph{Parallel Implementation}
One advantage of using a matrix-free implementation is the straight forward parallelization.
Differently to the matrix-based approach, the Jacobian $\mathbfcal{J}$ does not have to be distributed among the processors. While the parallelization relies on a domain decomposition, and $\mathbfcal{J}$ consists of four blocks, each corresponding to the whole spatial discretization, a parallel matrix-based implementation would require some restructuring of already existing codes. One would reorder the solution vector $\boldsymbol{X}$ such that $\ww$ and $\boldsymbol{\sigma}$ of each element are alternating in the solution vector. Most importantly, the assembling routines and the routines to evaluate the right hand side would have to be rewritten, which is a tedious task.

Instead, we can distribute $\ww$ and $\boldsymbol{\sigma}$ in a similar manner on the processors, such that for both variables the same domain decomposition is used. I.e. one processor contains the same segments of $\ww$ \emph{and} $\boldsymbol{\sigma}$, which are then assembled consecutively in the solution vector. With this, the parallelization of the spatial operator does not have to be changed, see~\cite{Krais2020} for an overview on the parallelization strategy for the FLEXI DGSEM code.

As the ILU(0) preconditioner requires the information of the whole matrix, it is not in line with the matrix-free approach. Differently, the $\text{BJ}_{\text{ext}}$ and $\text{BJ}^H_{\text{ext}}$ preconditioners require only the formation of element-wise small matrices which can be set up independently of each other. Additionally, the application of the preconditioner does not require any communication between different processors. It consists of independent evaluations of matrix-vector products
\begin{align*}
	P^{-1}\cdot\boldsymbol{X}=\begin{pmatrix}
	D_1\cdot\ww_1+E_1\cdot\boldsymbol{\sigma}_1\\
	\vdots\\
	D_{n_E}\cdot\ww_{n_E}+E_{n_E}\cdot\boldsymbol{\sigma}_{n_E}\\
	F_1\cdot\ww_1+G_1\cdot\boldsymbol{\sigma}_1\\
	\vdots\\
	F_{n_E}\cdot\ww_{n_E}+G_{n_E}\cdot\boldsymbol{\sigma}_{n_E}
	\end{pmatrix},
\end{align*}
and is hence well suited for parallel computations. A dedicated evaluation of the parallel performance of the novel method will be addressed in a future work.

\paragraph{Comparison of Matrix-Free and Matrix-Based Implementation}
In order to validate the matrix-free implementation and to illustrate the influence of using an approximation of the Jacobian-vector product on the required iterations, we compare matrix-free and matrix-based simulations. We use the same initialization for the Euler equations (Eq.~\eqref{eq:euler_convtest}) with the same discretization parameters as in Sec.~\ref{sec:NumResEuler}.
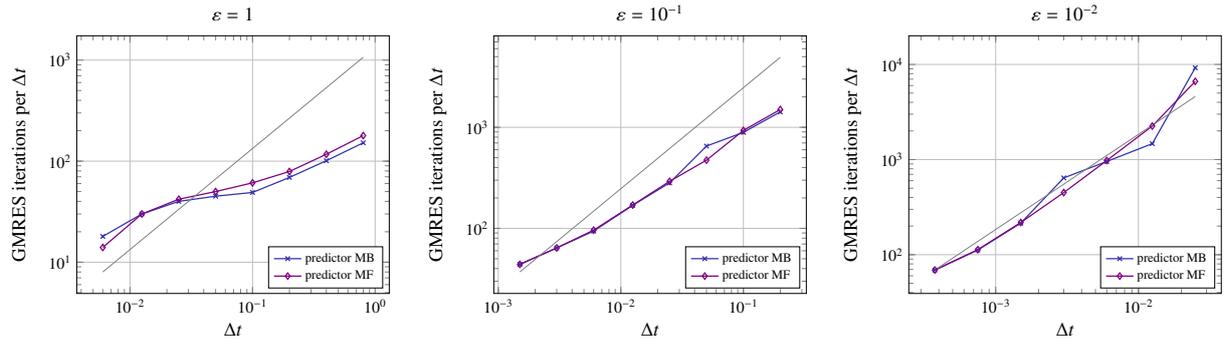
\begin{figure}[ht]
	\centering
	\begin{tabular}{ccc}
        \ifthenelse{\boolean{compilefromscratch}}
        {
            \tikzsetnextfilename{euler_gmres_1e0}
            \begin{tikzpicture}[scale=0.6]
            \begin{loglogaxis}[cycle list name=epslist,xlabel={$\dt$},ylabel={GMRES iterations per $\dt$} ,grid=major,legend style={at={(0.98,0.02)},anchor= south east,font=\footnotesize,fill opacity=1,text opacity=1},title={$\eps=1$},label style={font=\large},title style={font=\large},legend cell align={left}]
            \pgfplotsset{cycle list shift=4}
            \addplot table[skip first n=1,x expr={\thisrowno{0}}, y expr={\thisrowno{22}}] {./figures/euler/euler_condition_convergence_2d.csv};
            \addplot table[skip first n=1,x expr={\thisrowno{0}}, y expr={\thisrowno{20}}] {./figures/euler/euler_condition_convergence_2d.csv};
            \addplot[color=gray] coordinates {(0.006,8) (0.8,1067)};
            \legend{predictor MB,predictor MF}
            \end{loglogaxis}
            \end{tikzpicture}&
            \tikzsetnextfilename{euler_gmres_1e-1}
            \begin{tikzpicture}[scale=0.6]
            \begin{loglogaxis}[cycle list name=epslist,xlabel={$\dt$},ylabel={GMRES iterations per $\dt$} ,grid=major,legend style={at={(0.98,0.02)},anchor= south east,font=\footnotesize,fill opacity=1,text opacity=1},title={$\eps=10^{-1}$},label style={font=\large},title style={font=\large},legend cell align={left}]
            \pgfplotsset{cycle list shift=4}
            \addplot table[skip first n=1,x expr={\thisrowno{0}}, y expr={\thisrowno{11}}] {./figures/euler/euler_convergence_1e-1.csv};
            \addplot table[skip first n=1,x expr={\thisrowno{0}}, y expr={\thisrowno{9}}] {./figures/euler/euler_convergence_1e-1.csv};
            \addplot[color=gray] coordinates {(0.0015,37) (0.2,4933)};
            \legend{predictor MB,predictor MF}
            \end{loglogaxis}
            \end{tikzpicture}&
            \tikzsetnextfilename{euler_gmres_1e-2}
            \begin{tikzpicture}[scale=0.6]
            \begin{loglogaxis}[cycle list name=epslist,xlabel={$\dt$},ylabel={GMRES iterations per $\dt$} ,grid=major,legend style={at={(0.98,0.02)},anchor= south east,font=\footnotesize,fill opacity=1,text opacity=1},title={$\eps=10^{-2}$},label style={font=\large},title style={font=\large},legend cell align={left},ymax=2e4]
            \pgfplotsset{cycle list shift=4}
            \addplot table[skip first n=1,x expr={\thisrowno{0}}, y expr={\thisrowno{11}}] {./figures/euler/euler_convergence_1e-2.csv};
            \addplot table[skip first n=1,x expr={\thisrowno{0}}, y expr={\thisrowno{9}}] {./figures/euler/euler_convergence_1e-2.csv};
            \addplot[color=gray] coordinates {(0.000375,69) (0.025,4600)}; 
            \legend{predictor MB,predictor MF}
            \end{loglogaxis}
            \end{tikzpicture}
        }
        {
            \includegraphics{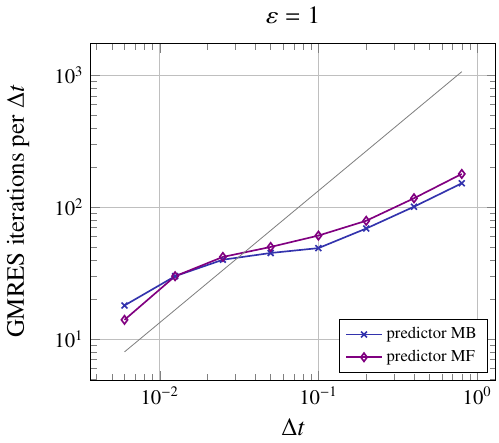}&
            \includegraphics{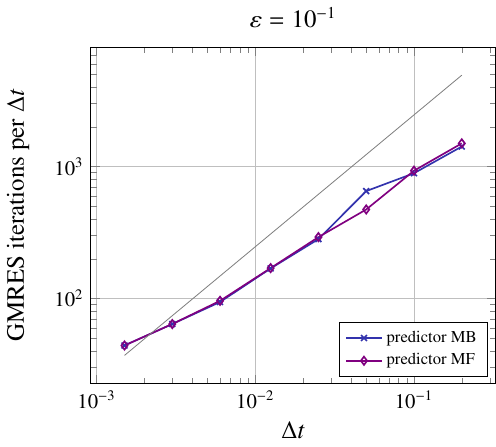}&
            \includegraphics{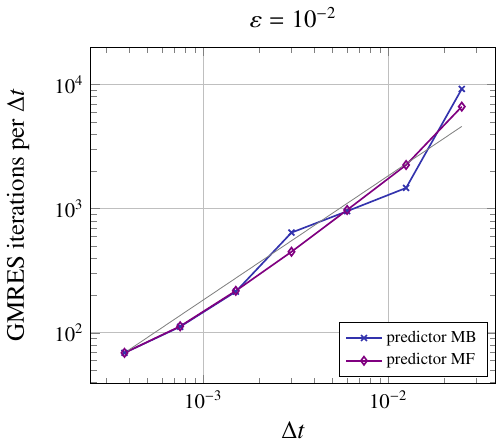}
        }
	\end{tabular}
	\caption{Influence on required GMRES iterations of choosing a matrix-based (MB) or a matrix-free (MF) discretization. To improve convergence, right $\text{BJ}_{\text{ext}}$ preconditioning is applied.}\label{fig:Euler_MF}
\end{figure}
In Fig.~\ref{fig:Euler_MF} the matrix-free and the matrix-based approach are compared regarding the required linear iterations. One can see that the required iterations are very similar. This illustrates that for this setup, a matrix-free approach can be chosen which allows for an efficient spatial parallelization. Consequently, the results in the remainder of this paper are obtained with the matrix-free approach choosing the $\text{BJ}_{\text{ext}}$ preconditioner.

\subsection{Experimental Order of Convergence}
\subsubsection{Linear Scalar Advection}
In order to validate the implementation and to illustrate the high order accuracy of the $\method{q}{\km}$ schemes, we perform a set of simulations with fixed spatial discretization where we vary the timestep size and the time discretization method. The spatial domain is discretized with $n_E=32^2$ elements with ${N}=7$. The final time of the simulation is set to $T_{\text{end}}=0.8$ and the convergence criteria for the GMRES method and Newton's method are set to $\eps_{\text{GMRES}}=10^{-5}$ and $\eps_{\text{Newton}}=10^{-12}$.

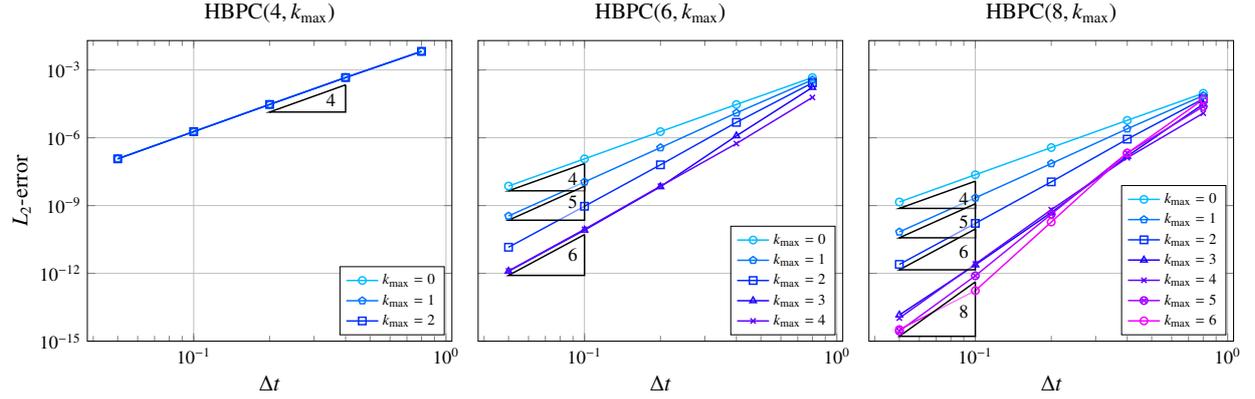
\begin{figure}[ht]
	\centering
	\setlength{\tabcolsep}{0.0em}
	\begin{tabular}{ccc}
        \ifthenelse{\boolean{compilefromscratch}}
        {
            \tikzsetnextfilename{l2_scaladv_4}
            \begin{tikzpicture}[scale=0.7]
            \begin{loglogaxis}[cycle list name=hierarchy,xlabel={$\dt$},ylabel={$L_2$-error} ,grid=major,legend style={at={(0.98,0.02)},anchor=south east,font=\footnotesize,fill opacity=1,text opacity=1},title={$\method{4}{\km}$},label style={font=\large},title style={font=\large},legend cell align={left},ymax=2e-2,ymin=1e-15]
            \addplot table[skip first n=0,x expr={\thisrowno{0}}, y expr={\thisrowno{1}}] {./figures/linadv/convtest_linadv_HBRK4.csv};
            \addplot table[skip first n=0,x expr={\thisrowno{0}}, y expr={\thisrowno{2}}] {./figures/linadv/convtest_linadv_HBRK4.csv};
            \addplot table[skip first n=0,x expr={\thisrowno{0}}, y expr={\thisrowno{3}}] {./figures/linadv/convtest_linadv_HBRK4.csv};
            \addplot[color=black,thick,fill = white,forget plot,fill opacity=0.5] coordinates {(0.4,2.2e-4) (0.4,0.5^4*2.2e-4) (0.2,0.5^4*2.2e-4) (0.4,2.2e-4)}; \node at (axis cs: 0.35, 4.5e-5){$4$};
            \legend{$\km=0$,$\km=1$,$\km=2$}
            \end{loglogaxis}
            \end{tikzpicture}&
            \tikzsetnextfilename{l2_scaladv_6}
            \begin{tikzpicture}[scale=0.7]
            \begin{loglogaxis}[cycle list name=hierarchy,xlabel={$\dt$},yticklabels={,,} ,grid=major,legend style={at={(0.98,0.02)},anchor=south east,font=\footnotesize,fill opacity=1,text opacity=1},title={$\method{6}{\km}$},label style={font=\large},title style={font=\large},legend cell align={left},ymax=2e-2,ymin=1e-15]
            \addplot table[skip first n=0,x expr={\thisrowno{0}}, y expr={\thisrowno{4}}] {./figures/linadv/convtest_linadv_HBRK4.csv};
            \addplot table[skip first n=0,x expr={\thisrowno{0}}, y expr={\thisrowno{5}}] {./figures/linadv/convtest_linadv_HBRK4.csv};
            \addplot table[skip first n=0,x expr={\thisrowno{0}}, y expr={\thisrowno{6}}] {./figures/linadv/convtest_linadv_HBRK4.csv};
            \addplot table[skip first n=0,x expr={\thisrowno{0}}, y expr={\thisrowno{7}}] {./figures/linadv/convtest_linadv_HBRK4.csv};
            \addplot table[skip first n=0,x expr={\thisrowno{0}}, y expr={\thisrowno{8}}] {./figures/linadv/convtest_linadv_HBRK4.csv};
            \addplot[color=black,thick,fill = white,forget plot,fill opacity=0.5] coordinates {(0.1,7.2e-8) (0.1,0.5^4*7.2e-8) (0.05,0.5^4*7.2e-8) (0.1,7.2e-8)}; \node at (axis cs: 0.09, 1.5e-8){$4$};
            \addplot[color=black,thick,fill = white,forget plot,fill opacity=0.5] coordinates {(0.1,7.2e-9) (0.1,0.5^5*7.2e-9) (0.05,0.5^5*7.2e-9) (0.1,7.2e-9)}; \node at (axis cs: 0.09, 1.5e-9){$5$};
            \addplot[color=black,thick,fill = white,forget plot,fill opacity=0.5] coordinates {(0.1,5.2e-11) (0.1,0.5^6*5.2e-11) (0.05,0.5^6*5.2e-11) (0.1,5.2e-11)}; \node at (axis cs: 0.09, 7e-12){$6$};
            \legend{$\km=0$,$\km=1$,$\km=2$,$\km=3$,$\km=4$}
            \end{loglogaxis}
            \end{tikzpicture}&
            \tikzsetnextfilename{l2_scaladv_8}
            \begin{tikzpicture}[scale=0.7]
            \begin{loglogaxis}[cycle list name=hierarchy,xlabel={$\dt$},yticklabels={,,} ,grid=major,legend style={at={(0.98,0.02)},anchor=south east,font=\footnotesize,fill opacity=1,text opacity=1},title={$\method{8}{\km}$},label style={font=\large},title style={font=\large},legend cell align={left},ymax=2e-2,ymin=1e-15]
            \addplot table[skip first n=0,x expr={\thisrowno{0}}, y expr={\thisrowno{9}}] {./figures/linadv/convtest_linadv_HBRK4.csv};
            \addplot table[skip first n=0,x expr={\thisrowno{0}}, y expr={\thisrowno{10}}] {./figures/linadv/convtest_linadv_HBRK4.csv};
            \addplot table[skip first n=0,x expr={\thisrowno{0}}, y expr={\thisrowno{11}}] {./figures/linadv/convtest_linadv_HBRK4.csv};
            \addplot table[skip first n=0,x expr={\thisrowno{0}}, y expr={\thisrowno{12}}] {./figures/linadv/convtest_linadv_HBRK4.csv};
            \addplot table[skip first n=0,x expr={\thisrowno{0}}, y expr={\thisrowno{13}}] {./figures/linadv/convtest_linadv_HBRK4.csv};
            \addplot table[skip first n=0,x expr={\thisrowno{0}}, y expr={\thisrowno{14}}] {./figures/linadv/convtest_linadv_HBRK4.csv};
            \addplot table[skip first n=0,x expr={\thisrowno{0}}, y expr={\thisrowno{15}}] {./figures/linadv/convtest_linadv_HBRK4.csv};
            \addplot[color=black,thick,fill = white,forget plot,fill opacity=0.5] coordinates {(0.1,1.2e-8) (0.1,0.5^4*1.2e-8) (0.05,0.5^4*1.2e-8) (0.1,1.2e-8)}; \node at (axis cs: 0.09, 2e-9){$4$};
            \addplot[color=black,thick,fill = white,forget plot,fill opacity=0.5] coordinates {(0.1,1.2e-9) (0.1,0.5^5*1.2e-9) (0.05,0.5^5*1.2e-9) (0.1,1.2e-9)}; \node at (axis cs: 0.09, 2e-10){$5$};
            \addplot[color=black,thick,fill = white,forget plot,fill opacity=0.5] coordinates {(0.1,9.2e-11) (0.1,0.5^6*9.2e-11) (0.05,0.5^6*9.2e-11) (0.1,9.2e-11)}; \node at (axis cs: 0.09, 9e-12){$6$};
            \addplot[color=black,thick,fill = white,forget plot,fill opacity=0.5] coordinates {(0.1,4.2e-13) (0.1,0.5^8*4.2e-13) (0.05,0.5^8*4.2e-13) (0.1,4.2e-13)}; \node at (axis cs: 0.09, 2e-14){$8$};
            \legend{$\km=0$,$\km=1$,$\km=2$,$\km=3$,$\km=4$,$\km=5$,$\km=6$}
            \end{loglogaxis}
            \end{tikzpicture}
        }
        {
            \includegraphics{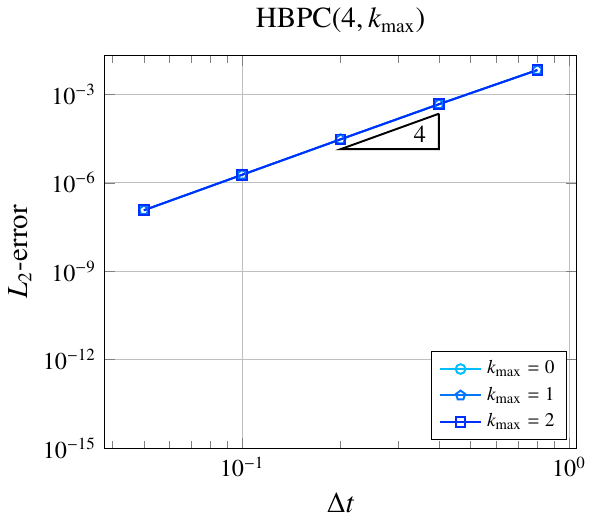}&
            \includegraphics{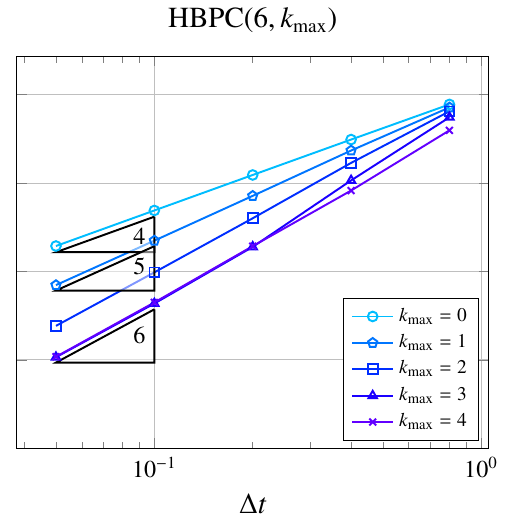}&
            \includegraphics{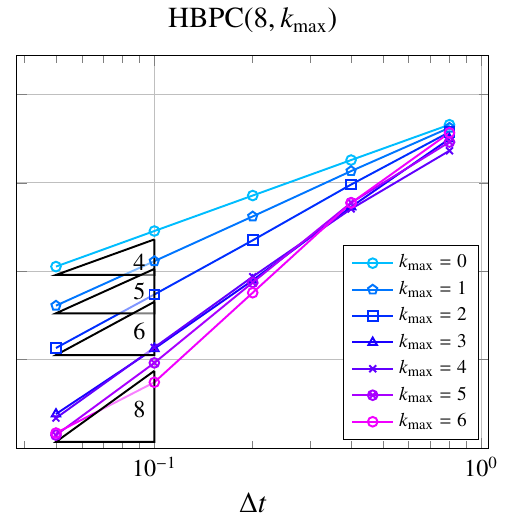}
        }
	\end{tabular}
	\caption{$L_2$-error for linear scalar advection of sine wave with $T_\text{end}=0.8$ for different $\method{q}{\km}$ methods. The spatial domain is discretized with $n_E=32^2$ elements with ${N}=7$. Note that the $L^2$-error of the initial spatial projection is approximately $1.1\cdot10^{-15}$.}\label{fig:LinAdv_convtest_HBRK4}
\end{figure}

The resulting $L_2$-errors of this series of simulations are visualized in Fig.~\ref{fig:LinAdv_convtest_HBRK4}. The figure shows that the temporal error decreases with the desired order of convergence for all considered $\method{q}{\km}$ schemes:
If the $\method{4}{\km}$ scheme is chosen, the corrector steps do not have any influence on the solution. This is an expected behavior as the limiting fourth order Hermite-Birkhoff Runge-Kutta scheme is already used as a predictor. For the $\method{6}{\km}$ and $\method{8}{\km}$ schemes one can see that starting with the order of the predictor, the schemes pick up one order of accuracy per correction step until the maximum order of the underlying quadrature rule is obtained. The only exception from this behavior is the $\method{8}{\km}$ scheme with $\km\ge4$. Here, for some refinement steps a slightly larger order of convergence as expected can be observed. One can see that using more correction steps than required to reach the maximum order of convergence can have a favorable influence on the accuracy of the method.

Note that the presented method is not limited to eighth order of accuracy in time. If one chooses another quadrature rule in Alg.~\ref{alg:mdpde} and a sufficient amount of corrector steps, in principle, arbitrary order can be achieved.

\subsubsection{Euler Equations}
A similar investigation is performed for the non-linear Euler equations with $\eps=1$. Again, the domain is discretized with $n_E=32^2$ elements with a polynomial degree of the ansatz functions of ${N}=7$ and $T_{\text{end}}=0.8$ is chosen as final time. The convergence tolerances are chosen to be $\eps_{\text{GMRES}}=10^{-5}$ and $\eps_{\text{Newton}}=10^{-12}$. 
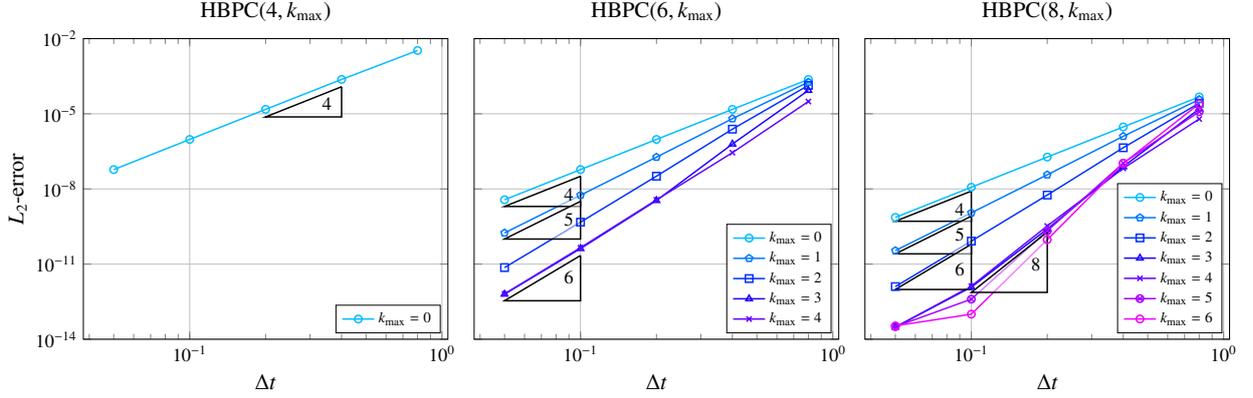
\begin{figure}[ht]
	\setlength{\tabcolsep}{0.0em}
	\centering
	\begin{tabular}{ccc}
        \ifthenelse{\boolean{compilefromscratch}}
        {
            \tikzsetnextfilename{l2_euler_4}
            \begin{tikzpicture}[scale=0.7]
            \begin{loglogaxis}[cycle list name=hierarchy,xlabel={$\dt$},ylabel={$L_2$-error} ,grid=major,legend style={at={(0.98,0.02)},anchor=south east,font=\footnotesize,fill opacity=1,text opacity=1},title={$\method{4}{\km}$},label style={font=\large},title style={font=\large},legend cell align={left},ymax=1e-2,ymin=1e-14]
            \addplot table[skip first n=0,x expr={\thisrowno{0}}, y expr={\thisrowno{1}+\thisrowno{2}+\thisrowno{3}+\thisrowno{4}}] {./figures/euler/convtest_euler_HBPC4_HBRK4.csv};
            \addplot[color=black,thick,fill = white,forget plot,fill opacity=0.5] coordinates {(0.4,1.2e-4) (0.4,0.5^4*1.2e-4) (0.2,0.5^4*1.2e-4) (0.4,1.2e-4)}; \node at (axis cs: 0.35, 2.5e-5){$4$};
            \legend{$\km=0$,$\km=1$,$\km=2$}
            \end{loglogaxis}
            \end{tikzpicture}&
            \tikzsetnextfilename{l2_euler_6}
            \begin{tikzpicture}[scale=0.7]
            \begin{loglogaxis}[cycle list name=hierarchy,xlabel={$\dt$},yticklabels={,,} ,grid=major,legend style={at={(0.98,0.02)},anchor=south east,font=\footnotesize,fill opacity=1,text opacity=1},title={$\method{6}{\km}$},label style={font=\large},title style={font=\large},legend cell align={left},ymax=1e-2,ymin=1e-14]
            \addplot table[skip first n=0,x expr={\thisrowno{0}}, y expr={\thisrowno{1}+\thisrowno{2}+\thisrowno{3}+\thisrowno{4}}] {./figures/euler/convtest_euler_HBPC6_HBRK4.csv};
            \addplot table[skip first n=0,x expr={\thisrowno{0}}, y expr={\thisrowno{5}+\thisrowno{6}+\thisrowno{7}+\thisrowno{8}}] {./figures/euler/convtest_euler_HBPC6_HBRK4.csv};
            \addplot table[skip first n=0,x expr={\thisrowno{0}}, y expr={\thisrowno{9}+\thisrowno{10}+\thisrowno{11}+\thisrowno{12}}] {./figures/euler/convtest_euler_HBPC6_HBRK4.csv};
            \addplot table[skip first n=0,x expr={\thisrowno{0}}, y expr={\thisrowno{13}+\thisrowno{14}+\thisrowno{15}+\thisrowno{16}}] {./figures/euler/convtest_euler_HBPC6_HBRK4.csv};
            \addplot table[skip first n=0,x expr={\thisrowno{0}}, y expr={\thisrowno{17}+\thisrowno{18}+\thisrowno{19}+\thisrowno{20}}] {./figures/euler/convtest_euler_HBPC6_HBRK4.csv};
            \addplot[color=black,thick,fill = white,forget plot,fill opacity=0.5] coordinates {(0.1,3.2e-8) (0.1,0.5^4*3.2e-8) (0.05,0.5^4*3.2e-8) (0.1,3.2e-8)}; \node at (axis cs: 0.09, 6.0e-9){$4$};
            \addplot[color=black,thick,fill = white,forget plot,fill opacity=0.5] coordinates {(0.1,3.2e-9) (0.1,0.5^5*3.2e-9) (0.05,0.5^5*3.2e-9) (0.1,3.2e-9)}; \node at (axis cs: 0.09, 6.0e-10){$5$};
            \addplot[color=black,thick,fill = white,forget plot,fill opacity=0.5] coordinates {(0.1,2.2e-11) (0.1,0.5^6*2.2e-11) (0.05,0.5^6*2.2e-11) (0.1,2.2e-11)}; \node at (axis cs: 0.09, 3e-12){$6$};
            \legend{$\km=0$,$\km=1$,$\km=2$,$\km=3$,$\km=4$}
            \end{loglogaxis}
            \end{tikzpicture}&
            \tikzsetnextfilename{l2_euler_8}
            \begin{tikzpicture}[scale=0.7]
            \begin{loglogaxis}[cycle list name=hierarchy,xlabel={$\dt$},yticklabels={,,} ,grid=major,legend style={at={(0.98,0.02)},anchor=south east,font=\footnotesize,fill opacity=1,text opacity=1},title={$\method{8}{\km}$},label style={font=\large},title style={font=\large},legend cell align={left},ymax=1e-2,ymin=1e-14]
            \addplot table[skip first n=0,x expr={\thisrowno{0}}, y expr={\thisrowno{1}+\thisrowno{2}+\thisrowno{3}+\thisrowno{4}}] {./figures/euler/convtest_euler_HBPC8_HBRK4.csv};
            \addplot table[skip first n=0,x expr={\thisrowno{0}}, y expr={\thisrowno{5}+\thisrowno{6}+\thisrowno{7}+\thisrowno{8}}] {./figures/euler/convtest_euler_HBPC8_HBRK4.csv};
            \addplot table[skip first n=0,x expr={\thisrowno{0}}, y expr={\thisrowno{9}+\thisrowno{10}+\thisrowno{11}+\thisrowno{12}}] {./figures/euler/convtest_euler_HBPC8_HBRK4.csv};
            \addplot table[skip first n=0,x expr={\thisrowno{0}}, y expr={\thisrowno{13}+\thisrowno{14}+\thisrowno{15}+\thisrowno{16}}] {./figures/euler/convtest_euler_HBPC8_HBRK4.csv};
            \addplot table[skip first n=0,x expr={\thisrowno{0}}, y expr={\thisrowno{17}+\thisrowno{18}+\thisrowno{19}+\thisrowno{20}}] {./figures/euler/convtest_euler_HBPC8_HBRK4.csv};
            \addplot table[skip first n=0,x expr={\thisrowno{0}}, y expr={\thisrowno{21}+\thisrowno{22}+\thisrowno{23}+\thisrowno{24}}] {./figures/euler/convtest_euler_HBPC8_HBRK4.csv};
            \addplot table[skip first n=0,x expr={\thisrowno{0}}, y expr={\thisrowno{25}+\thisrowno{26}+\thisrowno{27}+\thisrowno{28}}] {./figures/euler/convtest_euler_HBPC8_HBRK4.csv};
            \addplot[color=black,thick,fill = white,forget plot,fill opacity=0.5] coordinates {(0.1,8.2e-9) (0.1,0.5^4*8.2e-9) (0.05,0.5^4*8.2e-9) (0.1,8.2e-9)}; \node at (axis cs: 0.09, 1.5e-9){$4$};
            \addplot[color=black,thick,fill = white,forget plot,fill opacity=0.5] coordinates {(0.1,8.2e-10) (0.1,0.5^5*8.2e-10) (0.05,0.5^5*8.2e-10) (0.1,8.2e-10)}; \node at (axis cs: 0.09, 1.5e-10){$5$};
            \addplot[color=black,thick,fill = white,forget plot,fill opacity=0.5] coordinates {(0.1,6.2e-11) (0.1,0.5^6*6.2e-11) (0.05,0.5^6*6.2e-11) (0.1,6.2e-11)}; \node at (axis cs: 0.09, 6e-12){$6$};
            \addplot[color=black,thick,fill = white,forget plot,fill opacity=0.5] coordinates {(0.2,1.9e-10) (0.2,0.5^8*1.9e-10) (0.1,0.5^8*1.9e-10) (0.2,1.9e-10)}; \node at (axis cs: 0.18, 1e-11){$8$};
            \legend{$\km=0$,$\km=1$,$\km=2$,$\km=3$,$\km=4$,$\km=5$,$\km=6$}
            \end{loglogaxis}
            \end{tikzpicture}
        }
        {
            \includegraphics{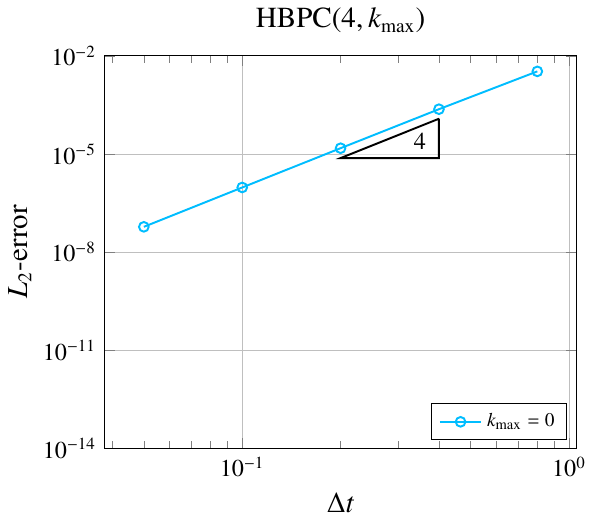}&
            \includegraphics{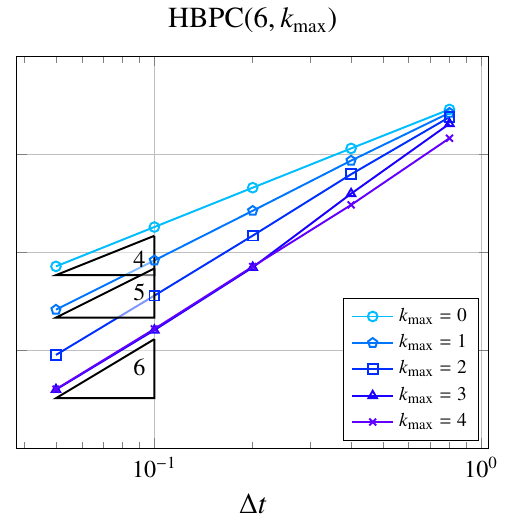}&
            \includegraphics{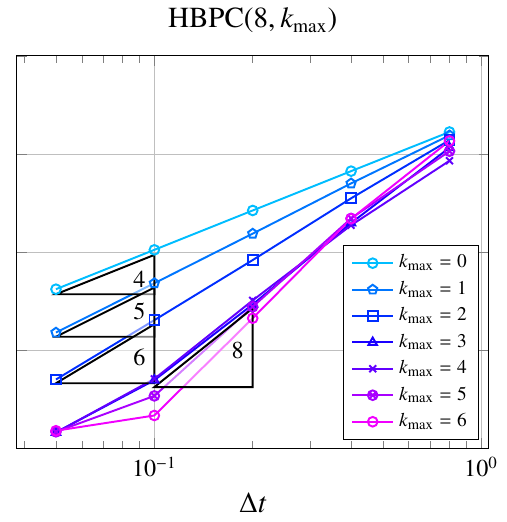}
        }
	\end{tabular}
	\caption{$L_2$-error for Euler equations with advection of sine wave with $T_\text{end}=0.8$ for different $\method{q}{\km}$ methods. The spatial domain is discretized with $n_E=32^2$ elements with ${N}=7$. Note that the $L^2$-error of the initial spatial projection is approximately $1.2\cdot10^{-15}$.}\label{fig:Euler_convtest_HBRK4}
\end{figure}

The total $L_2$-errors of these simulations of the test setup described in Eq.~\eqref{eq:euler_convtest} are visualized in Fig.~\ref{fig:Euler_convtest_HBRK4}. Again, the expected order of accuracy can be observed for all considered schemes in most cases. Similar as for the linear case,  the $\method{8}{\km}$ schemes shows a slightly larger order of convergence as one would expect  for some refinement steps and $\km\ge4$.

	\section{Application to the Navier-Stokes Equations}\label{sec:Applications}
	In this section, we extend the the presented method to handle the Navier-Stokes equations which are given by
\begin{align}\label{eq:NavierStokes}
	\vec{w}_t+\nabla_x\cdot\left(\vec{F}(\vec{w})-\vec{F}^v(\vec{w},\nabla_x\vec{w})\right)=0,\quad\text{with}\quad\vec{w}=\begin{pmatrix}\rho\\\rho\vec{v}\\E\end{pmatrix}\quad\text{and}\quad\vec{F}^v(\vec{w},\nabla_x\vec{w})=\begin{pmatrix}0\\\boldsymbol{\tau}\\\boldsymbol{\tau}\cdot\vec{v}+\vec{q},
	\end{pmatrix}.
\end{align}
They consist of the inviscid Euler flux $\vec{F}(\vec{w})$~\eqref{eq:Euler}, and an additional viscous flux $\vec{F}^v(\vec{w},\nabla_x\vec{w})$. The viscous flux depends on the viscous stress tensor 
\begin{align*}
	\boldsymbol{\tau}:=\mu\left(\nabla_x\vec{v}+(\nabla_x\vec{v})^T-\frac{2}{3}(\nabla_x\cdot\vec{v})\Id\right).
\end{align*}
Here, we have used dynamic viscosity $\mu$,  heat flux $\vec{q}=\lambda\nabla_x T$, thermal conductivity $\lambda=\frac{c_p\mu}{Pr}$, specific heat capacity $c_p=\frac{R\gamma}{\gamma-1}$ and the specific gas constant $R=\frac{1}{\gamma\eps^2}$. The temperature $T$ is defined via the ideal gas law $p=\rho RT$ and the fluid specific Prandtl number is set to $Pr=0.72$. The dynamic viscosity $\mu$ is selected differently for the considered testcases.

\subsection{Fully Discrete Scheme for the Navier-Stokes Equations}
From Eq.~\eqref{eq:NavierStokes} one can see that the viscous flux depends on the state vector $\vec{w}$ and its spatial gradients $\nabla_x\vec{w}$.
More specifically, it depends on the state vector and the gradients of velocity and temperature. 
For the discretization of this second order PDE system, we follow the BR2 lifting approach~\cite{BRMPS}. 
The idea of the lifting procedure is to rewrite the second order Navier-Stokes system as an extended system of first order PDEs, see e.g.~\cite{ABCM} for an overview on this technique.
For that purpose, one additional so-called lifting equation per spatial direction is introduced for the gradient vector $\vec{{d}}=(\vec{d}^1,\dots,\vec{d}^m)$. It contains the spatial gradients of the velocity and the temperature in the different spatial directions
\begin{align}\label{eq:LiftingEq}
	\vec{d}=\nabla_x\vec{w}_{\text{grad}},\quad\text{with}\quad\vec{w}_{\text{grad}}=\left(\vec{v},T\right)^T.
\end{align}

\subsubsection{Calculation of the First Temporal Derivative}
The equations to obtain the first temporal derivative $\vec{w}_t$ are then given by
\begin{align}\label{eq:NS_System}
	\begin{pmatrix}
	\vec{w}_t\\0
	\end{pmatrix}=
	\begin{pmatrix}
	\nabla_x\cdot\left(-\vec{F}(\vec{w})+\vec{F}^v(\vec{w},\vec{d})\right)\\
	\vec{d}-\nabla_x\vec{w}_{\text{grad}}
	\end{pmatrix}.
\end{align}
In order to obtain a discretization of the gradients $\vec{d}$, we follow the BR2 lifting approach~\cite{BRMPS}, where the weak form of Eq.~\eqref{eq:LiftingEq} is used to obtain a discretization for the gradients $\vec{d}$. Summing up, Eq.~\eqref{eq:NS_System} in weak formulation which gives the discrete spatial operator $\Rt_h=\Rt_h(\vec{w},\vec{d}(\vec{w}))$ is given by
\begin{align}
	\sum_{e=1}^{n_E}~\left(\vec{w}_t,\phi\right)_{\Omega_e}-\left(\vec{F}(\vec{w})-\vec{F}^v(\vec{w},\vec{d}),\nabla_x\phi\right)_{\Omega_e}+\left\langle\vec{F}^*(\vec{w}^L,\vec{w}^R)\cdot\vec{n}-{\vec{F}^v}^*(\vec{w}^L,\vec{w}^R,\vec{d}^L,\vec{d}^R)\cdot\vec{n},\phi\right\rangle_{\partial\Omega_e}&=0,\quad\forall\phi\in\Pi_N,\\\label{eq:Lifting}
	\text{with}\qquad\sum_{e=1}^{n_E}~\left(\vec{d},\boldsymbol{\phi}\right)_{\Omega_e}+\left(\vec{w}_{\text{grad}},\nabla_x\boldsymbol{\phi}\right)_{\Omega_e} - \left\langle\vec{w}^*_{\text{grad}}(\vec{w}_{\text{grad}}^L,\vec{w}_{\text{grad}}^R)\cdot\vec{n},\boldsymbol{\phi}\right\rangle_{\partial\Omega_e} &= 0,\quad\forall\boldsymbol{\phi}\in\Pi_N^2.
\end{align}
The viscous flux at the surfaces ${\vec{F}^{v}}^*$ and the lifting surface flux $\vec{w}_{\text{grad}}^*$ are given by the arithmetic means
\begin{align}\label{eq:Riemann_visc}
\begin{split}
	{\vec{F}^{v}}^*(\vec{w}^L,\vec{w}^R,\vec{d}^L,\vec{d}^R)=&\frac{1}{2}\left(\vec{F}^v(\vec{w}^L,\vec{d}^L)+\vec{F}^v(\vec{w}^R,\vec{d}^R)\right)\\
	\vec{w}_{\text{grad}}^*(\vec{w}_{\text{grad}}^L,\vec{w}_{\text{grad}}^R)=&\frac{1}{2}\left(\vec{w}_{\text{grad}}^L+\vec{w}_{\text{grad}}^R\right).
\end{split}
\end{align}
The gradients at the cell edges $\vec{d}^L$, $\vec{d}^R$ are \emph{not} obtained by straightforward evaluation of the solution polynomial of $\vec{d}$.
Instead, a local lifting equation is defined for each edge. The local lifting equations are similar as the global lifting equation given in Eq.~\eqref{eq:Lifting}, but with a modified  surface contribution. For each edge, only the surface integral contribution of the edge itself is taken into account and is scaled with a stabilization parameter. The obtained polynomial representation is then evaluated at the corresponding edge to obtain $\vec{d}^L$ or $\vec{d}^R$.
More details on the BR2 lifting procedure, especially the specific formulation for the DGSEM, can be found in~\cite[Eqs.~(3.81)-(3.86)]{HindenlangDiss}. Note that the implementation of the lifting procedure has been done in the strong form, which is identically to the weak form for the DGSEM~\cite{HindenlangDiss}. It is not indispensable to use the BR2 lifting procedure - other procedures would also be possible, see e.g.~\cite{ABCM} for an overview on different lifting procedures. Here, we choose the BR2 scheme because of its compact stencil as it has only dependencies on direct neighbors, see e.g.~\cite{Ortleb2020}.

\subsubsection{Calculation of the Second Temporal Derivative}
In order to obtain an equation for the discretization of the second temporal derivative, we pursue a similar approach as for the pure hyperbolic equation (see Eq.~\eqref{eq:second_derivative}) and extend the system by the lifting equation
\begin{align}\label{eq:NS_System2}
	\begin{pmatrix}
	\vec{w}_{tt}\\0
	\end{pmatrix}=
	\begin{pmatrix}
	\nabla_x\cdot\left(-\frac{\partial\vec{F}}{\partial\vec{w}}\vec{w}_t\right)+\nabla_x\cdot\left(\frac{\partial\vec{F}^v}{\partial\vec{w}}\vec{w}_t\right)+\nabla_x\cdot\left(\frac{\partial\vec{F}^v}{\partial\vec{d}}\vec{d}_t\right)\\
	\vec{d}_t-\nabla_x\left(\frac{\partial\vec{w}_{\text{grad}}}{\partial\vec{w}}\vec{w}_t\right)
	\end{pmatrix},
\end{align}
where the matrix $\frac{\partial\vec{w}_{\text{grad}}}{\partial\vec{w}}$ is given by the relation of velocity and momentum and by the ideal gas law.
We now make use of the auxiliary variable $\boldsymbol{\sigma}$ (see Eq.~\eqref{eq:Sigma}) and additionally introduce the auxiliary variable $\boldsymbol{\eta}:=\vec{d}_t$. With this, Eq.~\eqref{eq:NS_System2} reads
\begin{align}\label{eq:NS_System3}
	\begin{pmatrix}
	\vec{w}_{tt}\\0
	\end{pmatrix}=
	\begin{pmatrix}
	\nabla_x\cdot\left(-\frac{\partial\vec{F}}{\partial\vec{w}}\boldsymbol{\sigma}\right)+\nabla_x\cdot\left(\frac{\partial\vec{F}^v}{\partial\vec{w}}\boldsymbol{\sigma}\right)+\nabla_x\cdot\left(\frac{\partial\vec{F}^v}{\partial\vec{d}}\boldsymbol{\eta}\right)\\
	\boldsymbol{\eta}-\nabla_x\left(\frac{\partial\vec{w}_{\text{grad}}}{\partial\vec{w}}\boldsymbol{\sigma}\right)
	\end{pmatrix}.
\end{align}
The discretization of Eq.~\eqref{eq:NS_System3} follows the same steps as the discretization for the first temporal derivative (see Eq.~\eqref{eq:second_derivative} and Eq.~\eqref{eq:Euler_wtt} for a comparison). The discrete operator $\Rtt_h=\Rtt_h(\vec{w},\vec{d}(\vec{w}),\boldsymbol{\sigma},\boldsymbol{\eta}(\boldsymbol{\sigma},\vec{w}))$ to obtain the second temporal derivative is then defined by
\begin{align*}
	\sum_{e=1}^{n_E}~\left(\vec{w}_{tt},\phi\right)_{\Omega_e}-\left(\frac{\partial\vec{F(\vec{w})}}{\partial\vec{w}}\cdot\boldsymbol{\sigma},\nabla_x\phi\right)_{\Omega_e}&+
	\left(\frac{\partial\vec{F}^v(\vec{w},\vec{d})}{\partial\vec{w}}\cdot\boldsymbol{\sigma},\nabla_x\phi\right)_{\Omega_e}+
	\left(\frac{\partial\vec{F}^v(\vec{w},\vec{d})}{\partial\vec{d}}\cdot\boldsymbol{\eta},\nabla_x\phi\right)_{\Omega_e}\\
	&+\left\langle\frac{\partial\vec{F}^*(\vec{w}^L,\vec{w}^R)}{\partial\vec{w}^L}\boldsymbol{\sigma}^L\cdot\vec{n}+\frac{\partial\vec{F}^*(\vec{w}^L,\vec{w}^R)}{\partial\vec{w}^R}\boldsymbol{\sigma}^R\cdot\vec{n},\phi\right\rangle_{\partial\Omega_e}\\
	&-\left\langle\frac{1}{2}\left(\frac{\partial{\vec{F}^v}(\vec{w}^L,\vec{d}^L)}{\partial\vec{w}^L}\boldsymbol{\sigma}^L+\frac{\partial{\vec{F}^v}(\vec{w}^R,\vec{d}^R)}{\partial\vec{w}^R}\boldsymbol{\sigma}^R\right)\cdot\vec{n},\phi\right\rangle_{\partial\Omega_e}\\
	&-\left\langle\frac{1}{2}\left(\frac{\partial{\vec{F}^v}(\vec{w}^L,\vec{d}^L)}{\partial\vec{d}^L}\boldsymbol{\eta}^L+\frac{\partial{\vec{F}^v}(\vec{w}^R,\vec{d}^R)}{\partial\vec{d}^R}\boldsymbol{\eta}^R\right)\cdot\vec{n},\phi\right\rangle_{\partial\Omega_e}
	=0,\quad\forall\phi\in\Pi_N,
\end{align*}
with
\begin{align}\label{eq:Lifting_eta}
	\sum_{e=1}^{n_E}~\left(\boldsymbol{\eta},\boldsymbol{\phi}\right)_{\Omega_e} + \left(\left(\frac{\partial\vec{w}_{\text{grad}}}{\partial\vec{w}}\boldsymbol{\sigma}\right),\nabla_x\boldsymbol{\phi}\right)_{\Omega_e} - \left\langle\frac{1}{2}\left(\frac{\partial\vec{w}_{\text{grad}}^L}{\partial\vec{w}^L}\boldsymbol{\sigma}^L+\frac{\partial\vec{w}_{\text{grad}}^R}{\partial\vec{w}^R}\boldsymbol{\sigma}^R\right)\cdot\vec{n},\boldsymbol{\phi}\right\rangle_{\partial\Omega_e} &= 0,\quad\forall\boldsymbol{\phi}\in\Pi_N^2.
\end{align}
Note that we have directly used the simple structure of the viscous numerical flux and the lifting numerical flux (Eq.~\eqref{eq:Riemann_visc}) in the definition of $\Rtt_h$.
The values of $\boldsymbol{\eta}$ at the cell edges ($\boldsymbol{\eta}^{L/R}$) are obtained for each edge by the local forms of the lifting equation~\eqref{eq:Lifting_eta}.

\subsubsection{Solving the Non-Linear System}
Finally, the non-linear equation system which is solved by Newton's method, see Eq.~\eqref{eq:nonlinear_eqsys}, has to be modified accordingly. The matrix-vector product of the system matrix $\mathbfcal{J}$ and the increment vector $\Delta\boldsymbol{X}$ is then given by
\begin{align*}
	\mathbfcal{J}\Delta\boldsymbol{X}=\begin{pmatrix}\Id-\alpha_1\dt\left(\frac{\partial\Rt_h}{\partial\ww}+\frac{\partial\Rt_h}{\partial\vec{d}}\frac{\partial\vec{d}}{\partial\ww}\right)+\frac{\alpha_2\Delta t^2}{2}\left(\frac{\partial\Rtt_h}{\partial\ww}+\frac{\partial\Rtt_h}{\partial\vec{d}}\frac{\partial\vec{d}}{\partial\ww}+\frac{\partial\Rtt_h}{\partial\boldsymbol{\eta}}\frac{\partial\boldsymbol{\eta}}{\partial\ww}\right)&
	\frac{\alpha_2\Delta t^2}{2}\left(\frac{\partial\Rtt_h}{\partial\boldsymbol{\sigma}}+\frac{\partial\Rtt_h}{\partial\boldsymbol{\eta}}\frac{\partial\boldsymbol{\eta}}{\partial\boldsymbol{\sigma}}\right)\\
	-\left(\frac{\partial\Rt_h}{\partial\ww}+\frac{\partial\Rt_h}{\partial\vec{d}}\frac{\partial\vec{d}}{\partial\ww}\right)&
	\Id\end{pmatrix}\begin{pmatrix}\Delta\ww\\\Delta\boldsymbol{\sigma}\end{pmatrix}.
\end{align*}
Again, due to the definitions of $\boldsymbol{\sigma}$ and $\boldsymbol{\eta}$ a tedious but straightforward analysis reveals that 
\begin{align*}
	\frac{\partial\Rt_h}{\partial\vec{d}}\frac{\partial\vec{d}}{\partial\ww}\equiv\frac{\partial\Rtt_h}{\partial\boldsymbol{\eta}}\frac{\partial\boldsymbol{\eta}}{\partial\boldsymbol{\sigma}}.
\end{align*}
Similar as it has been done for the Euler equations, we neglect the Hessian contribution $\frac{\partial\Rtt_h}{\partial\ww}+\frac{\partial\Rtt_h}{\partial\vec{d}}\frac{\partial\vec{d}}{\partial\ww}+\frac{\partial\Rtt_h}{\partial\boldsymbol{\eta}}\frac{\partial\boldsymbol{\eta}}{\partial\ww}$ when constructing the preconditioner. In the following, the matrix-free approach (see Eq.~\eqref{eq:MatrixFree_nonlinear}) with the $BJ_{\text{ext}}$ preconditioner is used for the simulations.

\subsection{Numerical Validation}
For the numerical validation of the presented algorithm, we choose the same test setup as used for the Euler equations given in Eq.~\eqref{eq:euler_convtest}. Additionally, viscosity with $\mu=10^{-3}$ is taken into account.
Again, $n_E=32^2$ elements with ${N}=7$ are used for the spatial discretization. As the exact solution is no longer known, the solution of a simulation with a very small explicit timestep ($4^{th}$ order~\cite{CarpenterRK1994}, $\text{CFL}=0.1$) is taken as reference solution at $\Tend=0.8$.
The convergence tolerances are chosen to be $\eps_{\text{GMRES}}=10^{-3}$ and $\eps_{\text{Newton}}=10^{-10}$. We additionally use an absolute convergence tolerance for Newton's method $\|\Delta\ww\|_2\le10^{-12}$. This became necessary as the initial relative Newton norm for correction steps can sometimes be already very low, such that roundoff errors prevent convergence of the relative tolerance.
\begin{figure}[ht]
	\setlength{\tabcolsep}{0.0em}
	\centering
	\begin{tabular}{ccc}
	\ifthenelse{\boolean{compilefromscratch}}
	{
        \tikzsetnextfilename{l2_ns_4}
		\begin{tikzpicture}[scale=0.7]
		\begin{loglogaxis}[cycle list name=hierarchy,xlabel={$\dt$},ylabel={$L_2$-error} ,grid=major,legend style={at={(0.98,0.02)},anchor=south east,font=\footnotesize,fill opacity=1,text opacity=1},title={$\method{4}{\km}$},label style={font=\large},title style={font=\large},legend cell align={left},ymax=1e-2,ymin=1e-14]
		\addplot table[skip first n=0,x expr={\thisrowno{0}}, y expr={\thisrowno{1}+\thisrowno{2}+\thisrowno{3}+\thisrowno{4}}] {./figures/navierstokes/convtest_ns_HBPC4_HBRK4.csv};
		\addplot[color=black,thick,fill = white,forget plot,fill opacity=0.5] coordinates {(0.05,2e-5) (0.05,0.5^4*2e-5) (0.025,0.5^4*2e-5) (0.05,2e-5)}; \node at (axis cs: 0.045, 4e-6){$4$};
		\legend{$\km=0$,$\km=1$,$\km=2$}
		\end{loglogaxis}
		\end{tikzpicture}&
		\tikzsetnextfilename{l2_ns_6}
		\begin{tikzpicture}[scale=0.7]
		\begin{loglogaxis}[cycle list name=hierarchy,xlabel={$\dt$},yticklabels={,,} ,grid=major,legend style={at={(0.98,0.02)},anchor=south east,font=\footnotesize,fill opacity=1,text opacity=1},title={$\method{6}{\km}$},label style={font=\large},title style={font=\large},legend cell align={left},ymax=1e-2,ymin=1e-14]
		\addplot[color=black,thick,fill = white,forget plot,fill opacity=0.5] coordinates {(0.01,2.2e-9) (0.01,0.5^4*2.2e-9) (0.005,0.5^4*2.2e-9) (0.01,2.2e-9)}; \node at (axis cs: 0.009, 6.0e-10){$4$};
		\addplot[color=black,thick,fill = white,forget plot,fill opacity=0.5] coordinates {(0.01,2.9e-10) (0.01,0.5^5*2.9e-10) (0.005,0.5^5*2.9e-10) (0.01,2.9e-10)}; \node at (axis cs: 0.009,7.5e-11){$5$};
		\addplot[color=black,thick,fill = white,forget plot,fill opacity=0.5] coordinates {(0.01,4.2e-11) (0.01,0.5^6*4.2e-11) (0.005,0.5^6*4.2e-11) (0.01,4.2e-11)}; \node at (axis cs: 0.009, 1.6e-12){$6$};
		\addplot table[skip first n=1,x expr={\thisrowno{0}}, y expr={\thisrowno{1}+\thisrowno{2}+\thisrowno{3}+\thisrowno{4}}] {./figures/navierstokes/convtest_ns_HBPC6_HBRK4.csv};
		\addplot table[skip first n=1,x expr={\thisrowno{0}}, y expr={\thisrowno{5}+\thisrowno{6}+\thisrowno{7}+\thisrowno{8}}] {./figures/navierstokes/convtest_ns_HBPC6_HBRK4.csv};
		\addplot table[skip first n=1,x expr={\thisrowno{0}}, y expr={\thisrowno{9}+\thisrowno{10}+\thisrowno{11}+\thisrowno{12}}] {./figures/navierstokes/convtest_ns_HBPC6_HBRK4.csv};
		\addplot table[skip first n=1,x expr={\thisrowno{0}}, y expr={\thisrowno{13}+\thisrowno{14}+\thisrowno{15}+\thisrowno{16}}] {./figures/navierstokes/convtest_ns_HBPC6_HBRK4.csv};
		\addplot table[skip first n=1,x expr={\thisrowno{0}}, y expr={\thisrowno{17}+\thisrowno{18}+\thisrowno{19}+\thisrowno{20}}] {./figures/navierstokes/convtest_ns_HBPC6_HBRK4.csv};
		\legend{$\km=0$,$\km=1$,$\km=2$,$\km=3$,$\km=4$}
		\end{loglogaxis}
		\end{tikzpicture}
&       \tikzsetnextfilename{l2_ns_8}
		\begin{tikzpicture}[scale=0.7]
		\begin{loglogaxis}[cycle list name=hierarchy,xlabel={$\dt$},yticklabels={,,} ,grid=major,legend style={at={(0.98,0.02)},anchor=south east,font=\footnotesize,fill opacity=1,text opacity=1},title={$\method{8}{\km}$},label style={font=\large},title style={font=\large},legend cell align={left},ymax=1e-2,ymin=1e-14]
		\addplot[color=black,thick,fill = white,forget plot,fill opacity=0.5] coordinates {(0.01,4.2e-10) (0.01,0.5^4*4.2e-10) (0.005,0.5^4*4.2e-10) (0.01,4.2e-10)}; \node at (axis cs: 0.009, 1.2e-10){$4$};
		\addplot[color=black,thick,fill = white,forget plot,fill opacity=0.5] coordinates {(0.01,4.9e-11) (0.01,0.5^5*4.9e-11) (0.005,0.5^5*4.9e-11) (0.01,4.9e-11)}; \node at (axis cs: 0.009,1.3e-11){$5$};
		\addplot[color=black,thick,fill = white,forget plot,fill opacity=0.5] coordinates {(0.01,8.2e-12) (0.01,0.5^6*8.2e-12) (0.005,0.5^6*8.2e-12) (0.01,8.2e-12)}; \node at (axis cs: 0.009, 1.8e-12){$6$};
		\addplot[color=black,thick,fill = white,forget plot,fill opacity=0.5] coordinates {(0.025,6.9e-10) (0.025,0.5^8*6.9e-10) (0.0125,0.5^8*6.9e-10) (0.025,6.9e-10)}; \node at (axis cs: 0.022, 3e-11){$8$};
		\addplot table[skip first n=1,x expr={\thisrowno{0}}, y expr={\thisrowno{1}+\thisrowno{2}+\thisrowno{3}+\thisrowno{4}}] {./figures/navierstokes/convtest_ns_HBPC8_HBRK4.csv};
		\addplot table[skip first n=1,x expr={\thisrowno{0}}, y expr={\thisrowno{5}+\thisrowno{6}+\thisrowno{7}+\thisrowno{8}}] {./figures/navierstokes/convtest_ns_HBPC8_HBRK4.csv};
		\addplot table[skip first n=1,x expr={\thisrowno{0}}, y expr={\thisrowno{9}+\thisrowno{10}+\thisrowno{11}+\thisrowno{12}}] {./figures/navierstokes/convtest_ns_HBPC8_HBRK4.csv};
		\addplot table[skip first n=1,x expr={\thisrowno{0}}, y expr={\thisrowno{13}+\thisrowno{14}+\thisrowno{15}+\thisrowno{16}}] {./figures/navierstokes/convtest_ns_HBPC8_HBRK4.csv};
		\addplot table[skip first n=1,x expr={\thisrowno{0}}, y expr={\thisrowno{17}+\thisrowno{18}+\thisrowno{19}+\thisrowno{20}}] {./figures/navierstokes/convtest_ns_HBPC8_HBRK4.csv};
		\addplot table[skip first n=1,x expr={\thisrowno{0}}, y expr={\thisrowno{21}+\thisrowno{22}+\thisrowno{23}+\thisrowno{24}}] {./figures/navierstokes/convtest_ns_HBPC8_HBRK4.csv};
		\addplot table[skip first n=1,x expr={\thisrowno{0}}, y expr={\thisrowno{25}+\thisrowno{26}+\thisrowno{27}+\thisrowno{28}}] {./figures/navierstokes/convtest_ns_HBPC8_HBRK4.csv};
		\legend{$\km=0$,$\km=1$,$\km=2$,$\km=3$,$\km=4$,$\km=5$,$\km=6$}
		\end{loglogaxis}
		\end{tikzpicture}
    }
    {
        \includegraphics{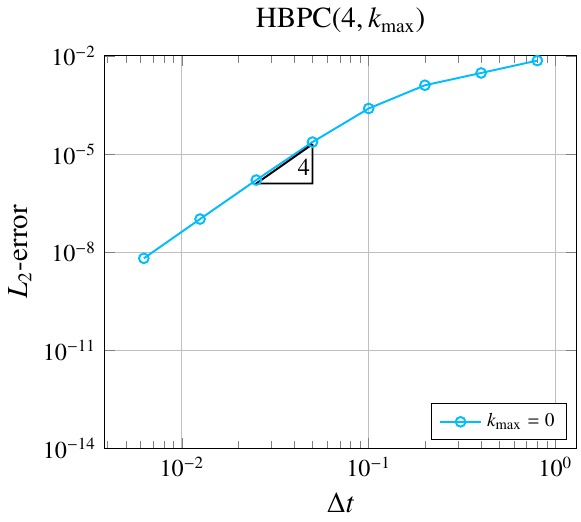}&
        \includegraphics{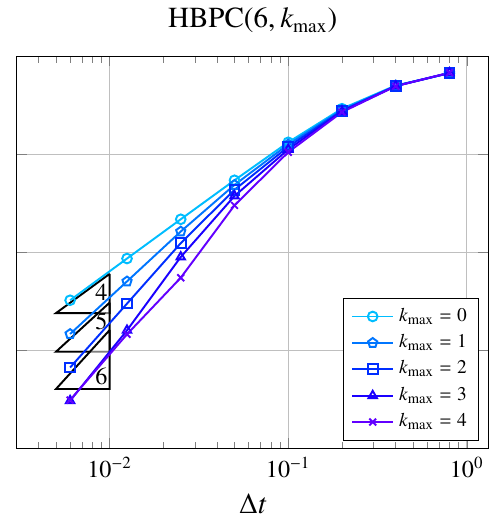}&
        \includegraphics{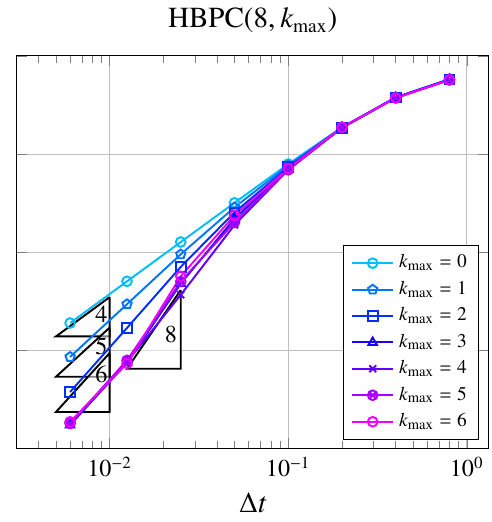}
    }
	\end{tabular}
	\caption{$L_2$-error for Navier-Stokes equations with advection and diffusion of sine wave with $T_\text{end}=0.8$ for different $\method{q}{\km}$ methods. The spatial domain is discretized with $n_E=32^2$ elements with ${N}=7$.}\label{fig:NavierStokes_convtest}
\end{figure}
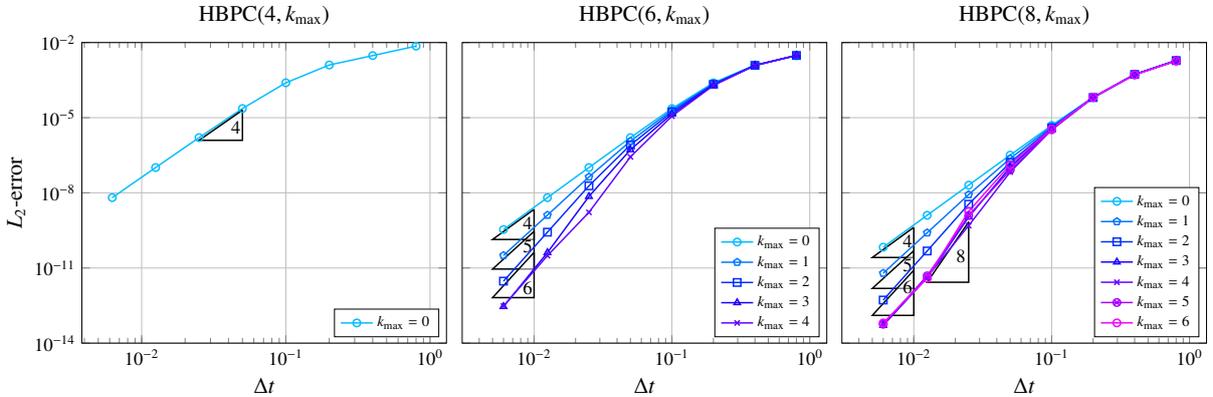

The resulting $L_2$-errors of the simulations are visualized in Fig.~\ref{fig:NavierStokes_convtest}. One can observe that the desired orders are reached for all considered $\method{q}{\km}$ schemes. For the $\method{6}{\km}$ schemes, performing more correction steps than required to reach the maximum order improves the achieved accuracy. For the $\method{8}{\km}$ schemes, the desired order cannot be observed during the last refinement step and $\km\ge3$ due to occurring roundoff errors.

A comparison with the Euler case without viscosity (Fig.~\ref{fig:Euler_convtest_HBRK4}) shows that one needs more temporal refinement steps to reach the asymptotic regime. This is most likely due to the very fast parabolic characteristics. When considering the setup with viscosity, the explicit timestep has to be reduced by approximately a factor of $5$ compared to the inviscid setting. This indicates the dominance of the fast parabolic characteristics and hence why the asymptotic regime is shifted to smaller timestep sizes.

\subsection{Illustrative Applications}
We finally show some illustrative applications of the novel method to typical flow problems.
\subsubsection{Lid-Driven Cavity}
The first example is the lid-driven cavity flow. The quadratic domain $\Omega=[0,1]^2$, which is discretized with $n_E=16^2$ elements, has wall boundary conditions at the left, right and lower boundary. At the upper boundary, a constant flow field
\begin{align*}
	\rho=1,~\vec{v}=(1,0)^T,~p=\frac{1}{\gamma},
\end{align*}
is prescribed. The domain is initialized with the same density and pressure, but with zero velocity. We select $\mu=2.5\cdot10^{-2}$, resulting in a Reynolds number of $Re=400$, and a reference Mach number of $\eps=0.1$.
The spatial and temporal domain are discretized with a sixth order method choosing $N=5$ and $\method{6}{2}$ with $\dt=0.02$. Tolerances of the implicit method are set to $\eps_{\text{Newton}}=10^{-3}$ and $\eps_{\text{GMRES}}=10^{-2}$.
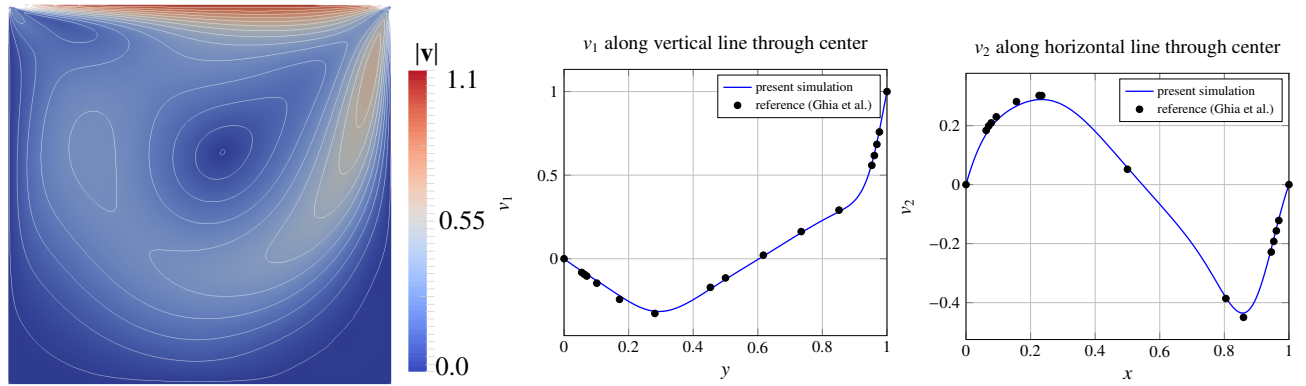
\begin{figure}[ht]
	\centering
	\setlength{\tabcolsep}{0.0em}
	\begin{tabular}{cccc}
		\includegraphics[height=0.17\paperheight]{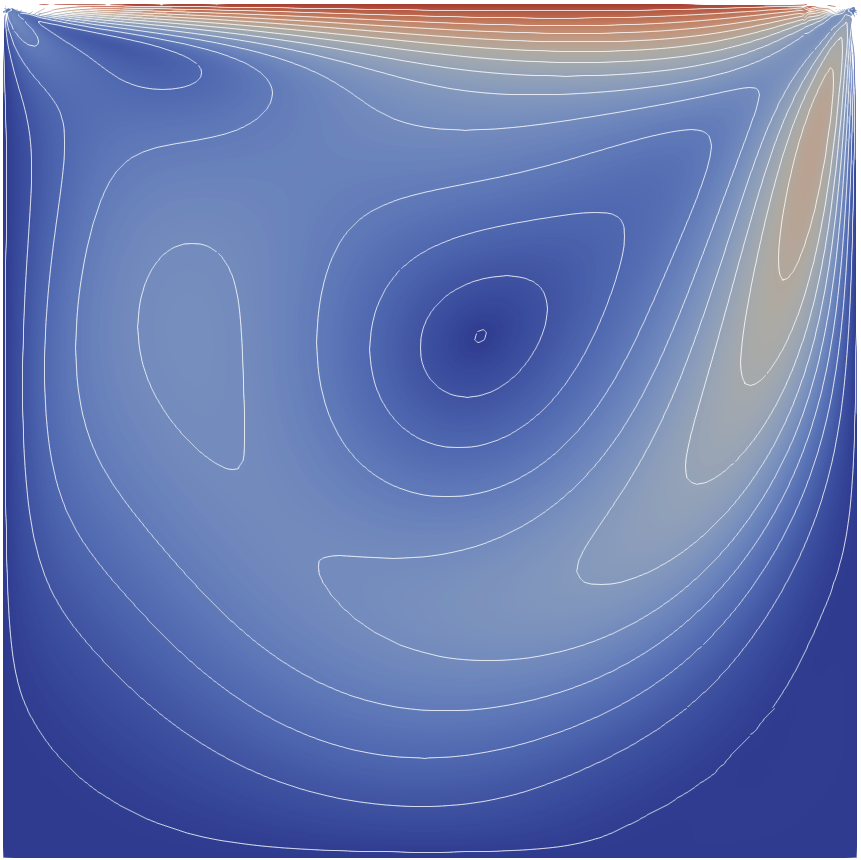}&
		\ifthenelse{\boolean{compilefromscratch}}{
            \tikzsetnextfilename{cavity_1}
            \begin{tikzpicture}[line cap=round,line join=round,x=1.15cm,y=1.15cm,axis/.style={->},scale=0.58]
            \node at (0 ,0) {\includegraphics[scale=0.24]{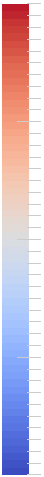}};
            \node at (0.1,3.3 ) {$|\vec{v}|$};
            \node at (0.8,2.8 ) {$1.1$};
            \node at (0.8,0.0 ) {$0.55$};
            \node at (0.6,-2.9) {$0.0$};
            \end{tikzpicture}&
            \tikzsetnextfilename{cavity_vertical}
            \begin{tikzpicture}[scale=0.62]
            \begin{axis}[xlabel={$y$},ylabel={$v_1$},grid=major,legend style={at={(0.98,0.98)},anchor=north east,font=\footnotesize,fill opacity=1,text opacity=1},title={$v_1$ along vertical line through center},label style={font=\large},title style={font=\large},legend cell align={left},xmin=0,xmax=1]
            \addplot[thick,blue] table[skip first n=0,x expr={\thisrowno{7}}, y expr={\thisrowno{1}}] {./figures/navierstokes/cavity_velx_vertical.csv};
            \addplot[only marks] table[skip first n=0,x expr={\thisrowno{1}}, y expr={\thisrowno{2}}] {./figures/navierstokes/Cavity_reference.csv};
            \legend{present simulation, reference (Ghia et al.)}
            \end{axis}
            \end{tikzpicture}&
            \tikzsetnextfilename{cavity_horizontal}
            \begin{tikzpicture}[scale=0.62]
            \begin{axis}[xlabel={$x$},ylabel={$v_2$},grid=major,legend style={at={(0.98,0.98)},anchor=north east,font=\footnotesize,fill opacity=1,text opacity=1},title={$v_2$ along horizontal line through center},label style={font=\large},title style={font=\large},legend cell align={left},xmin=0,xmax=1]
            \addplot[thick,blue] table[skip first n=0,x expr={\thisrowno{2}}, y expr={\thisrowno{7}}] {./figures/navierstokes/cavity_vely_horizontal.csv};
            \addplot[only marks] table[skip first n=0,x expr={\thisrowno{4}}, y expr={\thisrowno{5}}] {./figures/navierstokes/Cavity_reference.csv};
            \legend{present simulation, reference (Ghia et al.)}
            \end{axis}
            \end{tikzpicture}
        }
        {
            \includegraphics{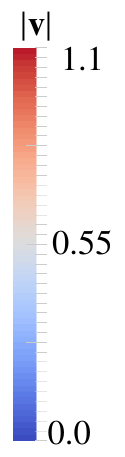}&
            \includegraphics{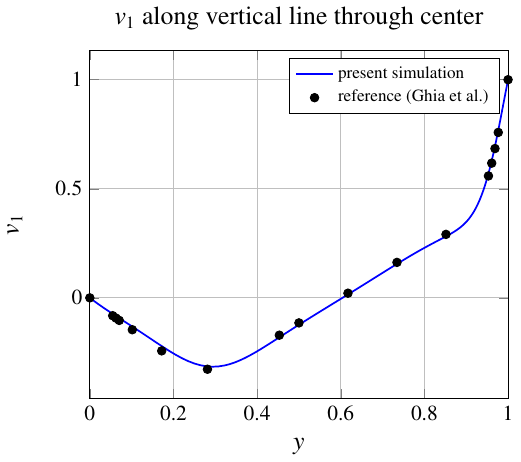}&
            \includegraphics{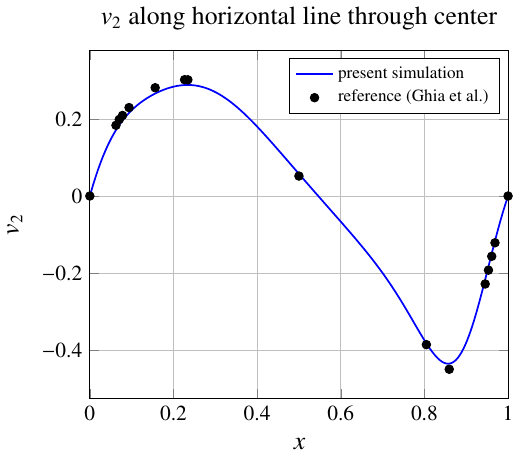}
        }
	\end{tabular}
	\caption{Velocity magnitude (coloring and isocontours) at steady state ($\Tend=50$) for lid-driven cavity flow problem (left). Additional validation is provided by comparison of present simulation with data from Ghia et. al~\cite{Ghia1982}: $v_1$-velocity along vertical line through geometrical center of cavity (middle) and $v_2$-velocity along horizontal line through geometrical center of cavity (right).}\label{fig:Cavity}
\end{figure}
An illustration of the steady state solution at $\Tend=50$ is depicted in Fig.~\ref{fig:Cavity} (left). The results of a comparison of our solution with the simulations from~\cite{Ghia1982}\footnote{Please note that in \cite[Table II]{Ghia1982}, we have left out the point 117, as this number seems to be erroneous.} is shown in Fig.~\ref{fig:Cavity} (middle and right). The present solution matches the reported results in literature very well.
As this is a steady state flow example, the high order accurate time discretization would not be necessary in this case. Nevertheless, this example illustrates that the novel method is capable of simulating such problems.

\subsubsection{Flow Around a Cylinder}
The next example is the two dimensional flow around a cylinder with a Reynolds number $Re_D=200$ and a reference Mach number $\eps=0.1$. The diameter of the cylinder is chosen to be $D=1$. Initially, the flow field is set to a constant state
\begin{align*}
	\rho=1,~\vec{v}=(1,0)^T,~p=\frac{1}{\gamma},
\end{align*}
and the viscosity is set to $\mu=5\cdot10^{-3}$. We use a cylindrical mesh with $n_E=1000$ elements and a sixth order discretization in space and time with $N=5$ and the $\method{6}{2}$ scheme. At the farfield, the initial condition is prescribed as Dirichlet boundary condition, and at the cylinder surface, wall boundaries are applied. A detailed description of the used mesh can be found in~\cite[Sec.~5.1.1]{Vangelatos}. The convergence tolerances for Newton's method and the linear solver are set to $\eps_{\text{Newton}}=10^{-5}$ and $\eps_{\text{GMRES}}=10^{-3}$, respectively.

For this setup, one can expect a vortex shedding behind the cylinder. One measure of the solution quality is the shedding frequency of the wake, which can be analyzed by the frequency of the lift forces at the cylinder. The frequency $f$ is typically related to the lift forces, the cylinder diameter $D$ and the freestream velocity $v_1$, which results in the so-called Strouhal number $\Sr=\frac{fD}{v_1}$.
We compare two different simulations with relatively large timesteps $\dt=0.1$ and $\dt=0.25$. These timestep sizes are approximately $200$ ($\dt=0.1$) and $500$ times ($\dt=0.25$) larger than it would be possible with a $4^{th}$ order explicit low-storage Runge-Kutta scheme~\cite{CarpenterRK1994}.
\begin{figure}[ht]
	\centering
	\setlength{\tabcolsep}{0.0em}
	\begin{tabular}{ccc}
        \includegraphics[height=0.17\paperheight]{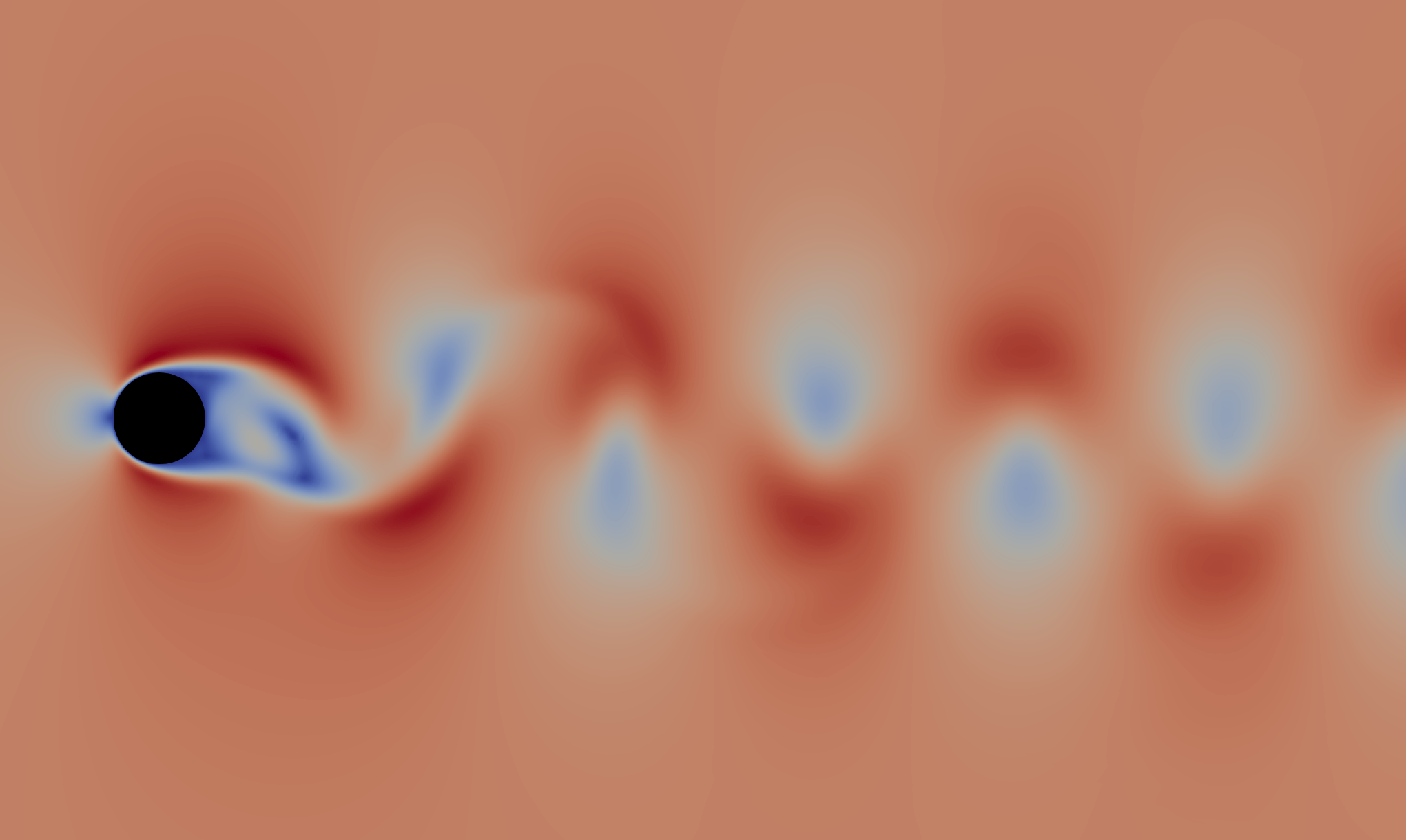}&
        \ifthenelse{\boolean{compilefromscratch}}
        {
            \tikzsetnextfilename{cylinder_1}
            \begin{tikzpicture}[line cap=round,line join=round,x=1.15cm,y=1.15cm,axis/.style={->},scale=0.58]
                \node at (0 ,0) {\includegraphics[scale=0.24]{./figures/navierstokes/bar_color_red.png}};
                \node at (0.1,3.3 ) {$|\vec{v}|$};
                \node at (0.8,2.8 ) {$1.4$};
                \node at (0.8,0.0 ) {$0.7$};
                \node at (0.6,-2.9) {$0.0$};
            \end{tikzpicture}
            \tikzsetnextfilename{cylinder_liftforces}
            \begin{tikzpicture}[scale=0.7]
                \begin{axis}[xlabel={$t$},ylabel={lift force},grid=major,legend style={at={(0.98,0.02)},anchor=south east,font=\footnotesize,fill opacity=1,text opacity=1},title={},label style={font=\large},title style={font=\large},legend cell align={left}]
                \addplot table[skip first n=0,x expr={\thisrowno{0}}, y expr={\thisrowno{2}}] {./figures/navierstokes/Cylinder_Re200_eps01_BodyForces_BC_cylinder.csv};
                \addplot table[skip first n=0,x expr={\thisrowno{0}}, y expr={\thisrowno{2}}] {./figures/navierstokes/Cylinder_Re200_eps01_dt025_BodyForces_BC_cylinder.csv};
                \addplot[thick,black] table[skip first n=0,x expr={\thisrowno{0}}, y expr={\thisrowno{2}}] {./figures/navierstokes/Cylinder_Re200_eps01_explicit_BodyForces_BC_cylinder.csv};
                \legend{$\dt=0.1$,$\dt=0.25$,explicit}
                \end{axis}
            \end{tikzpicture}
        }
        {   
            \includegraphics{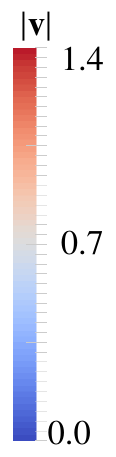}
            \includegraphics{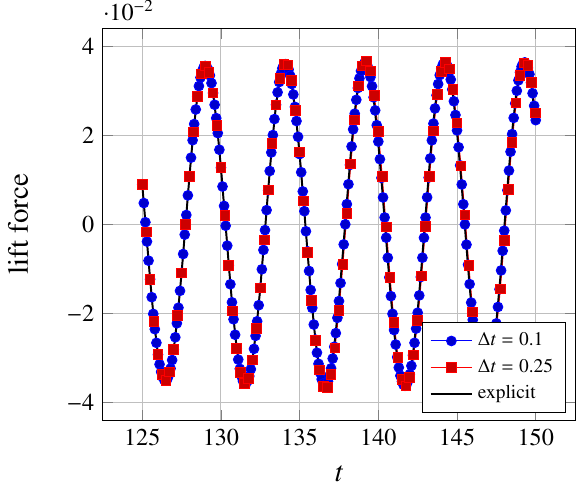}
        }
	\end{tabular}
	\caption{Velocity magnitude at $t=125$ of cylinder flow with $\eps=0.1$ and $Re_D=200$ (left) and temporal evolution of lift forces at cylinder surface (right) choosing different timestep sizes of the $\method{6}{2}$ scheme and explicit reference calculation.}\label{fig:Cylinder}
\end{figure}

Besides the instantaneous flow field at $t=125$, the temporal evolution of the lift forces from $t=125$ to $t=150$ are visualized in Fig.~\ref{fig:Cylinder}. One can see that the calculations with both timestep sizes are very similar and match the explicit reference calculation very well: The Strouhal number is $\Sr=0.1971$ for $\dt=0.1$ and $\Sr=0.1968$ for $\dt=0.25$. This is in very good agreement to what has been reported in literature, see e.g.~\cite{Rajani2009} ($\Sr=0.1957$) and~\cite{Meneghini2001} ($\Sr=0.196$).

\subsubsection{Taylor-Green-Vortex}
Finally, the three dimensional Taylor-Green-Vortex (TGV) is simulated. This illustrates that the method is also applicable to solve three dimensional problems. The TGV is a typical problem to study the transition to turbulence and its decay. It is initialized with large vortices, which then decompose into smaller vortices. When they are small enough they are dissipated by the viscosity of the fluid and their kinetic energy is transferred into internal energy.
A measure of this mechanism is the dissipation rate of the kinetic energy
\begin{align*}
	\frac{\partial E_{\text{kin}}}{\partial t}=\frac{2\mu}{\rho\|\Omega\|}\int_{\Omega}\nabla_x\vec{v}:\nabla_x\vec{v}d\vec{x},
\end{align*}
with the volume of the computational domain $\|\Omega\|$.
Similar as it has been done in~\cite{RSIMEXFullEuler}, we adopt the initialization of the TGV to the non-dimensional equations. The initial data are given by
\begin{align*}
	\rho=1,\quad
	\vec{v}=\begin{pmatrix}\cos(x)\cos(y)\cos(z)\\-\cos(x)\sin(y)\cos(z)\\0\end{pmatrix},\quad\text{and}\quad
	p=\frac{\rho}{\gamma}+\frac{\rho\eps^2}{16}\left(\cos(2x)+\cos(2y)\right)\left(\cos(2z)+2\right).
\end{align*}
We use $n_E=32^3$ elements with $N=3$ to discretize the periodic domain $\Omega=[0,2\pi]^3$. For the temporal discretization the $\method{4}{0}$ scheme with $\dt=0.25$ is chosen. The viscosity is set to $\mu=1.25\cdot10^{-3}$, resulting in a unit Reynolds number of $Re=800$.
\begin{figure}[ht]
	\centering
	\ifthenelse{\boolean{compilefromscratch}}
	{
        \tikzsetnextfilename{kineticenergy}
        \begin{tikzpicture}[scale=0.8]
                \begin{axis}[xlabel={$t$},ylabel={$-\frac{\partial E_{\text{kin}}}{\partial t}$},grid=major,legend style={at={(0.98,0.02)},anchor=south east,font=\footnotesize,fill opacity=1,text opacity=1},title={},label style={font=\large},title style={font=\large},legend cell align={left},xmin=0,xmax=10,ymin=0]
                \addplot table[skip first n=0,x expr={\thisrowno{0}}, y expr={\thisrowno{6}}] {./figures/navierstokes/TGV_Re800_HBRK4_dt025_TGVAnalysis.csv};
                \addplot[thick,red] table[skip first n=0,x expr={\thisrowno{0}}, y expr={\thisrowno{6}}] {./figures/navierstokes/TGV_Re800_N3explicit_TGVAnalysis.csv};
                \addplot[thick] table[skip first n=0,x expr={\thisrowno{0}}, y expr={\thisrowno{1}}] {./figures/navierstokes/TGV_Re800_DNS.csv};
                \legend{$\method{4}{0}$,explicit,DNS}
                \end{axis}
        \end{tikzpicture}
    }
    {
        \includegraphics{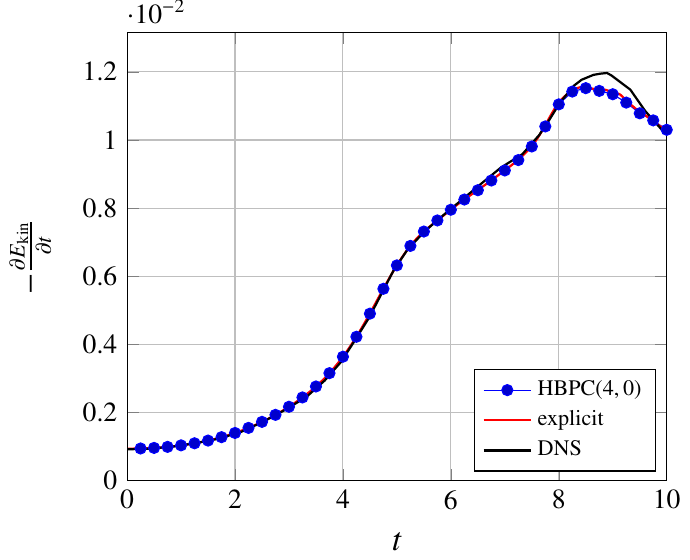}
    }
	\caption{Temporal evolution of kinetic energy dissipation rate for $Re=800$ TGV. The $\method{4}{0}$ scheme uses a timestep size of $\dt=0.25$ and the explicit $4^{th}$ order scheme~\cite{CarpenterRK1994} $\dt\approx1.25\cdot10^{-3}$. Both use a spatial resolution of $n_E=32^3$ elements with $N=3$.
	DNS data are taken from~\cite{Brachet1983}.}\label{fig:TGV}
\end{figure}
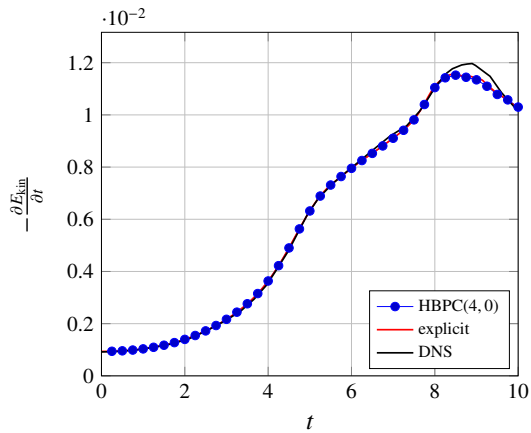
In Fig.~\ref{fig:TGV}, the dissipation rate of the kinetic energy is shown for the novel scheme, an explicit reference simulation and a DNS solution, taken from~\cite{Brachet1983}. One can see, that the DNS is matched very well and only the peak in the decay rate at $t\approx9$ is slightly underestimated. This is most probably caused by a too coarse spatial resolution.
The solution with the $\method{4}{0}$ and the explicit time discretization agree very good although the timestep size of the implicit method is approximately $200$ times larger than for the explicit method.

Summing up, the presented testcases illustrate the applicability of the novel scheme for the simulation of complex flow phenomena.

	\section{Conclusion and Outlook}\label{sec:Conclusion}
	In this work, we have shown how two-derivative deferred correction methods can be combined with the DGSEM to obtain a numerical method that is high order in space and time. It has been illustrated that, due to the implicit nature of the time discretization, very large timesteps can be used. For the discretization of the second temporal derivative, the approach of~\cite{SSJ2017} is adopted to handle non-linear equations.
Different options how to set up the non-linear system to be solved implicitly have been investigated. The preferred method comes with a matrix-free approach, a novel and relatively simple preconditioner and allows for a straight-forward parallelization of the spatial domain.
The flexibility of the novel approach in handling complex applications has been illustrated by the simulation of typical benchmark problems for the Navier-Stokes equations.

The most obvious future development for the novel scheme can be deduced from~\cite{Schutz2021}: the handling of a mixed implicit-explicit (IMEX) flux splitting, such as e.g.~\cite{toro2012flux,RSIMEXFullEuler} and a temporal parallelization. Both will be addressed in a future work.
In addition, improvements are possible in terms of efficiency and flexibility of the new method. The use of adaptive tolerances and timestep sizes can ease the setup of simulations.
We have observed that the choice of the numerical dissipation in Eq.~\eqref{eq:Riemann} is crucial for both stability and efficiency. A detailed study, also regarding asymptotic consistency and other numerical flux functions is worth pursuing.
Additionally, ideas from the ADER community, see e.g.~\cite{Toro2002,Busto2020}, might help to improve the flexibility in using different Riemann solvers.

	\section*{Acknowledgments}
	The authors would like to thank David Seal and Alexander Jaust for the discussions on multiderivative timestepping schemes.
	J.~Zeifang was funded by the Deutsche Forschungsgemeinschaft (DFG, German Research Foundation) - project no. 457811052. We acknowledge the Institute of Aerodynamics and Gas Dynamics at the University of Stuttgart and the VSC (Flemish Supercomputer Center) for providing computing resources.  The VSC is funded by the Research Foundation - Flanders (FWO) and the Flemish Government.

	\section*{Declaration of Competing Interest}
	The authors declare that they have no known competing financial interests or personal relationships that could have appeared to influence the work reported in this paper.
	
	\bibliographystyle{elsarticle-num}
	\bibliography{ListPaper}

\begin{thebibliography}{10}
\expandafter\ifx\csname url\endcsname\relax
  \def\url#1{\texttt{#1}}\fi
\expandafter\ifx\csname urlprefix\endcsname\relax\def\urlprefix{URL }\fi
\expandafter\ifx\csname href\endcsname\relax
  \def\href#1#2{#2} \def\path#1{#1}\fi

\bibitem{HaWa73}
E.~Hairer, G.~Wanner, Multistep-multistage-multiderivative methods for ordinary
  differential equations, Computing (Arch. Elektron. Rechnen) 11~(3) (1973)
  287--303.

\bibitem{LaxWend1960}
P.~Lax, B.~Wendroff, Systems of conservation laws, Communications on Pure and
  Applied Mathematics 13~(2) (1960) 217--237.

\bibitem{KastlungerWanner1972}
K.~Kastlunger, G.~Wanner, On {T}uran type implicit {R}unge-{K}utta methods,
  Computing 9 (1972) 317--325.

\bibitem{Christlieb2015}
A.~J. Christlieb, Y.~G\"u\c{c}l\"u, D.~C. Seal, The {P}icard integral
  formulation of weighted essentially nonoscillatory schemes, SIAM Journal on
  Numerical Analysis 53~(4) (2015) 1833--1856.

\bibitem{JiangShuZhang13}
Y.~Jiang, C.-W. Shu, M.~Zhang, An alternative formulation of finite difference
  weighted {ENO} schemes with {L}ax-{W}endroff time discretization for
  conservation laws, SIAM Journal on Scientific Computing 35~(2) (2013)
  A1137--A1160.

\bibitem{Moe2017}
S.~A. Moe, J.~A. Rossmanith, D.~C. Seal, Positivity-preserving discontinuous
  {G}alerkin methods with {L}ax--{W}endroff time discretizations, Journal of
  Scientific Computing 71~(1) (2017) 44--70.

\bibitem{Qiu2005}
J.~Qiu, M.~Dumbser, C.-W. Shu, The {discontinuous {{Galerkin}}} method with
  {Lax--Wendroff} type time discretizations, Computer Methods in Applied
  Mechanics and Engineering 194~(42-44) (2005) 4528--4543.

\bibitem{ZorioEtAl}
D.~Zor\'{\i}o, A.~Baeza, P.~Mulet, An approximate {Lax--Wendroff}-type
  procedure for high order accurate schemes for hyperbolic conservation laws,
  Journal of Scientific Computing 71 (2017) 246--273.

\bibitem{Titarev2002}
V.~A. Titarev, E.~F. Toro, {ADER}: Arbitrary high order {G}odunov approach,
  Journal of Scientific Computing 17~(1) (2002) 609--618.

\bibitem{Busto2020}
S.~Busto, S.~Chiocchetti, M.~Dumbser, E.~Gaburro, I.~Peshkov, High order {ADER}
  schemes for continuum mechanics, Frontiers in Physics 8 (2020) 32.

\bibitem{Seal13}
D.~Seal, Y.~G\"u\c{c}l\"u, A.~Christlieb, High-order multiderivative time
  integrators for hyperbolic conservation laws, Journal of Scientific Computing
  60 (2014) 101--140.

\bibitem{Ji2018}
X.~Ji, F.~Zhao, W.~Shyy, K.~Xu, A family of high-order gas-kinetic schemes and
  its comparison with {R}iemann solver based high-order methods, Journal of
  Computational Physics 356 (2018) 150--173.

\bibitem{He2020}
Z.~He, F.~Gao, B.~Tian, J.~Li, Implementation of finite difference weighted
  compact nonlinear schemes with the two-stage fourth-order accurate temporal
  discretization, Communications in Computational Physics 27 (2020) 1470--1484.

\bibitem{PanXuLiLi2016}
L.~Pan, K.~Xu, Q.~Li, J.~Li, An efficient and accurate two-stage fourth-order
  gas-kinetic scheme for the {E}uler and {N}avier-{S}tokes equations, Journal
  of Computational Physics 326 (2016) 197 -- 221.

\bibitem{Chouchoulis2021}
J.~Chouchoulis, J.~Sch{\"u}tz, J.~Zeifang, Jacobian-free explicit
  multiderivative {R}unge-{K}utta methods for hyperbolic conservation laws,
  arXiv preprint arXiv:2107.06633 (2021).

\bibitem{TCW14}
A.~Tsai, R.~Chan, S.~Wang, Two-derivative {{Runge-Kutta}} methods for {PDEs}
  using a novel discretization approach, Numerical Algorithms 65 (2014)
  687--703.

\bibitem{MultiDerHDG2015}
A.~Jaust, J.~Sch{\"u}tz, D.~C. Seal, Implicit multistage two-derivative
  discontinuous {G}alerkin schemes for viscous conservation laws, Journal of
  Scientific Computing 69 (2016) 866--891.

\bibitem{SSJ2017}
J.~Sch\"utz, D.~Seal, A.~Jaust, Implicit multiderivative collocation solvers
  for linear partial differential equations with discontinuous {Galerkin}
  spatial discretizations, Journal of Scientific Computing 73 (2017)
  1145--1163.

\bibitem{DissAlex}
A.~Jaust, Novel implicit unconditionally stable time-stepping for {DG}-type
  methods and related topics, Ph.D. thesis, Hasselt University (2018).

\bibitem{SealSchuetz19}
J.~Sch\"{u}tz, D.~Seal, An asymptotic preserving semi-implicit multiderivative
  solver, Applied Numerical Mathematics 160 (2021) 84--101.

\bibitem{Schutz2021}
J.~Sch{\"u}tz, D.~C. Seal, J.~Zeifang, Parallel-in-time high-order
  multiderivative {IMEX} solvers, arXiv preprint arXiv:2101.07846 (2021).

\bibitem{kopriva2009}
D.~A. Kopriva, Implementing spectral methods for partial differential
  equations: {Algorithms} for scientists and engineers, Springer Science \&
  Business Media, 2009.

\bibitem{ReHi73}
W.~Reed, T.~Hill, Triangular mesh methods for the neutron transport equation,
  Tech. rep., Los Alamos Scientific Laboratory (1973).

\bibitem{ShuReviewDG}
C.-W. Shu, A brief survey on discontinuous {Galerkin} methods in computational
  fluid dynamics, Advances in Mechanics 43 (2013) 541--554.

\bibitem{MunzDG12}
F.~Hindenlang, G.~Gassner, C.~Altmann, A.~Beck, M.~Staudenmaier, C.-D. Munz,
  Explicit discontinuous {Galerkin} methods for unsteady problems, Computers \&
  Fluids 61 (2012) 86--93.

\bibitem{Krais2020}
N.~Krais, A.~Beck, T.~Bolemann, H.~Frank, D.~Flad, G.~Gassner, F.~Hindenlang,
  M.~Hoffmann, T.~Kuhn, M.~Sonntag, et~al., {FLEXI}: A high order discontinuous
  {G}alerkin framework for hyperbolic--parabolic conservation laws, Computers
  \& Mathematics with Applications 81 (2021) 186--219.

\bibitem{Higham2002}
N.~J. Higham, Accuracy and stability of numerical algorithms, SIAM, 2002.

\bibitem{Vangelatos}
S.~Vangelatos, On the efficiency of implicit discontinuous {Galerkin} spectral
  element methods for the unsteady compressible {Navier-Stokes} equations,
  Ph.D. thesis, University of Stuttgart (2019).

\bibitem{Zeifang2020diss}
J.~Zeifang, A discontinuous Galerkin method for droplet dynamics in weakly
  compressible flows, Verlag Dr. Hut, 2020.

\bibitem{petsc3}
S.~Balay, J.~Brown, K.~Buschelman, V.~Eijkhout, W.~D. Gropp, D.~Kaushik, M.~G.
  Knepley, L.~C. McInnes, B.~F. Smith, H.~Zhang, {PETS}c users manual, Tech.
  Rep. ANL-95/11 - Revision 3.1, Argonne National Laboratory (2010).

\bibitem{KS17}
K.~Kaiser, J.~Sch\"utz, A high-order method for weakly compressible flows,
  Communications in Computational Physics 22~(4) (2017) 1150--1174.

\bibitem{RSIMEXFullEuler}
J.~Zeifang, J.~Sch\"utz, K.~Kaiser, A.~Beck,
  M.~Luk{\'a}{\v{c}}ov{\'a}-Medvid'ov{\'a}, S.~Noelle, A novel full-{Euler} low
  {Mach} number {IMEX} splitting, Communications in Computational Physics 27
  (2020) 292--320.

\bibitem{ZKBSM17}
J.~Zeifang, K.~Kaiser, A.~Beck, J.~Sch\"utz, C.-D. Munz, Efficient high-order
  discontinuous {Galerkin} computations of low {Mach} number flows,
  Communications in Applied Mathematics and Computational Science 13 (2018)
  243--270.

\bibitem{Franciolini2017}
M.~Franciolini, A.~Crivellini, A.~Nigro, On the efficiency of a matrix-free
  linearly implicit time integration strategy for high-order discontinuous
  {G}alerkin solutions of incompressible turbulent flows, Computers \& Fluids
  159 (2017) 276--294.

\bibitem{Knoll2004}
D.~A. Knoll, D.~E. Keyes, Jacobian-free {Newton--Krylov} methods: a survey of
  approaches and applications, Journal of Computational Physics 193 (2004)
  357--397.

\bibitem{BRMPS}
F.~Bassi, S.~Rebay, G.~Mariotti, S.~Pedinotti, M.~Savini, A high-order accurate
  discontinuous {Finite Element} method for inviscid and viscous turbomachinery
  flows, Proceedings of 2nd European Conference on Turbomachinery, Fluid
  Dynamics and Thermodynamics (1997) 99--108.

\bibitem{ABCM}
D.~N. Arnold, F.~Brezzi, B.~Cockburn, L.~D. Marini, Unified analysis of
  {discontinuous {{Galerkin}}} methods for elliptic problems, SIAM Journal on
  Numerical Analysis 39 (2002) 1749--1779.

\bibitem{HindenlangDiss}
F.~Hindenlang, Mesh curving techniques for high order parallel simulations on
  unstructured meshes, Ph.D. thesis, University of Stuttgart (2014).

\bibitem{Ortleb2020}
S.~Ortleb, A comparative {F}ourier analysis of discontinuous {G}alerkin schemes
  for advection--diffusion with respect to {BR1}, {BR2}, and local
  discontinuous {G}alerkin diffusion discretization, Mathematical Methods in
  the Applied Sciences 43~(13) (2020) 7841--7863.

\bibitem{CarpenterRK1994}
M.~Carpenter, C.~Kennedy, Fourth-order {$2N$}-storage {Runge-Kutta} schemes,
  Tech. rep., NASA Langley Research Center (1994).

\bibitem{Ghia1982}
U.~Ghia, K.~N. Ghia, C.~Shin, High-{R}e solutions for incompressible flow using
  the {N}avier-{S}tokes equations and a multigrid method, Journal of
  Computational Physics 48~(3) (1982) 387--411.

\bibitem{Rajani2009}
B.~Rajani, A.~Kandasamy, S.~Majumdar, Numerical simulation of laminar flow past
  a circular cylinder, Applied Mathematical Modelling 33~(3) (2009) 1228--1247.

\bibitem{Meneghini2001}
J.~Meneghini, F.~Saltara, C.~Siqueira, J.~Ferrari~Jr, Numerical simulation of
  flow interference between two circular cylinders in tandem and side-by-side
  arrangements, Journal of Fluids and Structures 15~(2) (2001) 327--350.

\bibitem{Brachet1983}
M.~E. Brachet, D.~I. Meiron, S.~A. Orszag, B.~Nickel, R.~H. Morf, U.~Frisch,
  Small-scale structure of the {T}aylor--{G}reen vortex, Journal of Fluid
  Mechanics 130 (1983) 411--452.

\bibitem{toro2012flux}
E.~F. Toro, M.~E. V{\'a}zquez-Cend{\'o}n, Flux splitting schemes for the
  {E}uler equations, Computers \& Fluids 70 (2012) 1--12.

\bibitem{Toro2002}
E.~Toro, V.~Titarev, Solution of the generalized {R}iemann problem for
  advection--reaction equations, Proceedings of the Royal Society of London.
  Series A: Mathematical, Physical and Engineering Sciences 458~(2018) (2002)
  271--281.

\end{thebibliography}
	
	\appendix
	\section{Butcher Tables of the Limiting Hermite-Birkhoff Runge-Kutta Methods}
We consider the following quadrature rules:
	\begin{itemize} 
		\item A fourth-order method ($q=4$) with two stages ($s=2$, one being fully explicit), which exactly corresponds to the method used in~\cite{SealSchuetz19}: 
		\begin{align}\label{eq:butcher4}
		c = \begin{pmatrix}0 \\ 1\end{pmatrix}, \quad
		B^{(1)} = 
		\begin{pmatrix}
		0 & 0 \\[0.5em]
		\frac 1 2& \frac 1 2
		\end{pmatrix}, \quad
		B^{(2)} = 
		\begin{pmatrix}
		0 & 0   \\[0.5em]
		\frac 1 {12} & \frac{-1}{12}
		\end{pmatrix}.
		\end{align}
		\item A sixth-order method  ($q=6$) with three stages ($s=3$, one being fully explicit), as also used in~\cite{Schutz2021,SSJ2017}:
		\begin{align}\label{eq:butcher6}
		c = \begin{pmatrix}0\\ \frac 1 2\\ 1\end{pmatrix}, \quad
		B^{(1)} = 
		\begin{pmatrix}
		0 & 0 & 0 \\[0.5em]
		\frac{101}{480} & \frac{8}{30} & \frac{55}{2400} \\[0.5em]
		\frac{7}{30} & \frac{16}{30} & \frac{7}{30} \\
		\end{pmatrix}, \quad
		B^{(2)} = 
		\begin{pmatrix}
		0 & 0 & 0  \\[0.5em]
		\frac{65}{4800} & -\frac{25}{600} &  -\frac{25}{8000} \\[0.5em]
		\frac{5} {300} & 0 &  -\frac{5}{300} 
		\end{pmatrix}.
		\end{align}
		\item An eigth-order method ($q=8$) with four stages ($s=4$, one being fully explicit), as also used in~\cite{Schutz2021}:
		\begin{align}\label{eq:butcher8}
		c = \begin{pmatrix}0\\[0.5em] \frac 1 3\\[0.5em] \frac 2 3\\[0.5em] 1\end{pmatrix}, \quad
		B^{(1)} = 
		\begin{pmatrix}
		0 & 0 & 0 & 0\\[0.5em]
		\frac{6893}{54432}& \frac{313}{2016}& \frac{89}{2016}& \frac{397}{54432} \\[0.5em]
		\frac{223}{1701}&    \frac{20}{63}&   \frac{13}{63}&   \frac{20}{1701} \\[0.5em]
		\frac{31}{224}&   \frac{81}{224}&  \frac{81}{224}&    \frac{31}{224} 
		\end{pmatrix},\quad
		B^{(2)} = 
		\begin{pmatrix}
		0 & 0 & 0 & 0  \\[0.5em]
		\frac{1283}{272160}& -\frac{851}{30240}& -\frac{269}{30240}& -\frac{163}{272160} \\[0.5em]
		\frac{43}{8505}&    -\frac{16}{945}&    -\frac{19}{945}&     -\frac{8}{8505} \\[0.5em]
		\frac{19}{3360}&    -\frac{9}{1120}&     \frac{9}{1120}&    -\frac{19}{3360}
		\end{pmatrix}.
		\end{align}
	\end{itemize}

\end{document}